\documentclass[11pt]{article}
\usepackage[utf8]{inputenc}
\usepackage[T1]{fontenc}
\usepackage[english]{babel}
\usepackage{indentfirst}
\usepackage{latexsym}
\usepackage[colorlinks=true, bookmarksnumbered=true, bookmarksopen=true,
bookmarksopenlevel=3, pdfstartview=FitH, linkcolor=cyan, pdfmenubar=true,
pdftoolbar=true, bookmarks=true,citecolor=cyan, urlcolor=magenta,
filecolor=cyan,menucolor=black,plainpages=false,pdfpagelabels]{hyperref}
\usepackage[lmargin=2.5cm,tmargin=2.5cm,bmargin=2.5cm,rmargin=2.5cm]%
{geometry}
\usepackage{amsbsy, amsmath,amsfonts,amssymb,amsthm}
\usepackage{times}
\usepackage{graphicx}
\usepackage{amsmath}
\usepackage{amsfonts}
\usepackage{amssymb}
\usepackage{mathabx}
\usepackage[usenames]{color}%
\def\fin { \vskip 0pt \hfill $\diamond$ \vskip 12pt}
\newtheorem{theorem}{Theorem}[section]
\newtheorem{definition}[theorem]{Definition}
\newtheorem{proposition}[theorem]{Proposition}
\newtheorem{corollary}[theorem]{Corollary}
\newtheorem{lemma}[theorem]{Lemma}
\newtheorem{remark}[theorem]{Remark}

\numberwithin{equation}{section}
\begin{document}

\title{On bilinear estimates and critical uniqueness classes for Navier-Stokes equations}
\author{\\{{Lucas C. F. Ferreira$^{1}$}{\thanks{L.C.F. Ferreira was partially supported by CNPq 308799/2019-4, Brazil. Email: lcff@ime.unicamp.br (corresponding author).}},\hspace{0.3cm}
		\ {Jhean E. P\'{e}rez-L\'{o}pez$^{2}$} {\thanks{J.E. P\'{e}rez-L\'{o}pez was partially supported by CNPq 152280/2016-2, Brazil, and Vicerrectoría de Investigación y Extensión, Project 3704, UIS, Colombia. Email: jhean.perez@uis.edu.co.}},\hspace{0.3cm}
		\ {Julio C. Valencia-Guevara$^{3}$ \vspace{0.2cm}}{\thanks{J. C. Valencia-Guevara. Email: jvalenciag@unsa.edu.pe.}}}\\\\{\small $^{1}$ State University of Campinas (Unicamp), IMECC-Department of Mathematics} \\{\small {Rua S\'{e}rgio Buarque de Holanda, 651, CEP 13083-859, Campinas, SP,
		 Brazil.\vspace{0.3cm}}}\\{\small $^{2}$ Industrial University of Santander (UIS), School of Mathematics}\\{\small {A.A. 678, Bucaramanga, Colombia.\vspace{0.3cm}}}\\{\small  $^{3}$ Universidad Nacional de San Agust\'{i}n de Arequipa (UNSA),}\\{\small Departamento Acad\'{e}mico de Matem\'{a}ticas, Calle Santa Catalina,}\\{\small No. 117, Arequipa, Peru.}}

\date{}
\maketitle

\begin{abstract}
We are concerned with bilinear estimates and uniqueness of mild solutions for
the Navier-Stokes equations in critical spaces. For that, we construct general
settings in which estimates for the bilinear term of the mild formulation hold
true without using auxiliary norms such as Kato time-weighted ones. We first
obtain necessary conditions in abstract critical spaces and then consider
further structures to obtain the estimates in general classes of Besov, Morrey
and Besov-Morrey spaces based on Banach spaces. Examples of applications are
provided in different spaces as well as for other PDEs. In particular, as far
as we know, the bilinear estimate and the uniqueness property obtained in the framework of Besov-weak-Herz
spaces are not available in the existing literature. The proofs are mainly
based on characterizations and estimates on the corresponding predual spaces.

\bigskip{} \noindent\textbf{Keywords:} Navier-Stokes equations; Bilinear
estimates; Uniqueness; Mild solutions; Critical spaces

\medskip{} \noindent\textbf{AMS MSC:} 35Q30; 76D05; 35A02; 76D03; 35C15; 42B35

\end{abstract}

\section{Introduction}

In this work we are concerned with the uniqueness of mild solutions for the
Navier-Stokes equations
\begin{align}
\partial_{t}u-\Delta u+\mathbb{P}(u\cdot\nabla)u  &  =0\hspace{1.2cm}%
x\in\mathbb{R}^{n}\ t>0,\label{eq:N-S-equation1}\\
\nabla\cdot u  &  =0\hspace{1.2cm}x\in\mathbb{R}^{n}\ t\geq
0,\label{eq:N-S-equation2}\\
u(x,0)  &  =u_{0}(x)\hspace{0.5cm}\,x\in\mathbb{R}^{n},
\label{eq:N-S-equation3}%
\end{align}
and the corresponding key estimate in critical spaces for the bilinear term%
\begin{equation}
\mathcal{B}(u,v)(t)=-\int_{0}^{t}\nabla U(t-s)\cdot\mathbb{P}(u\otimes
v)(s)\ ds, \label{eq:bilin-form-N-S1}%
\end{equation}
where $n\geq3$, $u=(u_{1},\cdots,u_{n})$ is a vector field in $\mathbb{R}^{n}%
$, the Leray projector $\mathbb{P}$ in matrix notation is given by
$(\mathbb{P})_{kj}=\delta_{kj}+\mathcal{R}_{k}\mathcal{R}_{j}$, where
$\mathcal{R}_{j}$ ($j=1,\ldots,n$) stands for the $j$-th Riesz transform, and
the family $\{U(t)\}_{t\geq0}$ is the heat semigroup, namely $U(t)f=\Phi
(t,x)\ast f$ with $\Phi$ denoting the heat kernel.

It is well known that if $u(x,t)$ is a solution of
\eqref{eq:N-S-equation1}-\eqref{eq:N-S-equation3}, then $u^{(\lambda
)}(x,t):=\lambda u(\lambda x,\lambda^{2}t)$ is also a solution with initial
data $\lambda u_{0}(\lambda x),$ for all $\lambda>0$. So, we have the scaling
map $u\mapsto u^{(\lambda)}$ and the one for the initial data
\begin{equation}
u_{0}(x)\mapsto\lambda u_{0}(\lambda x). \label{eq:scaling}%
\end{equation}
Related to this scaling, a Banach space $Z\subset\mathcal{S}^{\prime}(\mathbb{R}^{n}),$ where $\mathcal{S}^{\prime}(\mathbb{R}^{n})$ denotes the space of tempered distributions, is called critical for the Navier-Stokes equations whether its norm is invariant under (\ref{eq:scaling}) in the sense that $\Vert u_{0}\Vert_{Z}\approx\Vert\lambda u_{0}(\lambda x)\Vert_{Z}$, for all $u_{0}\in Z$ and $\lambda>0$, where the rescaling operations in $\lambda u_{0}(\lambda x)$ should be meant in the sense of distributions.

Considering certain classes of critical Banach spaces $Z,\,$our intent is to
develop bilinear estimates for $\mathcal{B}(u,v)$ in the natural critical
space $L^{\infty}((0,T);Z)$ for the flow associated to the Cauchy problem
(\ref{eq:N-S-equation1})-(\ref{eq:N-S-equation3}) with $u_{0}\in Z$ and
$\nabla\cdot u_{0}=0.$ Roughly speaking, under suitable conditions on $Z,$ we
obtain the bilinear estimate
\begin{equation}
\left\Vert \mathcal{B}(u,v)\right\Vert _{L^{\infty}((0,T);Z)}\leq K\left\Vert
u\right\Vert _{L^{\infty}((0,T);Z)}\left\Vert v\right\Vert _{L^{\infty
}((0,T);Z)}, \label{aux-bili-1}%
\end{equation}
for all $u,v\in L^{\infty}((0,T);Z)$ and $T\in(0,\infty]$. A feature of
(\ref{aux-bili-1}) is that it provides a control on the norm of $L^{\infty
}((0,T);Z)$ without employing auxiliary norms (e.g., Kato time-weighted
norms). As a consequence of (\ref{aux-bili-1}) and some further basic properties of spaces, we obtain the uniqueness of
solutions in the natural critical class $C(\left[0,T\right);\tilde{Z})$
for $T\in(0,\infty]$, regardless the size of solutions, where $\tilde{Z}$
stands for the maximal closed subspace of $Z$ in which the heat semigroup
$\left\{U(t)\right\}_{t\geq0}$ is strongly continuous (see \cite{Meyer},\cite{Lemarie1}).

In existence results of mild solutions for (\ref{eq:N-S-equation1}%
)-(\ref{eq:N-S-equation3}), it is relatively common the use of auxiliary norms
in addition to the natural norm $\left\Vert \cdot\right\Vert _{L^{\infty
}((0,T);Z)}$. In this direction, we have the so-called Kato approach, which
consists in a fixed point argument in a time-dependent critical Banach space
with an auxiliary norm of the type%

\[
\mathop{\sup}\limits_{t>0}t^{\rho}\left\Vert u\right\Vert _{Y}\text{ with
$\rho>0$ and a Banach space }Y\text{.}%
\]
The auxiliary norm is used to control some integral terms arising
in the estimates for the bilinear operator $\mathcal{B}(\cdot,\cdot)$.
Then, in general, that kind of approach itself prevents getting uniqueness
in the class $C(\left[0,T\right);\tilde{Z}),$ providing this property
in a strict subspace of it. In the sequel, referring the reader to
the review books \cite{Lemarie1,Lemarie3} and without making a complete
list, we mention some works that employ approaches with two (or more)
norms to estimate $\mathcal{B}(\cdot,\cdot)$ and obtain global-in-time
well-posedness of (\ref{eq:N-S-equation1})-(\ref{eq:N-S-equation3})
with divergence-free small data in critical spaces: there are results
in homogeneous Sobolev space $\dot{H}^{1/2}(\mathbb{R}^{3})$ \cite{FuKa},
Lebesgue space $L^{n}(\mathbb{R}^{n})$ \cite{Ka1}, Marcinkiewicz
space $L^{(n,\infty)}\left(\mathbb{R}^{n}\right)$ \cite{Barraza},
homogeneous Besov spaces $\dot{B}_{p,\infty}^{\frac{n}{p}-1}(\mathbb{R}^{n})$
with $p>n$ \cite{Cann1}, homogeneous Fourier-Besov-type spaces $F\dot{B}_{p,q}^{2-\frac{3}{p}}\left(\mathbb{R}^{3}\right)$
and $\mathcal{FN}_{p,\mu,q}^{2-\frac{3-\mu}{p}}\left(\mathbb{R}^{3}\right)$
\cite{de Almeida-Fer,IwabTakada,LeiLin}, Morrey spaces $\mathcal{M}_{p}^{n}(\mathbb{R}^{n})$
\cite{Giga,Ka2,Taylor-1}, homogeneous Besov-Morrey spaces $\dot{B}\mathcal{M}_{p,\infty}^{r,\frac{n}{r}-1}(\mathbb{R}^{n})$
with $r>n/2$ \cite{KoYa,Maz}, homogeneous Besov-weak-Herz spaces
$\dot{B}W\dot{K}_{p,q,r}^{\alpha,s}(\mathbb{R}^{n})$ \cite{LucJhe1},
$BMO^{-1}(\mathbb{R}^{n})$ \cite{KocTat}, among others. We also
quote \cite{Karch} where results were obtained in abstract critical
Banach spaces by means of the Kato approach (see also \cite{Weissler}).

In general, estimates in the form (\ref{aux-bili-1}) are more difficult
to obtain than those with auxiliary norms and involve more subtle
arguments. As far as we know, that kind of bilinear estimate has been
proved in $L^{(n,\infty)}\left(\mathbb{R}^{n}\right)$ \cite{Meyer,Yam},
pseudomeasure space $\mathcal{PM}^{n-1}\left(\mathbb{R}^{n}\right)$
\cite{CannKarch,Le-Jan}, Fourier-Besov spaces $F\dot{B}_{p,\infty}^{2-\frac{3}{p}}\left(\mathbb{R}^{3}\right)$
with $p>3$ \cite{KoniYone}, homogeneous weak-Herz space $W\dot{K}_{n,\infty}^{0}\left(\mathbb{R}^{n}\right)$
\cite{Tsutsui-1}, weak-Morrey spaces $\mathcal{M}_{(p,\infty)}^{n}\left(\mathbb{R}^{n}\right)$
with $2<p\leq n$ \cite{Lemarie2,Lucas}, $\dot{B}_{p,\infty}^{\frac{n}{p}-1}\left(\mathbb{R}^{n}\right)$
with $2\leq p<n$ \cite{Cann2}, and Besov-weak-Morrey $\dot{B}W\mathcal{M}_{p,\infty}^{r,\frac{n}{r}-1}(\mathbb{R}^{n})$
with $2<p\leq r<n$ and $n/2<r$ \cite{LucJhe2}.

First we construct an abstract framework of Banach spaces for \eqref{eq:N-S-equation1}-\eqref{eq:N-S-equation3}
in which the estimate (\ref{aux-bili-1}) holds true with a general
feature (see Lemma \ref{Lem: Estimativa integral espacio dual} and
Theorem \ref{Teo: Teo main. Estimativa integral en el dual.}). In
subsection \ref{Sec-Appl-1}, we present examples of applications
which recover some previous results. In the case of the weak-Herz
space $Z=W\dot{K}_{n,\infty}^{0}\left(\mathbb{R}^{n}\right)$, it
is worth noting that our construction provides a different proof for
estimate (\ref{aux-bili-1}) (see Remark \ref{Rem-weak-Herz}). In
fact, the work \cite{Tsutsui-1} follows the spirit of \cite{Lemarie1,Meyer}
while here we are inspired by arguments in \cite{Yam} (see details
more below). Technically speaking, the proofs of Lemma \ref{Lem: Estimativa integral espacio dual}
and Theorem \ref{Teo: Teo main. Estimativa integral en el dual.}
are relatively straightforward as long as we are putting the conditions
we need to carry out the predual approach. Their aim is to serve as
a basis for unifying some proofs when considering examples of spaces
with more structure. Also, Lemma \ref{Lem: Estimativa integral espacio dual}
depends on a scaling condition connecting two predual spaces $E$
and $E_{0}$, even though in the case of the Navier-Stokes equations
we should consider the dual $Z=E^{^{\prime}}$ as being critical,
that is, with the Navier-Stokes scaling $\sigma_{E^{\prime}}=1$ (see
Remark \ref{Rem-exist-2}). The more general condition on the predual
spaces allows to apply the theory to other nonlinear PDEs.

In Section \ref{sect:bilinear-estimate-besov-spaces} we introduce further
structure in the first framework in order to treat certain general classes of
Besov spaces where we obtain (\ref{aux-bili-1}) and provide examples (see
Theorem \ref{Teo. Estimativa bilinear espacios tipo besov} and subsection
\ref{Sec-Besov-weak-Herz}). Afterwards, working in the context of Banach
spaces $X$ of locally integrable functions similar to those of \cite{Lemarie2}%
, we consider the framework of $X$-Morrey spaces $\mathcal{M}[X]^{l}$ (Morrey
spaces based on $X$) and Besov-$X$-Morrey type spaces (Besov spaces based on
$\mathcal{M}[X]^{l}$). Then, we show estimates and properties in those spaces
and their preduals which ensure the needed conditions in order to apply the
general theory in Section \ref{sect:bilinear-estimate-general-spaces} and,
consequently, obtain (\ref{aux-bili-1}) (see Theorems
\ref{Teo:estim-bilin-cambio} and \ref{Teo-Est-BM-Spaces}).

For concreteness, our applications recover known bilinear estimates
in several spaces, such as $L^{n,\infty}\left(\mathbb{R}^{n}\right)$,
$W\dot{K}_{n,\infty}^{0}\left(\mathbb{R}^{n}\right)$, and $\mathcal{M}_{(p,\infty)}^{n}\left(\mathbb{R}^{n}\right)$,
as well as the Besov-type spaces $\dot{B}_{p,\infty}^{\frac{n}{p}-1}\left(\mathbb{R}^{n}\right)$
and $\dot{B}W\mathcal{M}_{p,\infty}^{l,\frac{n}{l}-1}\left(\mathbb{R}^{n}\right)$.
However, our estimate (\ref{aux-bili-1}) in the context of Besov-weak-Herz
spaces $\dot{B}W\dot{K}_{p,\infty,\infty}^{0,\frac{n}{p}-1}\left(\mathbb{R}^{n}\right)$
seems to be a new contribution to the existing literature (see subsection
\ref{Sec-Besov-weak-Herz}). As a byproduct (see Theorem \ref{Teo: condiciones para Unicidad} and its applications in Section \ref{Sec:unicidad}), adapting an argument
by \cite{Meyer} yields a uniqueness class for mild solutions of (\ref{eq:N-S-equation1})-(\ref{eq:N-S-equation3}),
namely the class
\begin{equation}
	C(\left[0,T\right);\tilde{Z})\text{ with }Z=\dot{B}W\dot{K}_{p,\infty,\infty}^{0,\frac{n}{p}-1}\left(\mathbb{R}^{n}\right)=\dot{B}\left[W\dot{K}_{p,\infty}^{0}\left(\mathbb{R}^{n}\right)\right]_{\infty}^{\frac{n}{p}-1}.\label{aux-uniqueness-1}
\end{equation}

In a certain sense, the present work provides bilinear estimates and uniqueness results in classes of critical spaces presenting functional structures compatible with (\ref{eq:bilin-form-N-S1}).
The main constructions encompass characterizations, basic properties and estimates on the predual spaces,
as well as careful use and choices of spaces with suitable interpolation properties, H\"{o}lder-type
inequality and heat semigroup estimates, among other ingredients. For that, we
are motivated by Yamazaki approach in \cite{Yam} as well as the works
\cite{Lucas,LucJhe2}. In fact, our approach can be seen as an adaptation of
that in \cite{Yam} to other spaces with a more intricate structure, such as
the Besov-weak-Herz space above.

Comparing with previous references, we point out that the approach to
obtaining bilinear estimates (\ref{aux-bili-1}) in abstract spaces developed
in \cite{Lemarie2} (see Proposition 4.1 therein) follows the spirit of
\cite{Meyer} by relying on the boundedness of the Riesz operator
$(-\Delta)^{-1/2}$, the Hardy-Littlewood maximal operator and the boundedness
of the pointwise product, and employing spaces of pointwise multipliers, among
others. Also, in order to obtain a $L^{\infty}$-estimate for $\mathcal{B}%
(u,v)(t),$ the pointwise estimate of the kernel of $e^{\Delta t}\mathbb{P}$ by $C(t^{2}+\left\vert x\right\vert ^{4})^{-1}$ plays a central role in \cite{Lemarie2, Meyer}. It is worth mentioning that some basic conditions on
the base space $X$, such as inclusion in $\mathcal{S}^{\prime}(\mathbb{R}^{n})$, translation
invariance, dilatation control (or scaling), and product estimates
(H\"{o}lder-type estimates), are common to other abstract frameworks employed
previously to analyze the Navier-Stokes equations; for example, adapted spaces
\cite{Cann3},\cite{Meyer}, shift-invariant spaces of distributions and local
measures \cite{Lemarie1}, and adequate Banach spaces \cite{Karch}.

Finally, we intend to show the versatility of the approach presented here by
analyzing a reaction-diffusion system with quadratic reaction term and a
nonlocal advection-diffusion system. In this part, the space $\mathcal{M}%
[X]^{l}$ is chosen to illustrate the theory. We obtain estimates in the spirit
of (\ref{aux-bili-1}) and, consequently, the uniqueness property in $C(\left[
0,T\right)  ;\tilde{Z})$ with $Z=\mathcal{M}[X]^{l}$. The estimates for the
corresponding bilinear terms are obtained for general conditions on $X$ (see
Theorems \ref{Theorem:bilinear-estim-B1} and \ref{Theorem:bilinear-estim-B2})
and afterwards applications are presented in the specific case of $X$ being a
Lorentz space. In order to analyze the nonlocal model, we need to obtain some
properties in $\mathcal{M}[X]^{l}$-spaces for fractional derivatives via the
Riesz potential. For that, we extend and adapt some ideas of
\cite{adams-libro}.

Let us describe the organization of this work. Section \ref{sec: preliminares} is devoted to recall some notations and give some preliminaries about interpolation theory in sequence spaces. In Section \ref{sect:bilinear-estimate-general-spaces}, we provide our first abstract framework and present the applications. In Section
\ref{sect:bilinear-estimate-besov-spaces}, keeping in mind the construction
carried out in Section \ref{sect:bilinear-estimate-general-spaces}, we present
the abstract construction in Besov-type spaces and give the applications. In
Sections \ref{sect:bilinear-in-morrey-spaces} and
\ref{sect:bilinear-estimate-in-besov-morrey}, we obtain the bilinear estimate
for $X$-Morrey spaces $Z=\mathcal{M}[X]^{l}$ and Besov-$X$-Morrey type spaces $Z=\dot{B}\left[
\mathcal{M}[X]^{l}\right]  _{\infty}^{\frac{n}{l}-1}$, respectively. The subject of Section \ref{Sec:unicidad} is a uniqueness result compatible with our abstract frameworks constructed in the previous sections. Finally, in Section \ref{section:related-models} we analyze other nonlinearities and PDEs.

\section{Preliminaries\label{sec: preliminares}}
First, let us set out some notations for the rest of the work. Given a real number $1\leq r\leq\infty,$ we denote by $r^\prime$ the Hölder conjugate of $r$, that is, the number such that $1=1/r+1/r^{\prime}.$ Also, if $E$ is a Banach space we denote its dual by $E^{\prime}$. Since all the spaces used in this work are defined in the whole space $\mathbb{R}^{n}$, we omit the symbol $\mathbb{R}^{n}$ in their notation. For example, $\mathcal{S}^{\prime}\left(\mathbb{R}^{n}\right)$ will be denoted by $\mathcal{S}^{\prime}$, $L^{p,\infty}(\mathbb{R}^n)$ by $L^{p,\infty}$, and so on. Given a Banach space $X$ and two $\lambda$-families of nonnegative functionals $F_{\lambda},G_{\lambda}:X\longrightarrow\mathbb{R}$, where $\lambda$ is a positive parameter, we use the notation $F_{\lambda}\approx G_{\lambda}, \forall\lambda>0$, to mean that there exist constants $a,\,b>0$ such that $aF_{\lambda}(f)\leq G_{\lambda}(f)\leq bF_{\lambda}(f)$ for all $\lambda>0$ and $f\in X$. For two norms $F,G$ in $X$, note that we recover the definition of equivalent norms.
	
	Now, we briefly recall some definitions and properties of sequence spaces. Let $s\in\mathbb{R}$, $1\leq r\leq\infty,$ and let $E$ be a
Banach space. The space $\dot{l}_{r}^{s}(E)$ is the set of all sequences
$a=(a_{k})_{k\in\mathbb{Z}}$ such that $a_{k}\in E$, for all $k,$ and
\[
\left\Vert a\right\Vert _{\dot{l}_{r}^{s}\left(  E\right)  }=\left(
\sum\limits_{k=-\infty}^{\infty}2^{ksr}\left\Vert a_{k}\right\Vert _{E}%
^{r}\right)  ^{1/r}<\infty.
\]

Most of the spaces that we consider throughout this paper present the
structure of sequence spaces $\dot{l}_{r}^{s}(E)$ with some suitable Banach
space $E$. \ The next two lemmas contain real interpolation properties for
that kind of spaces (see, e.g., \cite{BL}).

\begin{lemma}
\label{lem: Intepolacion en sequencias manteniendo fijo el espacio base.}
Suppose that $1\leq r_{0},r_{1},r\leq\infty$ and $s_{0}\neq s_{1}.$ Then, we
have that
\[
\left(  \dot{l}_{r_{0}}^{s_{0}}\left(  E\right)  ,\dot{l}_{r_{1}}^{s_{1}%
}\left(  E\right)  \right)  _{\theta,r}=\dot{l}_{r}^{s}\left(  E\right)  ,
\]
where $0<\theta<1$ and $s=(1-\theta)s_{0}+\theta s_{1}.$ If $s=s_{0}=s_{1}$,
it follows that
\[
\left(  \dot{l}_{r_{0}}^{s}\left(  E\right)  ,\dot{l}_{r_{1}}^{s}\left(
E\right)  \right)  _{\theta,r}=\dot{l}_{r}^{s}\left(  E\right)
\]
provided that
\[
\frac{1}{r}=\frac{1-\theta}{r_{0}}+\frac{\theta}{r_{1}}.
\]

\end{lemma}

\begin{lemma}
\label{Lem: Intepolacion en sequencias cambiando el espacio base} Let $1\leq
r_{0},r_{1}<\infty$ and $s_{0}\neq s_{1}.$ Then,
\[
\left(  \dot{l}_{r_{0}}^{s_{0}}\left(  E_{0}\right)  ,\dot{l}_{r_{1}}^{s_{1}%
}\left(  E_{1}\right)  \right)  _{\theta,r}=\dot{l}_{r}^{s}\left(  \left(
E_{0},E_{1}\right)  _{\theta,r}\right)  ,
\]
where $0<\theta<1,$ $s=(1-\theta)s_{0}+\theta s_{1}$ and $\frac{1}{r}%
=\frac{1-\theta}{r_{0}}+\frac{\theta}{r_{1}}.$
\end{lemma}

In the sequel we recall a basic duality property for sequence spaces (see
\cite{Hernandez-Yang}).

\begin{lemma}
\label{Lem: dualidad en espacios de sequencias} Assume that $s\in\mathbb{R}$
and $1\leq r<\infty$. Then, we have the duality relation
\[
\left(  \dot{l}_{r}^{s}\left(  E\right)  \right)  ^{\prime}=\left(  \dot
{l}_{r^{\prime}}^{-s}\left(  E^{\prime}\right)  \right)  .
\]

\end{lemma}

\section{Bilinear estimate in a general framework}

\label{sect:bilinear-estimate-general-spaces} In the present section, we
provide conditions on suitable Banach spaces which yield the bilinear estimate
(\ref{aux-bili-1}) for the Navier-Stokes equations. We start with a basic
definition that introduces a scaling property for Banach spaces of
distributions. Similar definitions can be found in \cite{Karch}%
, \cite{Lemarie2}.

\begin{definition}
\label{def:G-set} For each $\,f\in\mathcal{S}^{\prime}$ and $\lambda>0,$ let
$f_{\lambda}(\cdot)=f\left(  \lambda\cdot\right)  $. We denote by
$\mathcal{G}$ the class of all Banach spaces $F\subset\mathcal{S}^{\prime}$
for which there exists $\sigma_{F}\in\mathbb{R}$ such that $f_{\lambda}\in F$
and $\left\Vert f_{\lambda}\right\Vert _{F}\leq C\lambda^{\sigma_{F}%
}\left\Vert f\right\Vert _{F}$, for all $\lambda>0$ and $f\in F$, where $C>0$
is a constant.
\end{definition}

Replacing $\lambda$ with $\lambda^{-1}$, it is straightforward to check that
$\frac{1}{C}\lambda^{\sigma_{F}}\Vert f\Vert_{F}\leq\Vert f_{\lambda}\Vert
_{F}$. We point that in almost all cases considered by us, we have that
$\sigma_{F}<0$. If fact, we think in $\Vert\cdot\Vert_{F}$ as a way to measure
the volume below of the graphic of the function $f$, then it is at least
intuitive that, for $\sigma_{F}<0,$ the volume of $f_{\lambda}$ is being
spread out as $\lambda\rightarrow\infty$. This fact can be seen by the
inequality $\Vert f_{\lambda}\Vert_{F}\leq C\lambda^{\sigma_{F}}\Vert
f\Vert_{F}.$

Let $E$ and $E_{0}$ be Banach spaces satisfying the following conditions:

\begin{enumerate}
\item[\textbf{(H1)}] $E,E_{0}\in\mathcal{G}$.

\item[\textbf{(H2)}] The Riesz transform $\mathcal{R}_{j}:E_{0}^{\prime
}\rightarrow E_{0}^{\prime},$ defined by
\begin{equation}
\mathcal{R}_{j}(f)=c_{n}\text{P.V.}\int_{\mathbb{R}^{n}}\frac{y_{j}}%
{|y|^{n+1}}f(x-y)dy,\text{ for }j=1,\ldots,n\text{,}
\label{eq:riesz-transform}%
\end{equation}
is bounded, where $c_{n}=\frac{\Gamma\left(  \frac{n+1}{2}\right)  }%
{\pi^{(n+1)/2}}$.

\item[\textbf{(H3)}] For every $f,g\in E^{\prime}$, the product $f\cdot g\in
E_{0}^{\prime}$; and, for some universal constant $C>0$, we have the estimate
\begin{equation}
\left\Vert fg\right\Vert _{E_{0}^{\prime}}\leq C\left\Vert f\right\Vert
_{E^{\prime}}\left\Vert g\right\Vert _{E^{\prime}}. \label{holder-dual-1}%
\end{equation}

\item[\textbf{(H4)}] Denote by $\nabla U(t):E\rightarrow E_{0}$ the operator
\[
\nabla U(t)f=(\nabla\Phi(x,t))\ast f,
\]
where
\begin{equation}
\Phi(t,x)=(4\pi t)^{-n/2}e^{-\frac{|x|^{2}}{4t}} \label{eq:Heat Kernel}%
\end{equation}
is the heat kernel. We assume that $\nabla U(t)$ is a bounded operator and
that the estimate
\begin{equation}
\int_{0}^{\infty}s^{\frac{1}{2}\left(  \sigma_{E_{0}}-\sigma_{E}\right)
-\frac{1}{2}}\left\Vert \nabla U\left(  s\right)  f\right\Vert _{E_{0}}ds\leq
C\left\Vert f\right\Vert _{E} \label{eq:integral estimate}%
\end{equation}
holds true.
\end{enumerate}

With respect to the condition \textbf{(H3)}, we have the following remark.

\begin{remark}
\label{remark:scaling-condit-to-holder-type-inequa} Let $F,G,H\in\mathcal{G}$
be such that for every $f\in F$ and $g\in G$ we have that $f\cdot g\in H$ and
\[
\Vert f\cdot g\Vert_{H}\leq C\Vert f\Vert_{F}\Vert g\Vert_{G}.
\]
Then necessarily $\sigma_{H}=\sigma_{F}+\sigma_{G}$.
\end{remark}
\noindent In fact, by replacing $f,g$ with $f_{\lambda},g_{\lambda}$ in the previous
inequality, we arrive to
\[
\lambda^{\sigma_{H}-\left(  \sigma_{F}+\sigma_{G}\right)  }\Vert f\cdot
g\Vert_{H}\leq C\Vert f\Vert_{F}\Vert g\Vert_{G}.
\]
So, if $\sigma_{F}+\sigma_{G}-\sigma_{H}\neq0$, we can take either the limit
as $\lambda\rightarrow0$ or $\lambda\rightarrow\infty$ to get a contradiction.

\noindent Now, recall that if $\mathbf{U}:F\rightarrow G$ is a bounded
operator between two normed spaces, then the dual operator $\mathbf{U}%
^{\prime}:G^{\prime}\rightarrow F^{\prime},$ defined by
\[
\left\langle \mathbf{U}^{\prime}g^{\prime},f\right\rangle =\left\langle
g^{\prime},\mathbf{U}f\right\rangle ,
\]
is also bounded. Let $F,G\in\mathcal{G}$ be such that $\nabla
U(t):F\rightarrow G$ is bounded. It is well known that the dual operator of
$\nabla^{m}U(t):F\rightarrow G$ is given by $\left(  \nabla^{m}U(t)\right)
^{\prime}=(-1)^{m}\nabla^{m}U(t)$, for each $m\in\mathbb{N}_{0}$. The
following result employs this fact and is based on a duality argument in order
to show the $L^{\infty}$-boundedness of the bilinear form
(\ref{eq:bilin-form-N-S1}).

\noindent

\begin{lemma}
\label{Lem: Estimativa integral espacio dual}

Let $E$ and $E_{0}$ verify \textbf{(H1)} and \textbf{(H4)} and assume that
$\sigma_{E_{0}}-\sigma_{E}-1=0.$ Given $f\in L^{\infty}((0,\infty
);E_{0}^{\prime})$, define the linear functional $\mathcal{T}(f)\in E^{\prime
}$ by
\[
\left\langle \mathcal{T}(f),h\right\rangle =-\int_{0}^{\infty}\left\langle
\nabla U\left(  s\right)  f,h\right\rangle ds,\text{ for all }h\in E\text{.}%
\]
Then, there exists a constant $C>0$ such that
\begin{equation}
\left\Vert \mathcal{T}\left(  f\right)  \right\Vert _{E^{\prime}}\leq
C\mathop{\sup}\limits_{t>0}\left\Vert f\left(  t\right)  \right\Vert
_{E_{0}^{\prime}}, \label{aux-dual-est-2}%
\end{equation}
for all $f\in L^{\infty}((0,\infty);E_{0}^{\prime}).$
\end{lemma}

\noindent\textbf{Proof.} First, using duality, we obtain that%
\begin{align}
\left\Vert \mathcal{T}\left(  f\right)  \right\Vert _{E^{\prime}}  &
=\mathop{\sup}\limits_{\left\Vert h\right\Vert _{E}=1}\left\vert \left\langle
\mathcal{T}\left(  f\right)  ,h\right\rangle \right\vert \nonumber\\
&  \leq\mathop{\sup}\limits_{\left\Vert h\right\Vert _{E}=1}\int_{0}^{\infty
}\left\vert \left\langle \nabla U\left(  s\right)  f\left(  s\right)
,h\right\rangle \right\vert ds\nonumber\\
&  =\mathop{\sup}\limits_{\left\Vert h\right\Vert _{E}=1}\int_{0}^{\infty
}\left\vert \left\langle f\left(  s\right)  ,\nabla U\left(  s\right)
h\right\rangle \right\vert ds\nonumber\\
&  \leq\mathop{\sup}\limits_{\left\Vert h\right\Vert _{E}=1}\int_{0}^{\infty
}\left\Vert f\left(  s\right)  \right\Vert _{E_{0}^{\prime}}\left\Vert \nabla
U\left(  s\right)  h\right\Vert _{E_{0}}ds\nonumber\\
&  \leq\mathop{\sup}\limits_{t>0}\left\Vert f\left(  t\right)  \right\Vert
_{E_{0}^{\prime}}\mathop{\sup}\limits_{\left\Vert h\right\Vert _{E}=1}\int
_{0}^{\infty}\left\Vert \nabla U\left(  s\right)  h\right\Vert _{E_{0}}ds.
\label{aux-dual-est-1}%
\end{align}
Now, the condition \textbf{(H4)} and $\sigma_{E_{0}}-\sigma_{E}-1=0$ lead us
to%
\begin{align*}
\text{R.H.S. of (\ref{aux-dual-est-1})}  &  =\mathop{\sup}\limits_{t>0}%
\left\Vert f\left(  t\right)  \right\Vert _{E_{0}^{\prime}}%
\mathop{\sup}\limits_{\left\Vert h\right\Vert _{E}=1}\int_{0}^{\infty}%
s^{\frac{1}{2}\left(  \sigma_{E_{0}}-\sigma_{E}\right)  -\frac{1}{2}%
}\left\Vert \nabla U\left(  s\right)  h\right\Vert _{E_{0}}ds\\
&  \leq C\mathop{\sup}\limits_{t>0}\left\Vert f\left(  t\right)  \right\Vert
_{E_{0}^{\prime}}\mathop{\sup}\limits_{\left\Vert h\right\Vert _{E}%
=1}\left\Vert h\right\Vert _{E}\\
&  =C\mathop{\sup}\limits_{t>0}\left\Vert f\left(  t\right)  \right\Vert
_{E_{0}^{\prime}},
\end{align*}
which, together with (\ref{aux-dual-est-1}), yield the desired estimate. \fin

With the estimate (\ref{aux-dual-est-2}) in hand, we are in position to show
the bilinear estimate for (\ref{eq:bilin-form-N-S1}).

\begin{theorem}
\label{Teo: Teo main. Estimativa integral en el dual.}

Let $E$ and $E_{0}$ verify \textbf{(H1)}-\textbf{(H4)} and assume that
$\sigma_{E_{0}}-\sigma_{E}-1=0\,.$ Then, we have the bilinear estimate%
\begin{equation}
\mathop{\sup}\limits_{0<t<T}\left\Vert \mathcal{B}(u,v)(t)\right\Vert
_{E^{\prime}}\leq K\mathop{\sup}\limits_{0<t<T}\left\Vert u\left(  t\right)
\right\Vert _{E^{\prime}}\mathop{\sup}\limits_{0<t<T}\left\Vert v\left(
t\right)  \right\Vert _{E^{\prime}}, \label{eq:Estimativa bilinear}%
\end{equation}
for all $u,v\in L^{\infty}((0,T);E^{\prime})$ and $T\in(0,\infty]$, where
$K>0$ is a universal constant.
\end{theorem}

\begin{remark}
\label{Rem-exist-2} In the case of Navier-Stokes equations, the conditions in
Lemma \ref{Lem: Estimativa integral espacio dual} lead to the space
$E^{\prime}$ to be critical for those equations, that is, $\sigma_{E^{\prime}%
}=1$. However, we prefer to maintain the conditions that force to take
$\sigma_{E^{\prime}}=1$ separately and independently of the others, such as
the balance scaling condition $\sigma_{E_{0}}-\sigma_{E}-1=0.$ A reason is
that Lemma \ref{Lem: Estimativa integral espacio dual} is useful by itself and
depends on the latter condition which gives us more versatility and
applicability of the theory to other evolution PDEs (see Section
\ref{section:related-models}).
\end{remark}

\noindent\textbf{Proof of Theorem
\ref{Teo: Teo main. Estimativa integral en el dual.}.}

Let $0<T\leq\infty$ and $t\in\left(  0,T\right)  .$ The bilinear term
$\mathcal{B}(u,v)$ can be written as
\[
\mathcal{B}(u,v)(t)=-\int\limits_{0}^{t}\nabla_{x}U\left(  t-s\right)
\mathbb{P}f(\cdot,s)ds=\mathcal{T}(f_{t}),
\]
where $f_{t}(x,s)$ is defined by
\begin{align*}
f_{t}(\cdot,s)  &  =\mathbb{P}(u\otimes v)\left(  \cdot,t-s\right)
,\,\mbox{a.e.}\,s\in\left(  0,t\right)  ,\\
f_{t}(\cdot,s)  &  =0,\,\text{a.e. }s\in\left(  t,\infty\right)  .
\end{align*}
It follows from Lemma \ref{Lem: Estimativa integral espacio dual} that%

\begin{equation}
\left\Vert \mathcal{B}(u,v)(t)\right\Vert _{E^{\prime}}=\left\Vert \mathcal{T}\left(
f_{t}\right)  \right\Vert _{E^{\prime}}\leq C\mathop{\sup}\limits_{s>0}%
\left\Vert f_{t}\left(  s\right)  \right\Vert _{E_{0}^{\prime}}.
\label{aux-bili-proof-1}%
\end{equation}
Moreover, using \textbf{(H2)} and \textbf{(H3)}, we obtain that
\begin{align}
\mathop{\sup}\limits_{0<s<T}\left\Vert f_{t}(s)\right\Vert _{E_{0}^{\prime}}
&  \leq C\mathop{\sup}\limits_{0<s<t<T}\left\Vert (u\otimes v)\left(
\cdot,t-s\right)  \right\Vert _{E_{0}^{\prime}}\nonumber\\
&  \leq C\mathop{\sup}\limits_{0<s<t<T}\left\Vert u\left(  \cdot,t-s\right)
\right\Vert _{E^{\prime}}\left\Vert v\left(  \cdot,t-s\right)  \right\Vert
_{E^{\prime}}\nonumber\\
&  \leq C\mathop{\sup}\limits_{0<s<T}\left\Vert u\left(  \cdot,s\right)
\right\Vert _{E^{\prime}}\mathop{\sup}\limits_{0<s<T}\left\Vert v\left(
\cdot,s\right)  \right\Vert _{E^{\prime}}. \label{aux-bili-proof-2}%
\end{align}
Estimate (\ref{eq:Estimativa bilinear}) follows by inserting
(\ref{aux-bili-proof-2}) into (\ref{aux-bili-proof-1}). \fin

The next result gives sufficient conditions to ensure the condition
\textbf{(H4)}. In fact, inequality (\ref{eq:integral estimate}) can be seen as
one of the main ingredients to obtain the bilinear estimate for the
Navier-Stokes equations. In general lines, it is a consequence of a careful
use of real interpolation tools.

\begin{lemma}
\label{Lem: condiciones para estimativa integral} Let $E,E_{0},E_{1},E_{2}%
\in\mathcal{G}$ be Banach spaces verifying the following conditions:

\begin{enumerate}
\item[\textbf{(H4$^\prime$)}] $\left\Vert \nabla U\left(  s\right)  f\right\Vert
_{E_{0}}\leq Cs^{\frac{1}{2}\left(  \sigma_{E_{i}}-\sigma_{E_{0}}\right)
-\frac{1}{2}}\left\Vert f\right\Vert _{E_{i}},$ for $i=1,2$.

\item[\textbf{(H4$^{\prime \prime}$)}] Assume that $z_{1}$ and $z_{2}$ satisfy $0<z_{1}<1<z_{2}$
and $\frac{1}{z_{i}}=\frac{1}{2}\left(  \sigma_{E}-\sigma_{E_{i}}\right)  +1$,
for $i=1,2$.

\item[\textbf{(H4$^{\prime \prime \prime}$)}] If $\theta\in\left(  0,1\right)  $ is such that
$1=\frac{1-\theta}{z_{1}}+\frac{\theta}{z_{2}}$ with $z_{1},z_{2}$ as in
\textbf{(H4$^{\prime \prime}$)}, then $E\hookrightarrow\left(  E_{1},E_{2}\right)  _{\theta,1}$.

Then, the integral estimate (\ref{eq:integral estimate}) holds true.
\end{enumerate}
\end{lemma}

\begin{remark}
\label{Rem-lemma-H5-H7}Note that, in \textbf{(H4$^{\prime \prime}$)}, the order taken on the
spaces $E_{1}$ and $E_{2}$ is not relevant, but we keep it for convenience of exposition.
\end{remark}

\textbf{Proof of Lemma \ref{Lem: condiciones para estimativa integral}. }The
proof extends some ideas in \cite{Yam} and we sketch it for the reader
convenience (see also \cite[Lemma 5.1]{Lucas}). Considering $h_{f}\left(
s\right)  =s^{\frac{1}{2}\left(  \sigma_{E_{0}}-\sigma_{E}\right)  -\frac
{1}{2}}\left\Vert \nabla U\left(  s\right)  f\right\Vert _{E_{0}}$ and using
\textbf{(H4$^{\prime \prime}$)}, we can estimate
\[
h_{f}\left(  s\right)  \leq Cs^{\frac{1}{2}\left(  \sigma_{E_{i}}-\sigma
_{E}\right)  -1}\left\Vert f\right\Vert _{E_{i}}=Cs^{-\frac{1}{z_{i}}%
}\left\Vert f\right\Vert _{E_{i}}\,\,i=1,2.
\]
It follows that $h_{f}\in L^{z_{1},\infty}\cap L^{z_{2},\infty}$ and
\[
\left\Vert h_{f}\right\Vert _{L^{z_{i},\infty}}\leq C\left\Vert f\right\Vert
_{E_{i}}.
\]
Using interpolation and \textbf{(H4$^{\prime \prime}$)}, we conclude that
\[
\left\Vert h_{f}\right\Vert _{L^{1}}=\left\Vert h_{f}\right\Vert _{\left(
L^{z_{1},\infty},L^{z_{2},\infty}\right)  _{\theta,1}}\leq C\left\Vert
f\right\Vert _{\left(  E_{1},E_{2}\right)  _{\theta,1}}\leq C\left\Vert
f\right\Vert _{E},
\]
which is (\ref{eq:integral estimate}).

\fin

\subsection{Applications}

\label{Sec-Appl-1}In this section we present two cases where Theorem
\ref{Teo: Teo main. Estimativa integral en el dual.} applies in a more or less
direct way.

\subsubsection{Bilinear estimate on weak-$L^{p}$ spaces}

By taking $E=L^{\frac{n}{n-1},1}$ and $E_{0}=L^{\frac{n}{n-2},1},$ we obtain
the bilinear estimate (\ref{eq:Estimativa bilinear}) in the space $E^{\prime
}=L^{n,\infty}$ (see \cite{Yam}). In fact, \textbf{(H1)}-\textbf{(H3)} are
well-known and we check \textbf{(H4$^\prime$)}-\textbf{(H4$^{\prime \prime \prime}$ )}. For this, take
$E_{i}=L^{(p_{i},1)}$ with $1<p_{2}<\frac{n}{n-1}<p_{1}<\frac{n}{n-2}$, then
it follows \textbf{(H4$^\prime$)} and \textbf{(H4$^{\prime \prime}$)}. Choosing $\theta\in(0,1)$ such
that $\frac{n-1}{n}=\frac{1-\theta}{p_{1}}+\frac{\theta}{p_{2}}$, we have that
$L^{p,1}=(L^{(p_{1},1)},L^{(p_{2},1)})_{\theta,1}$, by interpolation theory,
and then it follows \textbf{(H4$^{\prime \prime \prime}$ )}.

\subsubsection{Bilinear estimate on Lorentz-Herz spaces}

Now we show that the Lorentz-Herz spaces $E=\dot{K}_{\frac{n}{n-1},1,1}^{0}$
and $E_{0}=\dot{K}_{\frac{n}{n-2},1,1}^{0}$ fulfill the conditions
\textbf{(H1)}-\textbf{(H4)}. Therefore, we obtain the bilinear estimate
(\ref{eq:Estimativa bilinear}) in the space $E^{\prime}=W\dot{K}_{n,\infty
}^{0}$.

\begin{remark}
\label{Rem-weak-Herz}It is worthy to comment that, in comparison with
\cite{Tsutsui-1}, this approach provides a new proof of
(\ref{eq:Estimativa bilinear}) in $E^{\prime}=W\dot{K}_{n,\infty}^{0}.$ In
fact, the proof in \cite{Tsutsui-1} relies on the so-called Meyer approach
\cite{Meyer}.
\end{remark}

Let $1\leq p\leq\infty$ and $1\leq d,q\leq\infty,$ with $d=\infty$ if
$p=\infty.$ Recall that the Lorentz-Herz space $\dot{K}_{\left(  p,d\right)
,q}^{\alpha}$ is defined as a class of functions $f$ such that
\[
\left\Vert f\right\Vert _{\dot{K}_{\left(  p,d\right)  ,q}^{\alpha}%
}=\left\Vert 2^{k\alpha}\left\Vert f\right\Vert _{L^{p,d}\left(  A_{k}\right)
}\right\Vert _{l^{q}\left(  \mathbb{Z}\right)  }<\infty,
\]
where the sets $A_{k}$ are defined as
\[
A_{k}=\left\{  x\in\mathbb{R}^{n};2^{k-1}\leq x<2^{k}\right\}  .
\]

Next we collect some properties about Lorentz-Herz spaces that will be useful
in the sequel. For more details, see \cite{Tsutsui-1, LucJhe1}.

\begin{enumerate}
\item[\textbf{(a)}] $\left\Vert f(\lambda\cdot)\right\Vert _{\dot{K}
_{p,d,q}^{0}}\approx\lambda^{-\frac{n}{p}}\left\Vert f\right\Vert _{\dot
{K}_{p,d,q}^{0}},$ $\forall\lambda>0.$%
\label{Propiedad a) espacios de Lorentz Herz}

\item[\textbf{(b)}] $\left(  \dot{K}_{p^{\prime},d^{\prime},q^{\prime}}%
^{0}\right)  ^{\prime}=\dot{K}_{p,d,q}^{0},$ for $1<p<\infty$ and
$1<d,q\leq\infty.$\label{Propiedad b) espacios de Lorentz Herz Dualidad}

\item[\textbf{(c)}] The Riesz transform is bounded in $\dot{K}_{p,d,q}^{0},$
for $1<p<\infty$ and $1\leq d,q\leq\infty.$

\item[\textbf{(d)}] We have the H\"{o}lder-type inequality
\label{Propiedad c) espacios de Lorentz Herz Desigualdad de Holder}
\[
\left\Vert fg\right\Vert _{\dot{K}_{p,d,q}^{0}}\leq C\left\Vert f\right\Vert
_{\dot{K}_{p_{1},d_{1},q_{1}}^{0}}\left\Vert g\right\Vert _{\dot{K}%
_{p_{2},d_{2},q_{2}}^{0}},
\]
where $\frac{1}{p}=\frac{1}{p_{1}}+\frac{1}{p_{2}}$, $\frac{1}{d}=\frac
{1}{d_{1}}+\frac{1}{d_{2}}$ and $\frac{1}{q}=\frac{1}{q_{1}}+\frac{1}{q_{2}}$.

\item[\textbf{(e)}] $\left\Vert \nabla U\left(  s\right)  f\right\Vert
_{\dot{K}_{p,d,q}^{0}}\leq Cs^{\frac{1}{2}\left(  \frac{n}{p_{0}}-\frac{n}%
{p}\right)  -\frac{1}{2}}\left\Vert f\right\Vert _{\dot{K}_{p_{0},d,q}^{0}},$
$\forall s>0$, for $1<p_{0}\leq p<\infty$ and $1\leq d,q\leq\infty
$.\label{Propiedad e) espacios de Lorentz Herz Estimativa nucleo del calor}

\item[\textbf{(f)}] $\left\Vert \psi\ast f\right\Vert _{L^{\infty}}\leq
C(\psi)\left\Vert f\right\Vert _{\dot{K}_{p,d,q}^{0}},$ for $\psi
\in\mathcal{S}$%
.\label{Propiedad f) espacios de Lorentz Herz convolucion desde Linfty hasta Lorentz Herz}%

\item[\textbf{(g)}] $\left\Vert \psi\ast f\right\Vert _{\dot{K}_{p,d,q}^{0}%
}\leq C\max\left\{  \left\Vert \psi\right\Vert _{L^{1}},\left\Vert \left\vert
\cdot\right\vert ^{n}\psi\right\Vert _{L^{\infty}}\right\}  \left\Vert
f\right\Vert _{\dot{K}_{p,d,q}^{0}},$ for $\psi\in L^{1}(\mathbb{R}^{n})\cap
L^{\infty}(\mathbb{R}^{n},\left\vert \cdot\right\vert ^{n}dx)$%
.\label{Propiedad g) espacios de Lorentz Herz convolucion desde Lorentz manteniendo Lorentz}

\item[\textbf{(h)}] $\left\Vert \Delta_{0}f\right\Vert _{\dot{K}_{p,d,q}^{0}
}\leq C\left\Vert f\right\Vert _{\dot{K}_{p_{0},d,q}^{0}},$ for $1<p_{0}\leq
p<\infty$ and $1\leq d,q\leq\infty$. (see \eqref{eq:OperadorDeltaj} for the definition of the operator $\Delta_{0}$) \label{Propiedad e) espacios de Lorentz Herz Estimativa Delta cero}
\end{enumerate}

In what follows, we are going to prove (\ref{eq:integral estimate}) in the
context of Lorentz-Herz spaces. We start with the lemma below.

\begin{lemma}
\label{Lem: Desigualdade de interpolacao Lorentz-Herz} Let $1<p,p_{1}%
,p_{2}\leq\infty$ and $1\leq d,d_{1},d_{2},q_{1},q_{2}\leq\infty$ be such that
$\frac{1}{p}=\frac{1-\theta}{p_{1}}+\frac{\theta}{p_{2}}$ with $\theta
\in\left(  0,1\right)  .$ If $X_{1}$ and $X_{2}$ are Banach spaces and $T$ is
a continuous linear operator such that
\begin{equation}
T:\dot{K}_{\left(  p_{1},d_{1}\right)  ,q_{1}}^{0}\longrightarrow
X_{1}\,\,\,\mbox{and}\,\,\,T:\dot{K}_{\left(  p_{2},d_{2}\right)  ,q_{2}}%
^{0}\longrightarrow X_{2}, \label{eq: aux operador}%
\end{equation}
with the operator norms $C_{1}$ and $C_{2},$ respectively. Then
\[
T:\dot{K}_{\left(  p,d\right)  ,1}^{0}\longrightarrow\left(  X_{1}%
,X_{2}\right)  _{\theta,d},
\]
with the operator norm bounded by $\tilde{C}=\left(  C_{1}\right)  ^{1-\theta
}\left(  C_{2}\right)  ^{\theta}.$
\end{lemma}

\noindent\textbf{Proof.} For every $k\in\mathbb{Z}$, it follows from
\eqref{eq: aux operador} that
\begin{equation}
T:L^{p_{1},d_{1}}\left(  A_{k}\right)  \longrightarrow X_{1}%
\,\,\,\mbox{and}\,\,\,T:L^{p_{2},d_{2}}\left(  A_{k}\right)  \longrightarrow
X_{2}. \label{eq:aux operator in Lpi,di}%
\end{equation}
Denoting $N_{1}^{k}$ and $N_{2}^{k}$ the respective operator norms in
\eqref{eq:aux operator in Lpi,di}, we get $N_{1}^{k}\leq C_{1}$ and $N_{2}%
^{k}\leq C_{2}$. Next, using interpolation in Lorentz spaces (see \cite{BL}),
we obtain that
\begin{equation}
T\,:\,L^{p,d}\left(  A_{k}\right)  \longrightarrow\left(
X_{1},X_{2}\right)  _{\theta,d}, \label{eq:operator in Lp,d}%
\end{equation}
with operator norm $N^{k}\leq\left(  N_{1}^{k}\right)  ^{1-\theta}\left(
N_{2}^{k}\right)  ^{\theta}\leq\left(  C_{1}\right)  ^{1-\theta}\left(
C_{2}\right)  ^{\theta}=\tilde{C}$.

Now consider the dense subspace of $\dot{K}_{\left(  p,d\right)  ,1}^{0}$
given by

\[
Y=\left\{  f^{m}=\sum\limits_{k=-m}^{m}\chi_{A_{k}}f;m\in\mathbb{N}%
\,\mbox{and}\,f\in\dot{K}_{\left(  p,d\right)  ,1}^{0}\right\}  .
\]
For $f^{m}$ in $Y,$ we have that
\begin{align*}
\left\Vert T\left(  f^{m}\right)  \right\Vert _{\left(  X_{0}%
,X_{1}\right)  _{\theta,d}}  &  =\sum\limits_{k=-m}^{m}\left\Vert
T\left(  \chi_{A_{k}}f\right)  \right\Vert _{\left(  X_{0}%
,X_{1}\right)  _{\theta,d}}\leq\sum\limits_{k=-m}^{m}\tilde{C}\left\Vert
f\right\Vert _{L^{p,d}\left(  A_{k}\right)  }\\
&  =\sum\limits_{k=-m}^{m}\tilde{C}\left\Vert f^{m}\right\Vert _{L^{p,d}%
\left(  A_{k}\right)  }\leq\tilde{C}\left\Vert f^{m}\right\Vert _{\dot
{K}_{\left(  p,d\right)  ,1}^{0}},
\end{align*}
and then we conclude the proof by density. \fin

\begin{remark}
Note that, under the conditions of Lemma
\ref{Lem: Desigualdade de interpolacao Lorentz-Herz}, we have the continuous inclusion%

\begin{equation}
\dot{K}_{\left(  p,d\right)  ,1}^{0}\hookrightarrow\left(  \dot{K}_{\left(
p_{1},d_{1}\right)  ,q_{1}}^{0},\dot{K}_{\left(  p_{2},d_{2}\right)  ,q_{2}%
}^{0}\right)  _{\theta,d}.
\label{eq:inclusion lorentz herz en interpolado de lorentz herz}%
\end{equation}
In particular,%

\begin{equation}
\dot{K}_{\left(  p,1\right)  ,1}^{0}\hookrightarrow\left(  \dot{K}_{\left(
p_{1},1\right)  ,1}^{0},\dot{K}_{\left(  p_{2},1\right)  ,1}^{0}\right)
_{\theta,1}.
\label{eq:inclusion lorentz herz en interpolado de lorentz herz d=00003D00003D00003D00003D1}%
\end{equation}

\end{remark}

Now, we turn to complete the proof of the bilinear estimate (\ref{eq:Estimativa bilinear}) in
$W\dot{K}_{n,\infty}^{0}$. Let $p_{1}$ and $p_{2}$ be such that $0<p_{1}%
<\frac{n}{n-1}<p_{2}<\frac{n}{n-2}$ and define $z_{i}$ as $\frac{1}{z_{i}%
}=\frac{1}{2}\left(  -\frac{n}{\frac{n}{n-1}}-\left(  -\frac{n}{p_{i}}\right)
\right)  +1,$ for $i=1,2$. Note that $0<z_{2}<1<z_{1}<\infty$, and then we can
take $\theta\in\left(  0,1\right)  $ in such a way that $1=\frac{\theta}%
{z_{1}}+\frac{1-\theta}{z_{2}}$ and $\frac{1}{\frac{n}{n-1}}=\frac{\theta
}{p_{1}}+\frac{1-\theta}{p_{2}}$. Consider the spaces $E=\dot{K}_{\frac
{n}{n-1},1,1}^{0}$, $E_{0}=\dot{K}_{\frac{n}{n-2},1,1}^{0}$, $E_{1}=\dot
{K}_{p_{1},1,1}^{0}$ and $E_{2}=\dot{K}_{p_{2},1,1}^{0}$ with $p_{1}$ and
$p_{2}$ as before. Using these choices, we get directly \textbf{(H4$^{\prime \prime}$)}. The
property $\mathbf{(e)}$ implies \textbf{(H4$^\prime$)} and the inclusion
\eqref{eq:inclusion lorentz herz en interpolado de lorentz herz d=00003D00003D00003D00003D1}
gives \textbf{(H4$^{\prime \prime \prime}$ )}. By Lemma \ref{Lem: condiciones para estimativa integral},
we get (\ref{eq:integral estimate}) and then \textbf{(H4)}. Conditions\textbf{
(H1)}, \textbf{(H2)} and \textbf{(H3)} follow from properties \textbf{(a)},
\textbf{(c)} and \textbf{(d)}, respectively.

With the above properties in hand and applying Theorem
\ref{Teo: Teo main. Estimativa integral en el dual.}, we obtain the bilinear
estimate (\ref{eq:Estimativa bilinear}) with $E^{\prime}=\dot{K}_{\left(
n,\infty\right)  ,\infty}^{0}:=W\dot{K}_{n,\infty}^{0}$; that is, we reobtain
the bilinear estimate proved in \cite{Tsutsui-1} by using a different method.

\section{Bilinear estimates in Besov-type spaces}

\label{sect:bilinear-estimate-besov-spaces}

In what follows, the functions $\varphi$ and $\Theta$ are radially symmetric and
satisfy (see, e.g., \cite{BL})
\[
\varphi\in C_{c}^{\infty}\left(  \mathbb{R}^{n}\backslash\left\{  0\right\}
\right)  ,\,\mbox{supp}\,\varphi\subset\left\{  x\,;\,\frac{3}{4}%
\leq\left\vert x\right\vert \leq\frac{8}{3}\right\}  ,
\]%
\[
\sum\limits_{j\in\mathbb{Z}}\varphi_{j}(\xi)=1,\,\forall\xi\in\mathbb{R}%
^{n}\backslash\left\{  0\right\}  ,\,\mbox{where}\,\varphi_{j}(\xi
):=\varphi\left(  2^{-j}\xi\right)  ,
\]
and
\[
\Theta\left(  \xi\right)  =\left\{
\begin{array}
[c]{l}%
\sum\limits_{j\leq-1}\varphi_{j}(\xi),\,\mbox{if}\,\xi\neq0\\
1,\,\,\,\,\,\,\,\,\,\,\,\,\,\,\,\,\,\,\,\,\,\,\,\mbox{if}\,\xi=0.
\end{array}
\right.
\]
Thus, $\Theta\in C_{c}^{\infty}\left(  \mathbb{R}^{n}\right)  ,$
$\mbox{supp}\,\Theta\subset\left\{  x\,;\,\left\vert x\right\vert \leq\frac
{4}{3}\right\}  $ and
\[
\Theta\left(  \xi\right)  +\sum\limits_{j\geq0}\varphi_{j}(\xi)=1,\,\forall
\xi\in\mathbb{R}^{n}.
\]
Moreover, we also define $\tilde{\varphi}_{j}=\sum\limits_{\left\vert
k-j\right\vert \leq1}\varphi_{k}$, $\tilde{\Theta}=\Theta+\varphi$ and
$\tilde{D}_{j}=D_{j-1}\cup D_{j}\cup D_{j+1}$, where $D_{j}=\left\{
x;\frac{3}{4}2^{j}\leq\left\vert x\right\vert \leq\frac{8}{3}2^{j}\right\}  $
and $j\in\mathbb{Z}$. So, we have that $\tilde{\varphi}_{j}\equiv1$ on $D_{j}%
$, $\varphi_{j}=\tilde{\varphi}_{j}\varphi_{j}$ for all $j\in\mathbb{Z}$ and
$\Theta=\tilde{\Theta}\Theta$.

The operators $\Delta_{j}$ and $S_{j}$ are defined as%

\begin{align}
	\Delta_{j}f=\varphi_{j}(D)f & =\left(\varphi_{j}\widehat{f}\right)\widecheck{\ }\,\,\,\,\forall j\in\mathbb{Z},\label{eq:OperadorDeltaj}\\
	S_{j}f=\Theta_{j}(D)f & =\left(\Theta_{j}\widehat{f}\right)\widecheck{\ }\,\,\,\,\forall j\in\mathbb{Z},\label{eq:OperadorSj}
\end{align}
where $\Theta_{j}\left(  \xi\right)  =\Theta\left(  2^{-j}\xi\right)  .$ It is
not hard to show the cancellation identities
\[
\Delta_{j}\Delta_{k}f=0\,\,\mbox{if}\,\,\left\vert j-k\right\vert \geq2,\text{
and }\Delta_{j}\left(  S_{k-2}g\Delta_{k}f\right)
=0\,\,\mbox{if}\,\,\left\vert j-k\right\vert \geq3.
\]

Now we recall a general definition of Besov-type spaces based on a Banach
space $F$. This definition has been used by some authors, see e.g.
\cite{Furi-Lema-1,Lemarie1,Maz}.

\begin{definition}
\label{Def: Definicion espacios Besov baseado en E} Let $F\subset
\mathcal{S}^{\prime}$ be a Banach space, $1\leq r\leq\infty$ and
$s\in\mathbb{R}$. The homogeneous Besov-$F$ space, denoted by $\dot{B}\left[
F\right]  _{r}^{s},$ is defined as%

\[
\dot{B}\left[  F\right]  _{r}^{s}=\left\{  f\in\mathcal{S}^{\prime}%
(\mathbb{R}^{n})/\mathcal{P};\,\left\Vert f\right\Vert _{\dot{B}\left[
F\right]  _{r}^{s}}<\infty\right\}  ,
\]
where $\mathcal{P}$ stands for the set of all polynomials and%

\begin{equation}
\left\Vert f\right\Vert _{\dot{B}\left[  F\right]  _{r}^{s}}:=%
\begin{cases}
\left(  \sum\limits_{j\in\mathbb{Z}}2^{jsr}\left\Vert \Delta_{j}f\right\Vert
_{F}^{r}\right)  ^{\frac{1}{r}} & \mbox{if}\,\,r<\infty\\
\mathop{\sup}\limits_{j\in\mathbb{Z}}\,2^{js}\left\Vert \Delta_{j}f\right\Vert
_{F} & \mbox{if}\,\,r=\infty.
\end{cases}
\label{eq: Norma Besov-X homogeneo}%
\end{equation}

\end{definition}

The homogeneous Besov-$F$ space is Banach provided that the underlying space
$F$ verifies some basic conditions (see \cite{Karch},\cite{Maz} for related
ones). This is the subject of the next lemma.

\begin{lemma}
\label{Lem: Besov-X Banach space of distributions} Let $F$ be a Banach space
verifying the conditions

\begin{description}
\item[\textbf{(B1)}] $F\in\mathcal{G}$ ;

\item[\textbf{(B2)}] $\left\Vert \Delta_{0}f\right\Vert _{L^{\infty}}\leq
C\left\Vert f\right\Vert _{F}.$
\end{description}

Then, $\dot{B}\left[  F\right]  _{r}^{s}$ is a Banach space and $\dot
{B}\left[  F\right]  _{r}^{s}\hookrightarrow\mathcal{S}^{\prime}/\mathcal{P}$.
\end{lemma}

\textbf{Proof.} Using $\Delta_{j}f=\Delta_{0}\left(  f\left(  \frac{\cdot
}{2^{j}}\right)  \right)  \left(  2^{j}\cdot\right)  $, we can see that
\begin{align*}
\left\Vert \Delta_{j}f\right\Vert _{L^{\infty}}  &  =\left\Vert \Delta
_{0}\left(  f\left(  \frac{\cdot}{2^{j}}\right)  \right)  \left(  2^{j}%
\cdot\right)  \right\Vert _{L^{\infty}}\leq C\left\Vert \Delta_{0}\left(
f\left(  \frac{\cdot}{2^{j}}\right)  \right)  \right\Vert _{L^{\infty}}\\
&  \leq C\left\Vert f\left(  \frac{\cdot}{2^{j}}\right)  \right\Vert _{F}\leq
C2^{-j\sigma_{F}}\left\Vert f\right\Vert _{F}.
\end{align*}
The same arguments lead us to
\[
\left\Vert \Delta_{j}f\right\Vert _{L^{\infty}}=\left\Vert \tilde{\Delta}%
_{j}\Delta_{j}f\right\Vert _{L^{\infty}}\leq C2^{-j\sigma_{F}}\left\Vert
\Delta_{j}f\right\Vert _{F},
\]
and then $\dot{B}\left[  F\right]  _{r}^{s}\hookrightarrow\dot{B}_{\infty
,r}^{\gamma}$ with $\gamma=s+\sigma_{F}$. The rest of the proof follows the
same ideas in \cite[Lemma 4.2]{Karch} and \cite[Lemma 2.21]{Maz}. \fin

\begin{remark}
Assume that $F,G\in\mathcal{G}$ present the property$\ \left\Vert \Delta
_{0}f\right\Vert _{G}\leq C\left\Vert f\right\Vert _{F},$ $\forall f\in F$.
Proceeding as in the previous proof, we obtain that
\begin{align*}
\left\Vert \Delta_{j}f\right\Vert _{G}  &  =\left\Vert \Delta_{0}\left(
f\left(  \frac{\cdot}{2^{j}}\right)  \right)  \left(  2^{j}\cdot\right)
\right\Vert _{G}\leq C2^{j\sigma_{G}}\left\Vert \Delta_{0}\left(  f\left(
\frac{\cdot}{2^{j}}\right)  \right)  \right\Vert _{G}\\
&  \leq C2^{j\sigma_{G}}\left\Vert f\left(  \frac{\cdot}{2^{j}}\right)
\right\Vert _{F}\leq C2^{j\left(  \sigma_{G}-\sigma_{F}\right)  }\left\Vert
f\right\Vert _{F}.
\end{align*}

\end{remark}

Additionally, consider the following condition on the space $F$ which is an
extension of the convolution property in shift-invariant spaces of
distributions. In fact, taking $M=L^{1}(\mathbb{R}^{n})$ below, we obtain the
property in \cite[Proposition 4.1]{Lemarie1}.

\begin{description}
\item[\textbf{(B3)}] $\left\Vert \psi\ast f\right\Vert _{F}\leq C\left\Vert
\psi\right\Vert _{M}\left\Vert f\right\Vert _{F},$ for all $f\in F$ and
$\psi\in\mathcal{S}$, where the norm $\left\Vert \cdot\right\Vert _{M}$ is
such that $\left\Vert 2^{jn}\psi(2^{j}\cdot)\right\Vert _{M}\approx\left\Vert
\psi\right\Vert _{M}$ and $\left\Vert \psi\right\Vert _{M}\leq C\text{$\max$%
}\left\{  \left\Vert \psi\right\Vert _{L^{1}},\left\Vert \left\vert
\cdot\right\vert ^{n}\psi\right\Vert _{L^{\infty}}\right\}  $.
\end{description}

\begin{lemma}
\label{Lem: aux Operador pseudo X} Let $F\subset\mathcal{S}^{\prime}$ be a
Banach space verifying \textbf{(B3)}, $m\in\mathbb{R},$ and $P$ a $C^{\infty
}$-function on $\tilde{D}_{j}$ such that $\left\vert \partial^{\rho}P\left(
\xi\right)  \right\vert \leq C2^{\left(  m-\left\vert \rho\right\vert \right)
j},$ $\forall\xi\in\tilde{D}_{j}$ and $\forall\rho\in\mathbb{N}_{0}^{n}$ with
$\left\vert \rho\right\vert \leq n$. Then,%

\[
\left\Vert \left(  P\widehat{f}\right)\widecheck{\ }\right\Vert _{F}\leq
C2^{jm}\left\Vert f\right\Vert _{F},
\]
for all $f\in F$ such that $\mbox{supp}\hat{f}\subset D_{j}$.
\end{lemma}

\noindent\textbf{Proof.} Define $K(x)=\left(  P(\cdot)\tilde{\varphi}%
_{j}(\cdot)\right)  \widecheck{}(x)$, where $\tilde{\varphi}_{j}=\varphi
_{j-1}+\varphi_{j}+\varphi_{j+1}$. Since $\mbox{supp}\widehat{f}\subset D_{j}$, it
follows that $P\left(  \xi\right)  \widehat{f}\left(  \xi\right)  =P\left(
\xi\right)  \tilde{\varphi}_{j}(\xi)\widehat{f}\left(  \xi\right)  $. So, $\left(
P\widehat{f}\right)  \widecheck{}=\left(  P\tilde{\varphi}_{j}\widehat{f}\right)  \widecheck
{}=K\ast f$. Using \textbf{(B3)}, we obtain that
\[
\left\Vert \left(  P\widehat{f}\right)  \widecheck{\ }\right\Vert _{F}\leq\left\Vert
K\right\Vert _{M}\left\Vert f\right\Vert _{F}.
\]
We claim that $\left\Vert K\right\Vert _{M}\leq C2^{jm}$. For this, we choose
$N\in\mathbb{N}$ such that $n/2<N\leq\left[  n/2\right]  +1$ in order to get%

\begin{align*}
\left\Vert K\right\Vert _{L^{1}}  &  =\int\limits_{D\left(  0,2^{-j}\right)
}K\left(  y\right)  +\int\limits_{\left\vert y\right\vert \geq2^{-j}}K\left(
y\right) \\
&  \leq\left(  \int\limits_{D\left(  0,2^{-j}\right)  }1\right)  ^{1/2}\left(
\int\limits_{D\left(  0,2^{-j}\right)  }\left\vert K\left(  y\right)
\right\vert ^{2}\right)  ^{1/2}+\left(  \int\limits_{\left\vert y\right\vert
\geq2^{-j}}\left\vert y\right\vert ^{-2N}\right)  ^{1/2}\left(  \int
\limits_{\left\vert y\right\vert \geq2^{-j}}\left\vert y\right\vert
^{2N}\left\vert K\left(  y\right)  \right\vert ^{2}\right)  ^{1/2}\\
&  \leq C2^{-j\left(  \frac{n}{2}\right)  }\left\Vert P\tilde{\varphi}%
_{j}\right\Vert _{L^{2}}+C2^{-j\left(  -N+\frac{n}{2}\right)  }\sum
\limits_{\left\vert \rho\right\vert =N}\left\Vert \left(  \cdot\right)
^{\rho}K\right\Vert _{L^{2}}\\
&  \leq C2^{-j\left(  \frac{n}{2}\right)  }\left\Vert P\tilde{\varphi}%
_{j}\right\Vert _{L^{2}}+C2^{-j\left(  -N+\frac{n}{2}\right)  }\sum
\limits_{\left\vert \rho\right\vert =N}\left\Vert \partial^{\rho}\left(
P\tilde{\varphi}_{j}\right)  \right\Vert _{L^{2}}\\
&  \leq C2^{jm}+C2^{jm}\leq C2^{jm}.
\end{align*}
On the other hand, we also have the estimate
\begin{align*}
\left\Vert \left\vert \cdot\right\vert ^{n}K\right\Vert _{L^{\infty}}  &
\leq\sum\limits_{\left\vert \rho\right\vert =n}\left\Vert \left(
\cdot\right)  ^{\rho}K\right\Vert _{L^{\infty}}\leq\sum\limits_{\left\vert
\rho\right\vert =n}\left\Vert \partial^{\rho}\left(  P\tilde{\varphi}%
_{j}\right)  \right\Vert _{L^{1}}\\
&  \leq C\sum\limits_{\left\vert \rho\right\vert =n}2^{j(m-k)}2^{j(n-k)}%
2^{-jn}\leq C2^{jm},
\end{align*}
which gives the claim.\fin

\begin{corollary}
\label{Teo: Operador pseudo Besov X} Let $F\subset\mathcal{S}^{\prime}$ be a
Banach space verifying \textbf{(B3)}, $1\leq r\leq\infty$, and $m,s\in
\mathbb{R}$. Assume that $P\in C^{\infty}\left(  \mathbb{R}^{n}\backslash
\{0\}\right)  $ satisfies $\left\vert \partial^{\rho}P\left(  \xi\right)
\right\vert \leq C\left\vert \xi\right\vert ^{\left(  m-\left\vert
\rho\right\vert \right)  },$ for all $\rho\in\mathbb{N}_{0}^{n}$ and
$\left\vert \rho\right\vert \leq n$. Then, we have the estimate
\[
\left\Vert P\left(  D\right)  f\right\Vert _{\dot{B}\left[  F\right]
_{r}^{s-m}}\leq C\left\Vert f\right\Vert _{\dot{B}\left[  F\right]  _{r}^{s}%
},
\]
for all $f\in\dot{B}\left[  F\right]  _{r}^{s}$, where $C>0$ is a universal constant.
\end{corollary}

\noindent\textbf{Proof.} Note that, for every $j\in\mathbb{Z}$ and $\xi
\in\tilde{D}_{j}$, we have that $\left\vert \xi\right\vert ^{m-\left\vert
\rho\right\vert }\leq C2^{j\left(  m-\left\vert \rho\right\vert \right)  }.$
So, $\left\vert \partial^{\rho}P\left(  \xi\right)  \right\vert \leq
C2^{j\left(  m-\left\vert \beta\right\vert \right)  }.$ On the other hand,
since $\mbox{supp}\widehat{\Delta_{j}f}\subset D_{j}$, we can use Lemma
\ref{Lem: aux Operador pseudo X} to obtain
\[
\left\Vert \Delta_{j}\left(  P\left(  D\right)  f\right)  \right\Vert
_{F}=\left\Vert P\left(  D\right)  \left(  \Delta_{j}f\right)  \right\Vert
_{F}\leq C2^{jm}\left\Vert \Delta_{j}f\right\Vert _{F}.
\]
Multiplying by $2^{j\left(  s-m\right)  }$ and taking the $l^{r}$-norm yield
the result. \fin

The lemma below contains an auxiliary result that relates the sequence spaces
defined in Section \ref{sec: preliminares} to the Besov-type spaces. This
result has already been showed in some particular cases (see \cite{BL},
\cite{Maz}).

\begin{lemma}
\label{Lem: retract}Let $s\in\mathbb{R}$, $1\leq r\leq\infty$, and let $F$ be
a Banach space verifying \textbf{(B1)} and \textbf{(B3)}. Then, the space
$\dot{B}\left[  F\right]  _{r}^{s}$ is a retract of the sequence space
$\dot{l}_{r}^{s}(F)$.
\end{lemma}

\noindent\textbf{Proof. }We define the operator $\mathcal{I}$ as $\left(
\mathcal{I}\left(  f\right)  \right)  _{j}:=\Delta_{j}f$, for all
$f\in\mathcal{S}^{\prime}$. It is clear that $\mathcal{I}:\dot{B}\left[
F\right]  _{r}^{s}\longrightarrow\dot{l}_{r}^{s}(F)$ is bounded. We also
define the operator $\mathcal{K}$ as $\mathcal{K}(\beta)=\sum
\limits_{k=-\infty}^{\infty}\tilde{\Delta}_{j}\beta_{j}$, $\forall
\beta=\left(  \beta_{j}\right)  \in\dot{l}_{r}^{s}(F)$. Since $\varphi
_{j}=\tilde{\varphi}_{j}\varphi_{j}$, it follows that $\Delta_{j}%
\mathcal{K}(\beta)=\Delta_{j}\beta_{j}$, and
\[
\left\Vert \mathcal{K}(\beta)\right\Vert _{\dot{B}\left[  F\right]  _{r}^{s}%
}=\left(  \sum\limits_{j\in\mathbb{Z}}2^{jsr}\left\Vert \Delta_{j}\beta
_{j}\right\Vert _{F}^{r}\right)  ^{\frac{1}{r}}\leq C\left(  \sum
\limits_{j\in\mathbb{Z}}2^{jsr}\left\Vert \beta_{j}\right\Vert _{F}%
^{r}\right)  ^{\frac{1}{r}}\leq C\left\Vert \beta\right\Vert _{\dot{l}_{r}%
^{s}(F)}.
\]
Therefore, $\mathcal{K}:\dot{l}_{r}^{s}(F)\longrightarrow\dot{B}\left[
F\right]  _{r}^{s}$ is bounded. Moreover, $\mathcal{K}\circ\mathcal{I}$ is the
identity in $\dot{B}\left[  F\right]  _{r}^{s}$, as desired. \fin

\begin{corollary}
\label{Cor: Interpolacion entre Besov-X para dar Besov-X.} Let $F$ satisfy the
conditions of Lemma \ref{Lem: retract} and assume that $1\leq r_{0},r_{1}%
\leq\infty$ and $s_{0},s_{1}\in\mathbb{R}$ with $s_{0}\neq s_{1}.$ Then,
\[
\left(  \dot{B}\left[  F\right]  _{r_{0}}^{s_{0}},\dot{B}\left[  F\right]
_{r_{1}}^{s_{1}}\right)  _{\theta,r}=\dot{B}\left[  F\right]  _{r}^{s},\text{
}\forall r\in\lbrack1,\infty],
\]
where $\theta\in(0,1)$ and $s=\left(  1-\theta\right)  s_{0}+\theta s_{1}.$ In
the case $s_{0}=s_{1}=s,$ we have that
\[
\left(  \dot{B}\left[  F\right]  _{r_{0}}^{s},\dot{B}\left[  F\right]
_{r_{1}}^{s}\right)  _{\theta,r}=\dot{B}\left[  F\right]  _{r}^{s}%
\]
provided that $\frac{1}{r}=\frac{1-\theta}{r_{0}}+\frac{\theta}{r_{1}}$.
\end{corollary}

\textbf{Proof.} It follows from Lemma
\ref{lem: Intepolacion en sequencias manteniendo fijo el espacio base.} that
\[
\left(  \dot{l}_{r_{0}}^{s_{0}}\left(  F\right)  ,\dot{l}_{r_{1}}^{s_{1}%
}\left(  F\right)  \right)  _{\theta,r}=\dot{l}_{r}^{s}\left(  F\right)  .
\]
Now we can conclude the proof by employing Lemma \ref{Lem: retract}.\fin

\begin{lemma}
Let $F_{0}$ and $F_{1}$ verify the conditions \textbf{(B1)}-\textbf{(B3)},
$\theta\in\left(  0,1\right)  $ and $1\leq r\leq\infty.$ Then, the
interpolation space $F=\left(  F_{0},F_{1}\right)  _{\theta,r}$ also verifies
\textbf{(B1)}-\textbf{(B3)} with $\sigma_{F}=\left(  1-\theta\right)
\sigma_{F_{0}}+\theta\sigma_{F_{1}}.$
\end{lemma}

\noindent\textbf{Proof.} Let $j\in\mathbb{Z}$ fixed and note that $\Delta
_{j}:F_{i}\longrightarrow L^{\infty}$ with operator norm bounded by
$C2^{j\sigma_{F_{i}}}$ $(i=0,1)$. Then, by interpolation, we obtain that
\[
\left\Vert \Delta_{j}f\right\Vert _{L^{\infty}}=\left\Vert \Delta
_{j}f\right\Vert _{\left(  L^{\infty},L^{\infty}\right)  _{\theta,r}}\leq
C2^{j\left(  \left(  1-\theta\right)  \sigma_{F_{0}}+\theta\sigma_{F_{1}%
}\right)  }\left\Vert f\right\Vert _{\left(  F_{0},F_{1}\right)  _{\theta,r}%
}.
\]
This proves \textbf{(B2)} for $F$. The property \textbf{(B3)} can be proved
in the same way.

In the sequel we show \textbf{(B1)}. The fact $F\subset\mathcal{S}^{\prime}$
follows trivially. On the other hand, for $\lambda>0,$ we have that
\begin{align*}
K\left(  t,f(\lambda\cdot)\right)   &  =\mathop{\inf}\limits_{f(\lambda
\cdot)=g_{0}+g_{1}}\left(  \left\Vert g_{0}\right\Vert _{F_{0}}+t\left\Vert
g_{1}\right\Vert _{F_{1}}\right)  \leq\mathop{\inf}\limits_{f(\lambda
\cdot)=f_{0}(\lambda\cdot)+f_{1}(\lambda\cdot)}\left(  \left\Vert
f_{0}(\lambda\cdot)\right\Vert _{F_{0}}+t\left\Vert f_{1}(\lambda
\cdot)\right\Vert _{F_{1}}\right) \\
&  \lesssim\mathop{\inf}\limits_{f=f_{0}+f_{1}}\left(  \lambda^{\sigma_{F_{0}%
}}\left\Vert f_{0}\right\Vert _{F_{0}}+t\lambda^{\sigma_{F_{1}}}\left\Vert
f_{1}\right\Vert _{F_{1}}\right)  =\lambda^{\sigma_{F_{0}}}%
\mathop{\inf}\limits_{f=f_{0}+f_{1}}\left(  \left\Vert f_{0}\right\Vert
_{F_{0}}+t\lambda^{\sigma_{F_{1}}-\sigma_{F_{0}}}\left\Vert f_{1}\right\Vert
_{F_{1}}\right) \\
&  =\lambda^{\sigma_{F_{0}}}K\left(  \lambda^{\sigma_{F_{1}}-\sigma_{F_{0}}%
}t,f\right)  .
\end{align*}
Thus,
\begin{align*}
\Phi_{\theta,r}\left(  K\left(  t,f(\lambda\cdot)\right)  \right)   &
=\left(  \int_{0}^{\infty}\left(  t^{-\theta}K\left(  t,f(\lambda
\cdot)\right)  \right)  ^{r}\frac{dt}{t}\right)  ^{1/r}\lesssim\left(
\int_{0}^{\infty}\left(  t^{-\theta}\lambda^{\sigma_{F_{0}}}K\left(
\lambda^{\sigma_{F_{1}}-\sigma_{F_{0}}}t,f\right)  \right)  ^{r}\frac{dt}%
{t}\right)  ^{1/r}\\
&  =\lambda^{\sigma_{F_{0}}}\left(  \int_{0}^{\infty}\left(  \left(
\lambda^{\sigma_{F_{0}}-\sigma_{F_{1}}}u\right)  ^{-\theta}K\left(
u,f\right)  \right)  ^{r}\frac{du}{u}\right)  ^{1/r}\\
&  =\lambda^{\left(  1-\theta\right)  \sigma_{F_{0}}+\theta\sigma_{F_{1}}%
}\left(  \int_{0}^{\infty}\left(  u^{-\theta}K\left(  u,f\right)  \right)
^{r}\frac{du}{u}\right)  ^{1/r}\\
&  =\lambda^{\left(  1-\theta\right)  \sigma_{F_{0}}+\theta\sigma_{F_{1}}}%
\Phi_{\theta,r}\left(  K\left(  t,f\right)  \right)  ,
\end{align*}
which yields $\left\Vert f(\lambda\cdot)\right\Vert _{F}\lesssim
\lambda^{\sigma_{F}}\left\Vert f(\lambda\cdot)\right\Vert _{F}$ with
$\sigma_{F}=\left(  1-\theta\right)  \sigma_{F_{0}}+\theta\sigma_{F_{1}},$ and
we are done. \fin

From the previous results, we have another one about interpolation of
Besov-type spaces.

\begin{corollary}
Let $s_{0},s_{1}\in\mathbb{R}$, $1\leq q_{0},q_{1}<\infty,$ $\theta\in(0,1),$
and let $F_{0},F_{1}$ be two Banach spaces verifying the conditions
\textbf{(B1)-(B3)}. Then
\[
\left(  \dot{B}\left[  F_{0}\right]  _{r_{0}}^{s_{0}},\dot{B}\left[
F_{1}\right]  _{r_{1}}^{s_{1}}\right)  _{\theta,r}=\dot{B}\left[  \left(
F_{0},F_{1}\right)  _{\theta,r}\right]  _{r}^{s},
\]
where $s=(1-\theta)s_{0}+\theta s_{1}$ and $\frac{1}{r}=\frac{1-\theta}{r_{0}%
}+\frac{\theta}{r_{1}}.$
\end{corollary}

\textbf{Proof. }It follows from Lemma
\ref{Lem: Intepolacion en sequencias cambiando el espacio base} that
\[
\left(  \dot{l}_{r_{0}}^{s_{0}}\left(  F_{0}\right)  ,\dot{l}_{r_{1}}^{s_{1}%
}\left(  F_{1}\right)  \right)  _{\theta,r}=\dot{l}_{r}^{s}\left(  \left(
F_{0},F_{1}\right)  _{\theta,r}\right)  .
\]
Using again Lemma \ref{Lem: retract}, we are done.\fin

The following theorem provides a characterization for the dual of Besov-type spaces.

\begin{theorem}
\label{Teo: Dualidad espacios de Besov-X} Let $s\in\mathbb{R},$ $1\leq
r<\infty$ and let $F$ be a Banach space verifying the conditions \textbf{(B1)}
and \textbf{(B3)} such that $F^{\prime}$ also verifies \textbf{(B1)} and
\textbf{(B3)}. Then
\[
\left(  \dot{B}\left[  F\right]  _{r}^{s}\right)  ^{\prime}=\dot{B}\left[
F^{\prime}\right]  _{r^{\prime}}^{-s}.
\]

\end{theorem}

\textbf{Proof. }From Lemma \ref{Lem: dualidad en espacios de sequencias}, we
have that%

\[
\left(  \dot{l}_{r}^{s}\left(  F\right)  \right)  ^{\prime}=\left(  \dot
{l}_{r^{\prime}}^{-s}\left(  F^{\prime}\right)  \right)  .
\]
Thus, the result follows by employing Lemma \ref{Lem: retract}. \fin

\begin{proposition}
\label{Prop: Estimativa para el nucleo del calor en Besov-Lorentz-Morrey y besov-Lorentz-Bolcos}
Let $F$ be a Banach space verifying \textbf{(B1)} and \textbf{(B3)}, $1\leq
r\leq\infty$ and $s,\sigma\in\mathbb{R}$. Then, there exists $C>0$ independent
of $t>0$ such that
\begin{equation}
\left\Vert U\left(  t\right)  f\right\Vert _{\dot{B}\left[  F\right]
_{r}^{\sigma}}\leq Ct^{\left(  s-\sigma\right)  /2}\left\Vert f\right\Vert
_{\dot{B}\left[  F\right]  _{r}^{s}},\text{ for all }f\in\dot{B}\left[
F\right]  _{r}^{s}.\text{ }
\label{eq: Estimativa para nucleo do calor s leq sigma BLM}%
\end{equation}
Moreover, if $s<\sigma$, then
\begin{equation}
\left\Vert U\left(  t\right)  f\right\Vert _{\dot{B}\left[  F\right]
_{1}^{\sigma}}\leq Ct^{\left(  s-\sigma\right)  /2}\left\Vert f\right\Vert
_{\dot{B}\left[  F\right]  _{\infty}^{s}},\,\,\text{for all }f\in\dot
{B}\left[  F\right]  _{\infty}^{s}.
\label{eq: Estimativa para nucleo do calor s menor que sigma BLM}%
\end{equation}

\end{proposition}

\noindent\textbf{Proof.} Note that, for every $\beta\in\mathbb{N}_{0}^{n}$,
there exists a polynomial $p_{\rho}(\cdot)$ of degree $\left\vert
\rho\right\vert $ such that
\[
\partial^{\rho}\exp\left(  -t\left\vert \xi\right\vert ^{2}\right)
=t^{\left\vert \rho\right\vert /2}p_{\rho}(\sqrt{t}\xi)\exp\left(
-t\left\vert \xi\right\vert ^{2}\right)  .
\]
Thus, for some $C>0,$ it follows that
\[
\left\vert \partial^{\rho}\exp\left(  -t\left\vert \xi\right\vert ^{2}\right)
\right\vert \leq Ct^{-m/2}\left\vert \xi\right\vert ^{-m-\left\vert
\rho\right\vert }.
\]
Now, using Corollary \ref{Teo: Operador pseudo Besov X}, we obtain
\[
\left\Vert U\left(  t\right)  f\right\Vert _{\dot{B}\left[  F\right]
_{r}^{s+m}}\leq Ct^{-m/2}\left\Vert f\right\Vert _{\dot{B}\left[  F\right]
_{r}^{s}}.
\]
Taking $m=\sigma-s,$ we conclude the estimate \eqref{eq: Estimativa para nucleo do calor s leq sigma BLM}.

Next we turn to prove
\eqref{eq: Estimativa para nucleo do calor s menor que sigma BLM}. For
$s<\sigma,$ from \eqref{eq: Estimativa para nucleo do calor s leq sigma BLM}
with $r=\infty$, we get
\[
\left\Vert U\left(  t\right)  f\right\Vert _{\dot{B}\left[  F\right]
_{\infty}^{2\sigma-s}}\leq Ct^{s-\sigma}\left\Vert f\right\Vert _{\dot
{B}\left[  F\right]  _{\infty}^{s}}%
\]
and
\[
\left\Vert U\left(  t\right)  f\right\Vert _{\dot{B}\left[  F\right]
_{\infty}^{s}}\leq C\left\Vert f\right\Vert _{\dot{B}\left[  F\right]
_{\infty}^{s}}.
\]
Using Corollary \ref{Cor: Interpolacion entre Besov-X para dar Besov-X.}, we
arrive at
\[
U(t):\,\dot{B}\left[  F\right]  _{\infty}^{s}\longrightarrow\left(  \dot
{B}\left[  F\right]  _{\infty}^{2\sigma-s},\dot{B}\left[  F\right]  _{\infty
}^{2\sigma-s}\right)  _{\frac{1}{2},1}=\dot{B}\left[  F\right]  _{1}^{\sigma
},
\]
with $\left\Vert U\left(  t\right)  \right\Vert _{\dot{B}\left[  F\right]
_{\infty}^{s}\longrightarrow\dot{B}\left[  F\right]  _{1}^{\sigma}}\leq
Ct^{\left(  s-\sigma\right)  /2}$, as required.\fin

In order to prove the bilinear estimate in the context of Besov-type spaces we
need a product estimate in such spaces. In the next lemma we establish some
conditions to ensure an appropriate product estimate.

\begin{lemma} \label{Lem: Estimativa producto}Let $s\in\mathbb{R}$
	and $F\in\mathcal{G}$ verify the following conditions:
	\begin{description}
		\item [{(B4)}] $F^{\prime}\in\mathcal{G}$ and verifies \textbf{(B3)}.
		\item [{(B5)}] There exist Banach spaces $G_{0},G_{1},G_{2}\in\mathcal{G}$
		such that
		
		$\left\Vert fg\right\Vert _{F^{\prime}}\le C$$\left\Vert f\right\Vert _{G_{0}}$$\left\Vert g\right\Vert _{G_{0}}$
		with $\sigma_{F^{\prime}}=2\sigma_{G_{0}}$,
		
		$\left\Vert \Delta_{0}f\right\Vert _{G_{0}}\le C\left\Vert f\right\Vert _{F^{\prime}},$
		
		$\left\Vert \Delta_{0}f\right\Vert _{G_{1}}\le C\left\Vert f\right\Vert _{F^{\prime}}$,
		
		$\left\Vert \Delta_{0}f\right\Vert _{F^{\prime}}\le C\left\Vert f\right\Vert _{G_{2}}$,
		
		$\left\Vert fg\right\Vert _{G_{2}}\le C$$\left\Vert f\right\Vert _{G_{1}}$$\left\Vert g\right\Vert _{G_{1}}$
		with $\sigma_{G_{2}}=2\sigma_{G_{1}}$.
		\item [{(B6)}] $\sigma_{G_{1}}-\sigma_{F^{\prime}}-s<0$ and $-\frac{\sigma_{F^{\prime}}}{2}-s>0$.
	\end{description}
	Then, there exists $\rho_{0}>0$ such that for all $\rho\in[0,\rho_{0})$ and for $s_{0}=\sigma_{F^{\prime}}+2s,$ we
	have the product estimate
	\[
	\left\Vert fg\right\Vert _{\dot{B}\left[F^{\prime}\right]_{\infty}^{s_{0}+\rho}}\leq C\left\Vert f\right\Vert _{\dot{B}\left[F^{\prime}\right]_{\infty}^{s}}\left\Vert g\right\Vert _{\dot{B}\left[F^{\prime}\right]_{\infty}^{s+\rho}},
	\]
	where $C>0$ is a universal constant.
	
\end{lemma}

\textbf{Proof. } Let $\rho_{0}$ be such that $-\frac{\sigma_{F^{\prime}}}{2}-s-\rho_{0}=0$. Then, for $\rho\in[0,\rho_{0})$ we have that ~$-\frac{\sigma_{F^{\prime}}}{2}-s-\rho>0$. Now, for every $j\in\mathbb{Z}$, we can write

\noindent
\begin{align}
	\Delta_{j}\left(fg\right) & =\sum\limits _{\left|k-j\right|\le4}\Delta_{j}\left(S_{k-2}f\Delta_{k}g\right)+\sum\limits _{\left|k-j\right|\le4}\Delta_{j}\left(S_{k-2}g\Delta_{k}f\right)+\sum\limits _{k\ge j-2}\Delta_{j}\left(\Delta_{k}f\tilde{\Delta}_{k}g\right)\nonumber \\
	& =I_{1}^{j}+I_{2}^{j}+I_{3}^{j}.\label{eq:aux descomposicion paraproducto}
\end{align}

\noindent We estimate $I_{1}^{j}$ as follows:
\begin{align*}
	\left\Vert I_{1}^{j}\right\Vert _{F^{\prime}} & \leq C\sum\limits _{\left\vert k-j\right\vert \leq4}\left\Vert S_{k-2}f\right\Vert _{G_{0}}\left\Vert \Delta_{k}g\right\Vert _{G_{0}}\\
	& \leq C\sum\limits _{\left\vert k-j\right\vert \leq4}\left(\sum\limits _{m\leq k-2}\left\Vert \Delta_{m}f\right\Vert _{G_{0}}\right)\left\Vert \Delta_{k}g\right\Vert _{G_{0}}\\
	& \leq C\sum\limits _{\left\vert k-j\right\vert \leq4}\left(\sum\limits _{m\leq k-2}2^{m\left(\sigma_{G_{0}}-\sigma_{F^{\prime}}\right)}\left\Vert \Delta_{m}f\right\Vert _{F^{\prime}}\right)2^{k\left(\sigma_{G_{0}}-\sigma_{F^{\prime}}\right)}\left\Vert \Delta_{k}g\right\Vert _{F^{\prime}}\\
	& \leq C\left\Vert f\right\Vert _{\dot{B}\left[F^{\prime}\right]_{\infty}^{s}}\left\Vert g\right\Vert _{\dot{B}\left[F^{\prime}\right]_{\infty}^{s+\rho}}\sum\limits _{\left\vert k-j\right\vert \leq4}\left(\sum\limits _{m\leq k-2}2^{m\left(-\frac{\sigma_{F^{\prime}}}{2}-s\right)}\right)2^{k\left(-\frac{\sigma_{F^{\prime}}}{2}-s-\rho\right)}.
\end{align*}
Thus,
\begin{equation}
	\left\Vert I_{1}^{j}\right\Vert _{F^{\prime}}\leq C\left\Vert f\right\Vert _{\dot{B}\left[F^{\prime}\right]_{\infty}^{s}}\left\Vert g\right\Vert _{\dot{B}\left[F^{\prime}\right]_{\infty}^{s+\rho}}2^{j\left(-s_{0}-\rho\right)}.\label{eq: aux producto I1 BWM}
\end{equation}
For $I_{2}^{j},$ we proceed similarly in order to obtain

\noindent
	\begin{align*}
		\left\Vert I_{2}^{j}\right\Vert _{F^{\prime}} & \leq C\left\Vert f\right\Vert _{\dot{B}\left[F^{\prime}\right]_{\infty}^{s}}\left\Vert g\right\Vert _{\dot{B}\left[F^{\prime}\right]_{\infty}^{s+\rho}}\sum\limits _{\left\vert k-j\right\vert \leq4}\left(\sum\limits _{m\leq k-2}2^{m\left(-\frac{\sigma_{F^{\prime}}}{2}-s-\rho\right)}\right)2^{k\left(-\frac{\sigma_{F^{\prime}}}{2}-s\right)},
	\end{align*}
which leads us to	
\begin{equation}
	\left\Vert I_{2}^{j}\right\Vert _{F^{\prime}}\leq C\left\Vert f\right\Vert _{\dot{B}\left[F^{\prime}\right]_{\infty}^{s}}\left\Vert g\right\Vert _{\dot{B}\left[F^{\prime}\right]_{\infty}^{s+\rho}}2^{j\left(-s_{0}-\rho\right)}.\label{eq: aux producto I2}
\end{equation}
Now we estimate $I_{3}^{j}$:
\begin{align*}
	\left\Vert I_{3}^{j}\right\Vert _{F^{\prime}} & \leq\sum\limits _{k\geq j-2}\left\Vert \Delta_{j}\left(\Delta_{k}f\tilde{\Delta}_{k}g\right)\right\Vert _{F^{\prime}}\\
	& \leq C\sum\limits _{k\geq j-2}2^{j\left(\sigma_{F^{\prime}}-\sigma_{G_{2}}\right)}\left\Vert \Delta_{k}f\tilde{\Delta}_{k}g\right\Vert _{G_{2}}\\
	& \leq C2^{j\left(\sigma_{F^{\prime}}-\sigma_{G_{2}}\right)}\sum\limits _{k\geq j-2}\left\Vert \Delta_{k}f\right\Vert _{G_{1}}\left\Vert \tilde{\Delta}_{k}g\right\Vert _{G_{1}}.
\end{align*}
Moreover,
\[
\left\Vert \Delta_{k}f\right\Vert _{G_{1}}\leq C2^{k\left(\sigma_{G_{1}}-\sigma_{F^{\prime}}\right)}\left\Vert \Delta_{k}f\right\Vert _{F^{\prime}}\leq C2^{k\left(\sigma_{G_{1}}-\sigma_{F^{\prime}}-s\right)}\left\Vert f\right\Vert _{\dot{B}\left[F^{\prime}\right]_{\infty}^{s}}.
\]
Similarly, we have that
\[
\left\Vert \tilde{\Delta}_{k}g\right\Vert _{G_{1}}\leq C2^{k\left(\sigma_{G_{1}}-\sigma_{F^{\prime}}\right)}\left\Vert \tilde{\Delta}_{k}g\right\Vert _{F^{\prime}}\leq C2^{k\left(\sigma_{G_{1}}-\sigma_{F^{\prime}}-s-\rho\right)}\left\Vert g\right\Vert _{\dot{B}\left[F^{\prime}\right]_{\infty}^{s+\rho}}.
\]
Thus,
\begin{align}
	\left\Vert I_{3}^{j}\right\Vert _{F^{\prime}} & \leq C\left\Vert f\right\Vert _{\dot{B}\left[F^{\prime}\right]_{\infty}^{s}}\left\Vert g\right\Vert _{\dot{B}\left[F^{\prime}\right]_{\infty}^{s+\rho}}2^{j\left(\sigma_{F^{\prime}}-\sigma_{G_{2}}\right)}\sum\limits _{k\geq j-2}2^{k2\left(\sigma_{G_{1}}-\sigma_{F^{\prime}}-s-\frac{\rho}{2}\right)}\nonumber \\
	& \leq C\left\Vert f\right\Vert _{\dot{B}\left[F^{\prime}\right]_{\infty}^{s}}\left\Vert g\right\Vert _{\dot{B}\left[F^{\prime}\right]_{\infty}^{s+\rho}}2^{j\left(\sigma_{F^{\prime}}-\sigma_{G_{2}}\right)}2^{j2\left(\sigma_{G_{1}}-\sigma_{F^{\prime}}-s-\frac{\rho}{2}\right)}\nonumber \\
	& \leq C\left\Vert f\right\Vert _{\dot{B}\left[F^{\prime}\right]_{\infty}^{s}}\left\Vert g\right\Vert _{\dot{B}\left[F^{\prime}\right]_{\infty}^{s+\rho}}2^{j2\left(\sigma_{G_{1}}-\frac{\sigma_{G_{2}}}{2}+\frac{\sigma_{F^{\prime}}}{2}-\sigma_{F^{\prime}}-s-\frac{\rho}{2}\right)}\nonumber \\
	& \leq C\left\Vert f\right\Vert _{\dot{B}\left[F^{\prime}\right]_{\infty}^{s}}\left\Vert g\right\Vert _{\dot{B}\left[F^{\prime}\right]_{\infty}^{s+\rho}}2^{j2\left(-\frac{\sigma_{F^{\prime}}}{2}-s-\frac{\rho}{2}\right)}\nonumber \\
	& \leq C\left\Vert f\right\Vert _{\dot{B}\left[F^{\prime}\right]_{\infty}^{s}}\left\Vert g\right\Vert _{\dot{B}\left[F^{\prime}\right]_{\infty}^{s+\rho}}2^{j\left(-s_{0}-\rho\right)}.\label{eq:aux producto I3}
\end{align}
Putting together \eqref{eq: aux producto I1 BWM}-\eqref{eq:aux producto I3},
and using \eqref{eq:aux descomposicion paraproducto}, we conclude
the proof. \fin

Finally, we prove our bilinear estimate in the context of Besov-type spaces.

\begin{theorem}
\label{Teo. Estimativa bilinear espacios tipo besov} Let $s\in\mathbb{R}$ and
$F\in\mathcal{G}$ verify \textbf{(B1)-(B6)}. Then, we obtain
\textbf{(H1)-(H4)} for the spaces $E=\dot{B}\left[  F\right]  _{1}^{-s}$ and
$E_{0}=\dot{B}\left[  F\right]  _{1}^{-s_{0}}$. As a consequence, we get the
bilinear estimate (\ref{eq:Estimativa bilinear}) with $E^{\prime}=\dot
{B}\left[  F^{\prime}\right]  _{\infty}^{-\sigma_{F^{\prime}}-1}$.
\end{theorem}

\noindent\textbf{Proof. }The inclusion $\dot{B}\left[  F\right]  _{r}%
^{-s}\subset\mathcal{S}^{\prime}$ follows by definition; moreover, a direct
calculation shows that $\left\Vert f\left(  \lambda\cdot\right)  \right\Vert
_{\dot{B}\left[  F\right]  _{r}^{-s}}\approx\lambda^{-s+\sigma_{F}}\left\Vert
f\right\Vert _{\dot{B}\left[  F\right]  _{r}^{s}}$, and then we obtain
\textbf{(H1)}. The condition \textbf{(H2)} follows directly by using Corollary
\ref{Teo: Operador pseudo Besov X}, and the condition \textbf{(H3)} follows by
Lemma \ref{Lem: Estimativa producto}. Finally, note that Proposition
\ref{Prop: Estimativa para el nucleo del calor en Besov-Lorentz-Morrey y besov-Lorentz-Bolcos}
gives the estimate
\[
t^{\left(  -s_{0}+\sigma_{F}-\left(  -s+\sigma_{F}\right)  \right)
/2-1/2}\left\Vert \nabla U\left(  t\right)  f\right\Vert _{\dot{B}\left[
F\right]  _{1}^{-s_{0}}}\leq t^{\left(  -s_{i}+\sigma_{F}-\left(
-s+\sigma_{F}\right)  \right)  /2-1}\left\Vert f\right\Vert _{\dot{B}\left[
F\right]  _{1}^{-s_{i}}};
\]
in other words, we have that
\[
t^{\left(  \sigma_{E_{0}}-\sigma_{E}\right)  /2-1/2}\left\Vert \nabla U\left(
t\right)  f\right\Vert _{E_{0}}\leq t^{\left(  \sigma_{E_{i}}-\sigma
_{E}\right)  /2-1}\left\Vert f\right\Vert _{E_{i}},
\]
with $E_{i}=\dot{B}\left[  F\right]  _{1}^{-s_{i}},$ $i=1,2$. Taking
$-s_{1}=-s+\epsilon$ and $-s_{2}=-s-\epsilon$ for $0<\epsilon\ll1$, and
defining $\frac{1}{z_{i}}=\frac{1}{2}\left(  \sigma_{E}-\sigma_{E_{i}}\right)
+1,$ we obtain the relations\noindent%
\begin{align*}
\frac{1}{z_{1}}  &  =\frac{1}{2}\left(  -s+\sigma_{F}-(-s+\epsilon+\sigma
_{F})\right)  +1=\frac{1}{2}\epsilon+1>1,\\
0  &  <\frac{1}{z_{2}}=\frac{1}{2}\left(  -s+\sigma_{F}-(-s-\epsilon
+\sigma_{F})\right)  +1=\frac{1}{2}\left(  -\epsilon\right)  +1<1.
\end{align*}
Therefore, $0<z_{1}<1<z_{2}<\infty$ and, for $\theta=1/2$, we have
$1=\frac{\theta}{z_{1}}+\frac{1-\theta}{z_{2}}$. Finally, since $-s=(1-\theta
)\left(  -s_{1}\right)  +\theta\left(  -s_{2}\right)  $, by using Corollary
\ref{Cor: Interpolacion entre Besov-X para dar Besov-X.}, we get
\[
E\hookrightarrow\left(  E_{1},E_{2}\right)  _{\theta,1}=E.
\]

The previous calculations show the conditions \textbf{(H4$^\prime$)-(H4$^{\prime \prime \prime}$ )}. Condition
\textbf{(H4)} follows by Lemma \ref{Lem: condiciones para estimativa integral}%
. Finally, the bilinear estimate (\ref{eq:Estimativa bilinear}) follows from
Theorems \ref{Teo: Teo main. Estimativa integral en el dual.} and
\ref{Teo: Dualidad espacios de Besov-X}.\fin

\subsection{Applications}

\label{Sec-Besov-weak-Herz}

\subsubsection{Bilinear estimate in Besov-weak-Herz spaces}

Using the previous results, we can obtain the bilinear estimate
(\ref{eq:Estimativa bilinear}) for the Besov-weak-Herz space $E^{\prime}%
=\dot{B}\left[  W\dot{K}_{p,\infty}^{0}\right]  _{\infty}^{\frac{n}{p}-1},$
where $\frac{n}{2}<p<n$. For that, consider $F$ as the Lorentz-Herz space
$\dot{K}_{\left(  p^{\prime},1\right)  ,1}^{0}$. The properties
\textbf{(B1)-(B4)} follow by \textbf{(a)}, \textbf{(f)}, \textbf{(g)} and
\textbf{(b)} in page \pageref{Propiedad a) espacios de Lorentz Herz}.
Moreover, taking $G_{0}=W\dot{K}_{2p,\infty}^{0}$, $G_{1}=W\dot{K}_{\tilde
{p},\infty}^{0}$, and $G_{2}=W\dot{K}_{\tilde{p}/2,\infty}^{0}$ with
$p<\tilde{p}<n$, and using the properties \textbf{(d)} and \textbf{(e)} in
page \pageref{Propiedad a) espacios de Lorentz Herz}, we obtain \textbf{(B5)}
and \textbf{(B6)}. Then, the bilinear estimate follows by Theorem
\ref{Teo. Estimativa bilinear espacios tipo besov}. As far as we know, this
estimate in $\dot{B}\left[  W\dot{K}_{p,\infty}^{0}\right]  _{\infty}%
^{\frac{n}{p}-1}$ has not been previously obtained in the literature, at least
we were unable to locate it. Furthermore, as pointed out in Introduction,
using heat-semigroup estimates in Besov-Lorentz-Herz spaces (see
\cite{LucJhe1}), estimate (\ref{eq:Estimativa bilinear}) in $\dot{B}\left[
W\dot{K}_{p,\infty}^{0}\right]  _{\infty}^{\frac{n}{p}-1}$, and an argument by
\cite{Meyer}, one can obtain the uniqueness of mild solutions in $C(\left[
0,T\right)  ;\tilde{Z})$ with $Z=\dot{B}\left[  W\dot{K}_{p,\infty}%
^{0}\right]  _{\infty}^{\frac{n}{p}-1}.$

\subsubsection{Bilinear estimate in Besov and Besov-weak-$L^{p}$ spaces}

In this subsection we reobtain the bilinear estimate in the Besov space
$\dot{B}_{p,\infty}^{\frac{n}{p}-1}$ found in \cite{Cann2}. Also, we obtain an
estimate in the Besov-weak-$L^{p}$ space $\dot{B}\left[  L^{p,\infty}\right]
_{\infty}^{\frac{n}{p}-1}$. In fact, for the Besov space $\dot{B}_{p,\infty}%
^{\frac{n}{p}-1}$ it is sufficient to consider $F=L^{p^{\prime}}$, $G_{0}=L^{2p}%
$, $G_{1}=L^{\tilde{p}},$ $G_{2}=L^{\frac{\tilde{p}}{2}}$ and $E^{\prime
}=B_{p,\infty}^{\frac{n}{p}-1}=\dot{B}\left[  F^{\prime}\right]  _{\infty
}^{-\sigma_{F^{\prime}}-1}$, where $\frac{n}{2}<p<\tilde{p}<n,$ $F^{\prime
}=L^{p}$ and $\sigma_{F^{\prime}}=-\frac{n}{p}.$ For the Besov-weak-$L^{p}$ space
$\dot{B}[L^{p,\infty}]_{\infty}^{\frac{n}{p}-1}$ is sufficient to consider
$F=L^{p^{\prime},1}$, $G_{0}=L^{2p,\infty}$, $G_{1}=L^{\tilde{p},\infty},$
$G_{2}=L^{\frac{\tilde{p}}{2},\infty}$ and $E^{\prime}=\dot{B}[L^{p,\infty}%
]_{\infty}^{\frac{n}{p}-1}=\dot{B}\left[  F^{\prime}\right]  _{\infty
}^{-\sigma_{F^{\prime}}-1}$, where $\frac{n}{2}<p<\tilde{p}<n,$ $F^{\prime
}=L^{p,\infty}$ and $\sigma_{F^{\prime}}=-\frac{n}{p}.$ With these
considerations, we can show \textbf{(B1)-(B6)} and then employ Theorem
\ref{Teo. Estimativa bilinear espacios tipo besov} in order to obtain the
bilinear estimate (\ref{eq:Estimativa bilinear}) with $E^{\prime}=\dot
{B}_{p,\infty}^{\frac{n}{p}-1}$ and $E^{\prime}=\dot{B}\left[  L^{p,\infty
}\right]  _{\infty}^{\frac{n}{p}-1}$, according to the corresponding case. The
details of checking for the properties \textbf{(B1)-(B6)} can be carry out as
in the previous case of Besov-weak-Herz spaces and are left to the reader.

\section{Bilinear estimate in $X$-Morrey spaces}

\label{sect:bilinear-in-morrey-spaces}

This section is devoted to obtaining bilinear estimates in the framework of
abstract Morrey spaces $\mathcal{M}\left[  X\right]  ^{l}$. For our approach,
we need to introduce and develop some core properties in block spaces based on
Banach spaces, which are the preduals of $\mathcal{M}\left[  X\right]  ^{l}.$
For that, we need to assume some natural properties on the base space $X$. A
set of equivalent conditions have already been considered in \cite[Appendix
A]{Lemarie2}. In fact, mainly in view of the block spaces, we assume a little
more than \cite{Lemarie2}, namely the property of passing the $X$-norm into
the integral.

From now on, we denote by $D(a,r)$ the ball in $\mathbb{R}^{n}$ with center
$a\in\mathbb{R}^{n}$ and radius $r>0$. Also, for a set $A$, the symbol $1_{A}$
stands for the characteristic function of $A$.

Let us introduce $\mathcal{H}\subset\mathcal{G}$ as being the class of all
Banach spaces $X$ of measurable functions such that

\begin{description}
\item[\textbf{(M1)}] $X\hookrightarrow L_{loc}^{1}$;

\item[\textbf{(M2)}] If $f\in X$ and $g\in L^{\infty}$, then $fg\in X$ and
$\left\Vert fg\right\Vert _{X}\leq C\left\Vert g\right\Vert _{L^{\infty}%
}\left\Vert f\right\Vert _{X}$;

\item[\textbf{(M3)}] If $f\in X$ and $x_{0}\in\mathbb{R}^{n}$, then
$\tau_{x_{0}}f\left(  x\right)  =f\left(  x-x_{0}\right)  \in X$ and
$\left\Vert \tau_{x_{0}}f\right\Vert _{X}\leq C\left\Vert f\right\Vert _{X}$;

\item[\textbf{(M4)}] For $g\in L_{loc}^{1}$ and $K:\mathbb{R}^{n}%
\times\mathbb{R}^{n}\rightarrow\mathbb{R}$ such that $K(\cdot,y)\in X$, for
every $y\in\mathbb{R}^{n}$, it follows that $\left\Vert \int
\limits_{\mathbb{R}^{n}}g\left(  y\right)  K(\cdot,y)dy\right\Vert _{X}\leq
C\int\limits_{\mathbb{R}^{n}}\left\vert g\left(  y\right)  \right\vert
\left\Vert K(\cdot,y)\right\Vert _{X}dy;$

\item[\textbf{(M5)}] Denoting by $D\left(  0,1\right)  $ the unit ball in
$\mathbb{R}^{n}$, we have that $1_{D\left(  0,1\right)  }\in X;$

\item[\textbf{(M6)}] If $f\in L_{loc}^{1}$ is such that $f1_{D(x_{0},R)}\in
X,$ for all $x_{0}\in\mathbb{R}^{n}$ and $R>0,$ and
$\mathop{\sup}\limits_{x_{0}\in\mathbb{R}^{n},R>0}\left\Vert 1_{D\left(
x_{0},R\right)  }f\right\Vert _{X}<\infty$, then $f\in X$ and
\[
\Vert f\Vert_{X}\leq C\mathop{\sup}\limits_{x_{0}\in\mathbb{R}^{n}%
,R>0}\left\Vert 1_{D\left(  x_{0},R\right)  }f\right\Vert _{X}.
\]

\end{description}

\bigskip

\bigskip Moreover, we say that $X\in\tilde{\mathcal{H}}$ if $X\in\mathcal{H}$
and the following condition \textbf{(M7)} is verified:

\begin{description}
\item[\textbf{(M7)}] The Riesz transforms are bounded on $X$.
\end{description}

We list some consequences of conditions \textbf{(M1})-\textbf{(M6)} in the
following result. Related properties can be found in \cite[Appendix
A]{Lemarie2}.

\begin{proposition}
\label{prop:conseq-condition-morrey-space} Let $X\in\mathcal{H}$, then the
following assertions hold true:

\begin{itemize}
\item[\textbf{(i)}] For all $x_{0}\in\mathbb{R}^{n}$ and $R>0$, we have
$1_{D(x_{0},R)}\in X$ and
\[
\Vert1_{D(x_{0},R)}\Vert_{X}\leq CR^{-\sigma_{X}}\Vert1_{D(0,1)}\Vert_{X}.
\]
In particular, it follows by \textbf{(M2)} that $1_{A}\in X,$ for all bounded
measurable set $A\subset\mathbb{R}^{n}$.

\item[\textbf{(ii)}] If $f\in X$, then
\begin{equation}
\mathop{\sup}\limits_{x_{0}\in\mathbb{R}^{n},R>0}\left\Vert 1_{D\left(
x_{0},R\right)  }f\right\Vert _{X}\leq C\Vert f\Vert_{X}.
\label{aux-desig-equiv-1000}%
\end{equation}

\item[\textbf{(iii)}] For all $x_{0}\in\mathbb{R}^{n}$ and $R>0$, we have
\begin{equation}
\int_{D(x_{0},R)}|f(x)|dx\leq C\Vert f\Vert_{X}R^{n+\sigma_{X}}.
\label{eq: Desigualdad auxiliar}%
\end{equation}
Moreover, $\sigma_{X}=-n$ implies that $X=L^{1}$, provided that $\Vert
\cdot\Vert_{X}$ is lower semicontinuous with respect to a.e. convergence, and
$\sigma_{X}=0$ implies $X=L^{\infty}$.
\end{itemize}
\end{proposition}

\textbf{Proof.} Item \textbf{\textit{(i)}} follows directly from the identity
$1_{D(x_{0},R)}=\tau_{x_{0}}\left(  (1_{D(0,1)})_{R^{-1}}\right)  $,
$X\in\mathcal{G}$, and conditions \textbf{(M3)},\textbf{(M5)}. Since
$1_{D(x_{0},R)}\in L^{\infty}$, condition \textbf{(M2)} implies
\textbf{\textit{(ii)}}. For item \textbf{\textit{(iii)}}, note that the
continuous inclusion in \textbf{(M1)} implies that
\[
\int\limits_{D\left(  0,1\right)  }\left\vert f\left(  x\right)  \right\vert
dx\leq C\left\Vert f\right\Vert _{X}.
\]
Thus, changing variables and using \textbf{(M3)}, we get
(\ref{eq: Desigualdad auxiliar}). Now, assume that $\sigma_{X}=-n$, it is
clear from \eqref{eq: Desigualdad auxiliar} that $X\hookrightarrow L^{1}$.
Moreover, we get from \textbf{(M4)} that
\[
\left\Vert \int_{\mathbb{R}^{n}}f(y)\frac{1}{R^{n}}1_{D(0,R)}%
(x-y)\ dy\right\Vert _{X}\leq\Vert f\Vert_{L^{1}},
\]
for all $f\in L^{1}$. By taking the limit when $R\rightarrow0$ and using the
lower semicontinuity on the norm, we get that $X=L^{1}$ with equivalent norms.
Now, if $\sigma_{X}=0$, by item \textbf{\textit{(iii)}} and the Lebesgue
Differentiation Theorem, it follows that $X\hookrightarrow L^{\infty}$. On the
other hand, given $f\in L^{\infty}$, \textbf{(M2)} implies that $f1_{D(x_{0}%
,R)}\in X$, for every $D(x_{0},R)$, and
\[
\Vert f1_{D(x_{0},R)}\Vert_{X}\leq C\Vert f\Vert_{L^{\infty}}\Vert
1_{D(x_{0},R)}\Vert_{X}\leq C\Vert f\Vert_{L^{\infty}},
\]
where above we have used item \textbf{\textit{(i)}}. Then, \textbf{(M6)}
implies that $f\in X$ with $\Vert f\Vert_{X}\leq C\Vert f\Vert_{L^{\infty}}$.\fin

\begin{remark}
\label{remark:restriction-on-scaling} Note that, in order to work with
nontrivial spaces, estimate \eqref{eq: Desigualdad auxiliar} implies that
$-n\leq\sigma_{X}\leq0,$ for all $X\in\mathcal{H}$. Also, item
\textbf{\textit{(ii)}} in Proposition \ref{prop:conseq-condition-morrey-space}
and \textbf{(M6)} imply $\Vert f\Vert_{X}\approx C\mathop{\sup}\limits_{x_{0}%
\in\mathbb{R}^{n},R>0}\left\Vert 1_{D\left(  x_{0},R\right)  }f\right\Vert
_{X}$ which was assumed by \cite{Lemarie2} as a hypothesis on $X$.
\end{remark}

Next, we define Morrey spaces and their preduals based on abstract Banach
spaces $X$. The aim is to ensure the requirements for applying the results
developed in Section \ref{sect:bilinear-estimate-general-spaces}. The
definition of generalized Morrey spaces in item \textbf{\textit{(ii)}} below
can be found in \cite{Lemarie2}.

\begin{definition}
Let $X\in\mathcal{H}$ and $1\leq l\leq\infty$. \label{def:X-block-X-morrey}

\begin{itemize}
\item[\textbf{(i)}] Assume that $\frac{n}{l}+\sigma_{X}\geq0$. We say that a
mensurable function $b$ is a $\left(  l,X\right)  $-block if there exist
$a\in\mathbb{R}^{n}$ and $\rho>0$ such that $\mbox{supp}\left(  b\right)
\subset D\left(  a,\rho\right)  $ and
\[
\rho^{\frac{n}{l}+\sigma_{X}}\left\Vert 1_{D\left(  a,\rho\right)
}b\right\Vert _{X}\leq1.
\]

We define the $(l,X)$-block space, denoted by $\mathcal{PD}\left[  X\right]
^{l}$, as
\begin{align*}
&  \mathcal{PD}\left[  X\right]  ^{l}=\\
&  \left\{  h\in\mathcal{S}^{\prime};\exists\{\alpha_{k}\}_{k\in\mathbb{N}%
}\subset\mathbb{R},\sum\limits_{k=1}^{\infty}\left\vert \alpha_{k}\right\vert
<\infty\mbox{ and }h=\sum\limits_{k=1}^{\infty}\alpha_{k}b_{k}%
,\mbox{ a.e. in }\mathbb{R}^{n},\mbox{ where }b_{k}\mbox{ is a }\left(
l,X\right)  \text{-}\mbox{block}\right\}  .
\end{align*}

\item[\textbf{(ii)}] Assume that $\frac{n}{l}+\sigma_{X}\leq0$. Define the
$X$-Morrey space $\mathcal{M}\left[  X\right]  ^{l}$ as
\[
\mathcal{M}\left[  X\right]  ^{l}=\left\{  g\in\mathcal{S}^{\prime
};\mathop{\sup}\limits_{x\in\mathbb{R}^{n},\rho>0}\rho^{\frac{n}{l}+\sigma
_{X}}\left\Vert 1_{D\left(  x,\rho\right)  }g\right\Vert _{X}<\infty\right\}
.
\]

\end{itemize}
\end{definition}

Following the ideas as the classic case $X=L^{p}(\mathbb{R}^{n})$, it is
possible to show that the space $\mathcal{PD}[X]^{l}$ is Banach with the norm
$\Vert h\Vert_{\mathcal{PD}[X]^{l}}:=\inf\left\{  {\displaystyle\sum
_{k}|\alpha_{k}|}\right\}  ,$ where the infimum is taken over the sequences
$(\alpha_{k})_{k\in\mathbb{N}}$ such that $h=\sum_{k=1}^{\infty}\alpha
_{k}b_{k}$. In the same way, it follows that the space $\mathcal{M}[X]^{l}$ is
Banach with the norm $\Vert g\Vert_{\mathcal{M}[X]^{l}}%
=\mathop{\sup}\limits_{x\in\mathbb{R}^{n},\rho>0}\rho^{\frac{n}{l}+\sigma_{X}%
}\left\Vert 1_{D\left(  x,\rho\right)  }g\right\Vert _{X}$.

\begin{lemma}
\label{Lem: espacio de bloques pertenece a la clase G} Let $X\in\mathcal{H}$,
then $\mathcal{PD}\left[  X\right]  ^{l}\in\mathcal{G}$.
\end{lemma}

\textbf{Proof.} It follows directly from the definition that $\mathcal{PD}%
\left[  X\right]  ^{l}\subset\mathcal{S}^{\prime}$. Recall that, by the
definition of $\mathcal{H}$, we have that $X\in\mathcal{G}$. Let $\lambda>0$
fixed, if $b$ is a $\left(  l,X\right)  $-block, then $\frac{b_{\lambda}%
}{C\lambda^{-\frac{n}{l}}}$ is a $\left(  l,X\right)  $-block, where $C>0$ is
the constant given at the definition of the family $\mathcal{G}$ (see
Definition \ref{def:G-set}). In fact, if $b$ is $\left(  l,X\right)  $-block,
then there exist $a\in\mathbb{R}^{n}$ and $\rho>0$ such that
\[
\rho^{\frac{n}{l}+\sigma_{X}}\left\Vert 1_{D\left(  a,\rho\right)
}b\right\Vert _{X}\leq1.
\]

Using the identity $1_{D\left(  x,r\right)  }b_{\lambda}=\left(  1_{D\left(
x,r\right)  }\left(  \lambda^{-1}\cdot\right)  b\right)  _{\lambda}=\left(
1_{D\left(  \lambda x,\lambda r\right)  }b\right)  _{\lambda}$, we arrive at
\[
\left\Vert 1_{D\left(  x,r\right)  }b_{\lambda}\right\Vert _{X}\leq
C\lambda^{\sigma_{X}}\left\Vert 1_{D\left(  \lambda x,\lambda r\right)
}b\right\Vert _{X}.
\]
Now, defining $\tilde{a}=a/\lambda$ and $\tilde{\rho}=\rho/\lambda$, we have
that
\[
\tilde{\rho}^{\frac{n}{l}+\sigma_{X}}\left\Vert 1_{D\left(  \tilde{a}%
,\tilde{\rho}\right)  }\frac{b_{\lambda}}{C\lambda^{-\frac{n}{l}}}\right\Vert
_{X}\leq C\tilde{\rho}^{\frac{n}{l}+\sigma_{X}}\lambda^{\sigma_{X}}%
\frac{\lambda^{\frac{n}{l}}}{C}\left\Vert 1_{D\left(  \lambda\tilde{a}%
,\lambda\tilde{\rho}\right)  }b\right\Vert _{X}=\rho^{\frac{n}{l}+\sigma_{X}%
}\left\Vert 1_{D\left(  a,\rho\right)  }b\right\Vert _{X}\leq1.
\]
Thus, if $h={\displaystyle\sum_{k=1}^{\infty}\alpha_{k}b_{k}}$, then
$h_{\lambda}={\displaystyle\sum_{k=1}^{\infty}\left(  C\lambda^{-\frac{n}{l}%
}\alpha_{k}\right)  \frac{(b_{k})_{\lambda}}{C\lambda^{-\frac{n}{l}}}}$, which
gives $\Vert h_{\lambda}\Vert_{\mathcal{PD}[X]^{l}}\leq C\lambda^{-\frac{n}%
{l}}\Vert h\Vert_{\mathcal{PD}[X]^{l}}$, as desired. \fin

In view of the previous proof, we can see that $\sigma_{\mathcal{PD}[X]^{l}%
}=-\frac{n}{l}$. As expected, we recover the duality result between
$(l,X)$-block spaces and their counterpart $(l^{\prime},X^{\prime})$-Morrey spaces.

\begin{theorem}
\label{Teo: Dualidad bloques-Morrey} Let $1\leq l\leq\infty$ and
$X\in\mathcal{H}$ be such that $X^{\prime}\in\mathcal{H}$ and $\frac{n}%
{l}+\sigma_{X}\geq0$. Then,
\begin{equation}
\left(  \mathcal{PD}\left[  X\right]  ^{l}\right)  ^{\prime}=\mathcal{M}%
\left[  X^{\prime}\right]  ^{l^{\prime}}, \label{aux-duality-Morrey-1}%
\end{equation}
where $l^{\prime}$ is the conjugate exponent of $l$.
\end{theorem}

\textbf{Proof. }Let $g\in\mathcal{M}\left[  X^{\prime}\right]  ^{l^{\prime}}$,
for $f\in\mathcal{PD}\left[  X\right]  ^{l}$ we have that
\begin{align*}
\left\vert \left\langle g,f\right\rangle \right\vert  &  \leq\sum
\limits_{k=1}^{\infty}\left\vert \alpha_{k}\right\vert |\left\langle
g,b_{k}\right\rangle |=\sum\limits_{k=1}^{\infty}\left\vert \alpha
_{k}\right\vert |\left\langle g,1_{D\left(  a_{k},\rho_{k}\right)  }%
b_{k}\right\rangle |\\
&  \leq\sum\limits_{k=1}^{\infty}\left\vert \alpha_{k}\right\vert \left\Vert
1_{D\left(  a_{k},\rho_{k}\right)  }g\right\Vert _{X^{\prime}}\left\Vert
1_{D\left(  a_{k},\rho_{k}\right)  }b_{k}\right\Vert _{X}\\
&  =\sum\limits_{k=1}^{\infty}\left\vert \alpha_{k}\right\vert \rho
_{k}^{-\left(  \frac{n}{l}+\sigma_{X}\right)  }\left\Vert 1_{D\left(
a_{k},\rho_{k}\right)  }g\right\Vert _{X^{\prime}}\rho_{k}^{\frac{n}{l}%
+\sigma_{X}}\left\Vert 1_{D\left(  a_{k},\rho_{k}\right)  }b_{k}\right\Vert
_{X}\\
&  \leq\sum\limits_{k=1}^{\infty}\left\vert \alpha_{k}\right\vert \rho
_{k}^{\frac{n}{l^{\prime}}+\sigma_{X^{\prime}}}\left\Vert 1_{D\left(
a_{k},\rho_{k}\right)  }g\right\Vert _{X^{\prime}}\leq\left\Vert g\right\Vert
_{\mathcal{M}\left[  X^{\prime}\right]  ^{l^{\prime}}}\sum\limits_{k=1}%
^{\infty}\left\vert \alpha_{k}\right\vert .
\end{align*}
Taking the infimum in the right-hand side of the above inequality, we obtain
\[
\left\vert \left\langle g,f\right\rangle \right\vert \leq\left\Vert
g\right\Vert _{\mathcal{M}\left[  X^{\prime}\right]  ^{l^{\prime}}}\left\Vert
f\right\Vert _{\mathcal{PD}\left[  X\right]  ^{l}},
\]
which implies that
\[
\mathcal{M}\left[  X^{\prime}\right]  ^{l^{\prime}}\hookrightarrow\left(
\mathcal{PD}\left[  X\right]  ^{l}\right)  ^{\prime}.
\]
Conversely, let $\Psi\in\left(  \mathcal{PD}\left[  X\right]  ^{l}\right)
^{\prime}$. Note that, for every $\rho>0$ and $f\in X,$ we have that
$\rho^{-\left(  \frac{n}{l}+\sigma_{X}\right)  }1_{D\left(  0,\rho\right)
}f\in\mathcal{PD}\left[  X\right]  ^{l}$ with $\left\Vert \rho^{-\left(
\frac{n}{l}+\sigma_{X}\right)  }1_{D\left(  0,\rho\right)  }f\right\Vert
_{\mathcal{PD}\left[  X\right]  ^{l}}\leq\left\Vert f\right\Vert _{X}$.
Defining
\[
T_{\rho}\left(  f\right)  =\Psi\left(  \rho^{-\left(  \frac{n}{l}+\sigma
_{X}\right)  }1_{D\left(  0,\rho\right)  }f\right)  ,
\]
we get that $T_{\rho}\in X^{\prime}\hookrightarrow L_{loc}^{1}$, and then
there exists $g_{\rho}\in X^{\prime}$ such that $T_{\rho}\left(  f\right)
=\left\langle g_{\rho},f\right\rangle $. Let $\left(  \rho_{k}\right)
_{k\in\mathbb{N}}$ be such that $\rho_{k}\rightarrow\infty$ and set $g\left(
x\right)  =g_{\rho_{k}}\left(  x\right)  $ if $x\in D(0,\rho_{k})$. We intend
to show that $g\in\mathcal{M}\left[  X^{\prime}\right]  ^{l^{\prime}}$. For
this, let $a\in\mathbb{R}^{n}$ and $r>0$, and take $k$ such that $D\left(
a,r\right)  \subset D\left(  0,\rho_{k}\right)  $. Thus, we have that
\begin{align*}
\rho^{\frac{n}{l^{\prime}}+\sigma_{X^{\prime}}}\left\Vert 1_{D\left(
a,r\right)  }g\right\Vert _{X^{\prime}}  &  =\rho^{\frac{n}{l^{\prime}}%
+\sigma_{X^{\prime}}}\mathop{\sup}\limits_{\left\Vert f\right\Vert _{X}%
=1}\left\vert \left\langle 1_{D\left(  a,r\right)  }g,f\right\rangle
\right\vert =\rho^{\frac{n}{l^{\prime}}+\sigma_{X^{\prime}}}%
\mathop{\sup}\limits_{\left\Vert f\right\Vert _{X}=1}\left\vert \left\langle
g_{\rho_{k}},1_{D\left(  a,r\right)  }f\right\rangle \right\vert \\
&  =\rho^{\frac{n}{l^{\prime}}+\sigma_{X^{\prime}}}%
\mathop{\sup}\limits_{\left\Vert f\right\Vert _{X}=1}\left\vert \left\langle
\Psi,1_{D\left(  a,r\right)  }f\right\rangle \right\vert \leq C\left\Vert
\Psi\right\Vert \mathop{\sup}\limits_{\left\Vert f\right\Vert _{X}%
=1}\left\Vert \rho^{-\left(  \frac{n}{l}+\sigma_{X}\right)  }1_{D\left(
a,r\right)  }f\right\Vert _{\mathcal{PD}\left[  X\right]  ^{l}}\\
&  \leq C\left\Vert \Psi\right\Vert \mathop{\sup}\limits_{\left\Vert
f\right\Vert _{X}=1}\left\Vert f\right\Vert _{X}\leq C\left\Vert
\Psi\right\Vert .
\end{align*}
From the above, we get that $\left\Vert g\right\Vert _{\mathcal{M}\left[
X^{\prime}\right]  ^{l^{\prime}}}\leq C\left\Vert \Psi\right\Vert ,$ and then
$\left(  \mathcal{PD}\left[  X\right]  ^{l}\right)  ^{\prime}\hookrightarrow
\mathcal{M}\left[  X^{\prime}\right]  ^{l^{\prime}}$. \fin

It is possible to show a H\"{o}lder-type inequality in Morrey spaces by
assuming that the base spaces also satisfy a H\"{o}lder-type inequality.

\begin{lemma}
\label{Lem: producto en morrey} Let $G,G_{0},G_{1}\in\mathcal{H}$ be such that
$f\in G_{0}$ and $g\in G_{1}$ imply $f\cdot g\in G$ with
\begin{equation}
\left\Vert f\cdot g\right\Vert _{G}\leq C\left\Vert f\right\Vert _{G_{0}%
}\left\Vert g\right\Vert _{G_{1}}, \label{aux-hip-Holder-1}%
\end{equation}
where the constant $C>0$ is independent of $f,g$. Assume also that $\frac
{1}{l}=\frac{1}{l_{0}}+\frac{1}{l_{1}}$. Then, for every $f\in\mathcal{M}%
[G_{0}]^{l_{0}}$ and $g\in\mathcal{M}[G_{1}]^{l_{1}}$, we have that $f\cdot
g\in\mathcal{M}[G]^{l}$ with the H\"{o}lder-type inequality
\[
\left\Vert fg\right\Vert _{\mathcal{M}\left[  G\right]  ^{l}}\leq C\left\Vert
f\right\Vert _{\mathcal{M}\left[  G_{0}\right]  ^{l_{0}}}\left\Vert
g\right\Vert _{\mathcal{M}\left[  G_{1}\right]  ^{l_{1}}},
\]
where $C$ is the same constant of (\ref{aux-hip-Holder-1}).
\end{lemma}

\textbf{Proof. } First, by Remark
\ref{remark:scaling-condit-to-holder-type-inequa}, recall that $\sigma
_{G}=\sigma_{G_{0}}+\sigma_{G_{1}}$. Since $1_{D\left(  x_{0},\rho_{0}\right)
}fg=1_{D\left(  x_{0},\rho_{0}\right)  }f1_{D\left(  x_{0},\rho_{0}\right)
}g,$ we have that
\[
\left\Vert 1_{D\left(  x_{0},\rho_{0}\right)  }fg\right\Vert _{G}\leq
C\left\Vert 1_{D\left(  x_{0},\rho_{0}\right)  }f\right\Vert _{G_{0}%
}\left\Vert 1_{D\left(  x_{0},\rho_{0}\right)  }g\right\Vert _{G_{1}},
\]
which leads us to%
\begin{align*}
\rho_{0}^{\frac{n}{l}+\sigma_{G}}\left\Vert 1_{D\left(  x_{0},\rho_{0}\right)
}fg\right\Vert _{G}  &  \leq C\rho_{0}^{\frac{n}{l_{0}}+\sigma_{G_{0}}%
}\left\Vert 1_{D\left(  x_{0},\rho_{0}\right)  }f\right\Vert _{G_{0}}\rho
_{0}^{\frac{n}{l_{1}}+\sigma_{G_{1}}}\left\Vert 1_{D\left(  x_{0},\rho
_{0}\right)  }g\right\Vert _{G_{1}}\\
&  \leq C\left\Vert f\right\Vert _{\mathcal{M}\left[  G_{0}\right]  ^{l_{0}}%
}\left\Vert g\right\Vert _{\mathcal{M}\left[  G_{1}\right]  ^{l_{1}}.}%
\end{align*}
The result follows by taking the supremum in the left-hand side of the above inequality.\fin

The result below is due to Lemarie-Rieusset \cite[Proposition B.1.]{Lemarie2}.

\begin{lemma}
\label{Lem: Limitacion de operador de Riesz en Morrey} Let $G\in
\tilde{\mathcal{H}}$ and $1\leq l\leq\infty$ be such that $\frac{n}{l}%
+\sigma_{G}\leq0$. Then, the Riesz transforms are bounded operators on
$\mathcal{M}\left[  G\right]  ^{l}.$
\end{lemma}

The subject of the next lemma is a basic estimate for convolution that allows
us to get, in particular, a heat semigroup inequality.

\begin{lemma}
\label{Lem: Aux limitacion Linfty para Bloques} Let $1\leq l\leq\infty$ and
$G\in\mathcal{H}$ be such that $\frac{n}{l}+\sigma_{G}\geq0,$ and let $b$ be a
$(l,G)$-block with $\mbox{supp}(b)\subset D\left(  x_{0},\rho_{0}\right)  $.
Then, we have that
\[
\left\Vert \psi\ast b\right\Vert _{L^{\infty}}\leq C,
\]
for every radially symmetric $\psi\in\mathcal{S}$ , where $C=C(\psi)>0$ is
independent of $x_{0}$ and $\rho_{0}$.
\end{lemma}

\textbf{Proof. }Defining $\mu\left(  \rho\right)  =\int_{D\left(
0,\rho\right)  }b(x)dx$, we proceed as follows
\[
\left(  \psi\ast b\right)  \left(  0\right)  =\int\limits_{\mathbb{R}^{n}}%
\psi\left(  y\right)  b\left(  y\right)  dy=\int\limits_{0}^{\infty}%
\psi\left(  \rho\right)  d\mu\left(  \rho\right)  =-\int\limits_{0}^{\infty
}\frac{\partial\psi\left(  \rho\right)  }{\partial\rho}\mu\left(  \rho\right)
d\rho.
\]
Moreover, by item \textbf{\textit{(iii)}} of Proposition
\ref{prop:conseq-condition-morrey-space}, we obtain that
\begin{align*}
\mu\left(  \rho\right)   &  =\int\limits_{D\left(  0,\rho\right)  }1_{D\left(
x_{0},\rho_{0}\right)  }\left(  x\right)  b(x)dx=\int\limits_{D\left(
0,\rho\right)  \cap D\left(  x_{0},\rho_{0}\right)  }b(x)dx\leq C\min\left\{
\rho^{n+\sigma_{G}},\rho_{0}^{n+\sigma_{G}}\right\}  \left\Vert 1_{D\left(
x_{0},\rho_{0}\right)  }b\right\Vert _{G}\\
&  \leq C\min\left\{  \rho^{n+\sigma_{G}}\rho_{0}^{-\left(  \frac{n}{l}%
+\sigma_{G}\right)  },\rho_{0}^{n-\frac{n}{l}}\right\}  .
\end{align*}
Thus,
\begin{align*}
\left\vert \left(  \psi\ast b\right)  \left(  0\right)  \right\vert  &  \leq
C\int\limits_{0}^{\infty}\left\vert \frac{\partial\psi\left(  \rho\right)
}{\partial\rho}\right\vert \min\left\{  \rho^{n+\sigma_{G}}\rho_{0}^{-\left(
\frac{n}{l}+\sigma_{G}\right)  },\rho_{0}^{n-\frac{n}{l}}\right\}  d\rho\\
&  =C\int\limits_{0}^{\rho_{0}}\left\vert \frac{\partial\psi\left(
\rho\right)  }{\partial\rho}\right\vert \rho^{n+\sigma_{G}}\rho_{0}^{-\left(
\frac{n}{l}+\sigma_{G}\right)  }d\rho+C\int\limits_{\rho_{0}}^{\infty
}\left\vert \frac{\partial\psi\left(  \rho\right)  }{\partial\rho}\right\vert
\rho_{0}^{n-\frac{n}{l}}d\rho\\
&  \leq C\int\limits_{0}^{\rho_{0}}\left\vert \frac{\partial\psi\left(
\rho\right)  }{\partial\rho}\right\vert \rho^{n-\frac{n}{l}}d\rho
+C\int\limits_{\rho_{0}}^{\infty}\left\vert \frac{\partial\psi\left(
\rho\right)  }{\partial\rho}\right\vert \rho^{n-\frac{n}{l}}d\rho\\
&  \leq C(\psi).
\end{align*}
For $x\in\mathbb{R}^{n},$ note that $\left(  \psi\ast b\right)  \left(
x\right)  =\left(  \psi\ast\tau_{x}\tilde{b}\right)  \left(  0\right)  $.
Then, using \textbf{(M3)} and proceeding as before, we get
\begin{equation}
\left\vert \left(  \psi\ast b\right)  \left(  x\right)  \right\vert \leq C.
\label{aux-conv-1}%
\end{equation}
Now we can complete the proof by taking the supremum over $x\in\mathbb{R}^{n}$
in (\ref{aux-conv-1}). \fin

As a consequence of the previous lemma, we have the following estimate in
$(l,X)$-block spaces.

\begin{corollary}
\label{Cor: limitacion norma Linft deltaf con norma de bloques.} Let $\psi
\in\mathcal{S}$ be radially symmetric and $1\leq l\leq\infty.$ Then, there
exists a constant $C>0$ depending on $\psi$ such that
\[
\left\Vert \psi\ast f\right\Vert _{L^{\infty}}\leq C\left\Vert f\right\Vert
_{\mathcal{PD}\left[  X\right]  ^{l}},
\]
for all $f\in\mathcal{PD}\left[  X\right]  ^{l}$.
\end{corollary}

\textbf{Proof. }We have that $\psi\ast f=\sum\limits_{k=1}^{\infty}\alpha
_{k}\psi\ast b_{k}$, where $f=\sum\limits_{k=1}^{\infty}\alpha_{k}b_{k}$ is a
decomposition in $(l,X)$-blocks for $f$. Using Lemma
\ref{Lem: Aux limitacion Linfty para Bloques}, we can estimate
\[
\left\Vert \psi\ast f\right\Vert _{L^{\infty}}\leq\sum\limits_{k=1}^{\infty
}\left\vert \alpha_{k}\right\vert \left\Vert \psi\ast b_{k}\right\Vert
_{L^{\infty}}\leq C\sum\limits_{k=1}^{\infty}\left\vert \alpha_{k}\right\vert
,
\]
which gives the desired estimate after taking the infimum over all the
decompositions of $f$. \fin

In the next two lemmas, we show the basic convolution estimates in $X$-Morrey
spaces useful for our purposes.

\begin{lemma}
\label{Lem: Limitacion convolucion norma Linfty por norma el Morrey}Let
$\psi\in\mathcal{S}$ be radially symmetric, $G\in\mathcal{H}$ and $1\leq
l\leq\infty.$ Then, there exists a constant $C>0$ depending on $\psi$ such
that
\begin{equation}
\left\Vert \psi\ast f\right\Vert _{L^{\infty}}\leq C\left\Vert f\right\Vert
_{\mathcal{M}\left[  G\right]  ^{l}}, \label{aux-conv-2}%
\end{equation}
for all $f\in\mathcal{M}[G]^{l}$.
\end{lemma}

\textbf{Proof. }Proceeding as in the proof of Lemma
\ref{Lem: Aux limitacion Linfty para Bloques}, we have that
\[
\left(  \psi\ast f\right)  \left(  0\right)  =-\int\limits_{0}^{\infty}%
\frac{\partial\psi\left(  \rho\right)  }{\partial\rho}\mu\left(  \rho\right)
d\rho,
\]
where $\mu\left(  \rho\right)  =\int_{D\left(  0,\rho\right)  }f(x)dx$. Again,
by Proposition \ref{prop:conseq-condition-morrey-space} \textbf{\textit{(iii)}%
}, it follows that
\begin{align*}
\left\vert \mu\left(  \rho\right)  \right\vert  &  \leq C\left\Vert
1_{D\left(  0,\rho\right)  }f\right\Vert _{G}\rho^{n+\sigma_{G}}%
=C\rho^{n+\sigma_{G}}\rho^{-\left(  \frac{n}{l}+\sigma_{G}\right)  }%
\rho^{\frac{n}{l}+\sigma_{G}}\left\Vert 1_{D\left(  0,\rho\right)
}f\right\Vert _{G}\\
&  \leq C\rho^{\frac{n}{l^{\prime}}}\left\Vert f\right\Vert _{\mathcal{M}%
\left[  G\right]  ^{l}}.
\end{align*}
Therefore,
\begin{align*}
\left\vert \left(  \psi\ast f\right)  \left(  0\right)  \right\vert  &  \leq
C\int\limits_{0}^{\infty}\left\vert \frac{\partial\psi\left(  \rho\right)
}{\partial\rho}\right\vert \rho^{\frac{n}{l^{\prime}}}\left\Vert f\right\Vert
_{\mathcal{M}\left[  G\right]  ^{l}}d\rho\\
&  \leq C\left\Vert f\right\Vert _{\mathcal{M}\left[  G\right]  ^{l}}%
\int\limits_{0}^{\infty}\left\vert \frac{\partial\psi\left(  \rho\right)
}{\partial\rho}\right\vert \rho^{\frac{n}{l^{\prime}}}d\rho\\
&  \leq C\left\Vert f\right\Vert _{\mathcal{M}\left[  G\right]  ^{l}}.
\end{align*}
For $x\in\mathbb{R}^{n},$ noting $\left(  \psi\ast f\right)  \left(  x\right)
=\left(  \psi\ast\tau_{x}f\right)  \left(  0\right)  $ and using
\textbf{(M3)}, we arrive at
\[
\left\vert \left(  \psi\ast f\right)  \left(  x\right)  \right\vert \leq
C\left\Vert f\right\Vert _{\mathcal{M}\left[  G\right]  ^{l}},
\]
which gives (\ref{aux-conv-2}) by taking the supremum over $x\in\mathbb{R}%
^{n}$. \fin

\begin{lemma}
\label{Lem: Limitacion convolucion norma Morrey por norma Morrey}Let
$G\in\mathcal{H}$ and $1\leq l\leq\infty$. Then%
\begin{equation}
\left\Vert \psi\ast g\right\Vert _{\mathcal{M}\left[  G\right]  ^{l}}\leq
C\left\Vert \psi\right\Vert _{L^{1}}\left\Vert g\right\Vert _{\mathcal{M}%
\left[  G\right]  ^{l}}, \label{aux-conv-3}%
\end{equation}
for all $g\in\mathcal{M}\left[  G\right]  ^{l}$ and $\psi\in L^{1}$, where
$C>0$ is a universal constant.
\end{lemma}

\textbf{Proof. }For $x_{0}\in\mathbb{R}^{n}$ and $\rho_{0}>0,$ we can express
\[
1_{D\left(  x_{0},\rho_{0}\right)  }\left(  x\right)  \left(  \psi\ast
g\right)  \left(  x\right)  =\int\limits_{\mathbb{R}^{n}}1_{D\left(
x_{0},\rho_{0}\right)  }\left(  x\right)  \psi\left(  y\right)  g\left(
x-y\right)  dy=\int\limits_{\mathbb{R}^{n}}\psi\left(  y\right)  1_{D\left(
x_{0}-y,\rho_{0}\right)  }\left(  x-y\right)  g\left(  x-y\right)  dy.
\]
Thus, condition \textbf{(M4)} with $K(x,y)=1_{D\left(  x_{0}-y,\rho
_{0}\right)  }\left(  x-y\right)  g\left(  x-y\right)  $ implies that
\begin{align*}
\left\Vert 1_{D\left(  x_{0},\rho_{0}\right)  }\left(  \psi\ast g\right)
\right\Vert _{G}  &  \leq\int\limits_{\mathbb{R}^{n}}\left\vert \psi\left(
y\right)  \right\vert \left\Vert \tau_{y}\left(  1_{D\left(  x_{0}-y,\rho
_{0}\right)  }g\right)  \right\Vert _{G}dy\leq C\int\limits_{\mathbb{R}^{n}%
}\left\vert \psi\left(  y\right)  \right\vert \left\Vert 1_{D\left(
x_{0}-y,\rho_{0}\right)  }g\right\Vert _{G}dy\\
&  \leq C\rho_{0}^{-\left(  \frac{n}{l}+\sigma_{G}\right)  }\left\Vert
\psi\right\Vert _{L^{1}}\left\Vert g\right\Vert _{\mathcal{M}\left[  G\right]
^{l}},
\end{align*}
and then
\[
\rho_{0}^{\frac{n}{l}+\sigma_{G}}\left\Vert 1_{D\left(  x_{0},\rho_{0}\right)
}\left(  \psi\ast g\right)  \right\Vert _{G}\leq C\left\Vert \psi\right\Vert
_{L^{1}}\left\Vert g\right\Vert _{\mathcal{M}\left[  G\right]  ^{l}}.
\]
Taking the supremum over $x_{0}\in\mathbb{R}^{n}$ and $\rho_{0}>0,$ we are
done. \fin

From duality and the previous lemma, we obtain the corollary below.

\begin{corollary}
\label{Cor: convolucion en Bloques} Let $X\in\mathcal{H}$ be such that
$X^{\prime}\in\mathcal{H}$ and $1\leq l\leq\infty$. Then
\[
\left\Vert \psi\ast f\right\Vert _{\mathcal{PD}\left[  X\right]  ^{l}}\leq
C\left\Vert \psi\right\Vert _{L^{1}}\left\Vert f\right\Vert _{\mathcal{PD}%
\left[  X\right]  ^{l}},
\]
for all $f\in\mathcal{PD}\left[  X\right]  ^{l}$ and $\psi\in L^{1}$, where
$C>0$ is a universal constant.
\end{corollary}

\textbf{Proof. }By (\ref{aux-duality-Morrey-1}) and the canonical embedding
into the bidual, we have the linear isometry (not surjective) $\mathcal{PD}%
\left[  X\right]  ^{l}\hookrightarrow\left(  \mathcal{M}\left[  X^{\prime
}\right]  ^{l^{\prime}}\right)  ^{\prime}.$ Using it and estimate
(\ref{aux-conv-3}), we proceed as follows%
\begin{align*}
\left\Vert \psi\ast f\right\Vert _{\mathcal{PD}\left[  X\right]  ^{l}}  &
=\mathop{\sup}\limits_{\left\Vert g\right\Vert _{M\left[  X^{\prime}\right]
^{l^{\prime}}}=1}\left\vert \left\langle \psi\ast f,g\right\rangle \right\vert
\\
&  =\mathop{\sup}\limits_{\left\Vert g\right\Vert _{M\left[  X^{\prime
}\right]  ^{l^{\prime}}}=1}\left\vert \left\langle f,\tilde{\psi}\ast
g\right\rangle \right\vert \\
&  \leq\left\Vert f\right\Vert _{\mathcal{PD}\left[  X\right]  ^{l}%
}\mathop{\sup}\limits_{\left\Vert g\right\Vert _{M\left[  X^{\prime}\right]
^{l^{\prime}}}=1}\left\Vert \psi\ast g\right\Vert _{\mathcal{M}\left[
X^{\prime}\right]  ^{l^{\prime}}}\\
&  \leq C\left\Vert \psi\right\Vert _{L^{1}}\left\Vert f\right\Vert
_{\mathcal{PD}\left[  X\right]  ^{l}}\mathop{\sup}\limits_{\left\Vert
g\right\Vert _{M\left[  X^{\prime}\right]  ^{l^{\prime}}}=1}\left\Vert
g\right\Vert _{\mathcal{M}\left[  X^{\prime}\right]  ^{l^{\prime}}}\\
&  \leq C\left\Vert \psi\right\Vert _{L^{1}}\left\Vert f\right\Vert
_{\mathcal{PD}\left[  X\right]  ^{l}},
\end{align*}
as required.\fin

Now, using interpolation theory, we can extend Lemma
\ref{Lem: Limitacion convolucion norma Morrey por norma Morrey} as follows.

\begin{lemma}
\label{Lem: Limitacion de convolucion norma Morrey Q, por norma Morrey N} Let
$\psi\in\mathcal{S}$ be radially symmetric and $N,Q\in\mathcal{H}$. Suppose
that $\left(  L^{\infty},N\right)  _{\theta,d}\hookrightarrow Q$ with
$\sigma_{Q}=\theta\sigma_{N}$ for some $0<\theta<1$ and $1\leq l,d\leq\infty.$
Then, there exists a constant $C>0$ depending on $\psi$ such that
\[
\left\Vert \psi\ast f\right\Vert _{\mathcal{M}\left[  Q\right]  ^{l/\theta}%
}\leq C\left\Vert f\right\Vert _{\mathcal{M}\left[  N\right]  ^{l}},
\]
for all $f\in\mathcal{M}\left[  N\right]  ^{l}$.
\end{lemma}

\textbf{Proof. }First, it follows from Lemma
\ref{Lem: Limitacion convolucion norma Linfty por norma el Morrey} that%

\begin{equation}
\left\Vert 1_{D\left(  x_{0},\rho_{0}\right)  }\psi\ast f\right\Vert
_{L^{\infty}}\leq C_{1}(\psi)\left\Vert f\right\Vert _{\mathcal{M}\left[
N\right]  ^{l}}. \label{aux-conv-2001}%
\end{equation}
On the other hand, by Lemma
\ref{Lem: Limitacion convolucion norma Morrey por norma Morrey}, we know that
\begin{equation}
\left\Vert 1_{D\left(  x_{0},\rho_{0}\right)  }\psi\ast f\right\Vert _{N}%
\leq\rho_{0}^{-\left(  \frac{n}{l}+\sigma_{N}\right)  }C_{2}\left\Vert
\psi\right\Vert _{L^{1}}\left\Vert f\right\Vert _{\mathcal{M}\left[  N\right]
^{l}}. \label{aux-conv-2002}%
\end{equation}
Using (\ref{aux-conv-2001})-(\ref{aux-conv-2002}) and interpolation, we arrive
at
\[
\left\Vert 1_{D\left(  x_{0},\rho_{0}\right)  }\psi\ast f\right\Vert _{Q}\leq
C\left\Vert 1_{D\left(  x_{0},\rho_{0}\right)  }\psi\ast f\right\Vert
_{\left(  L^{\infty},N\right)  _{\theta,d}}\leq\rho_{0}^{-\theta\left(
\frac{n}{l}+\sigma_{N}\right)  }C_{3}(\psi)\left\Vert f\right\Vert
_{\mathcal{M}\left[  N\right]  ^{l}},
\]
which yields the desired estimate after multiplying both sides by $\rho
_{0}^{\theta\left(  \frac{n}{l}+\sigma_{N}\right)  }$ and taking the supremum
over $x_{0}\in\mathbb{R}^{n}$ and $\rho_{0}>0.$\fin

Recall that $U(t)$ denotes the heat semigroup, which is defined for $t>0$ as
$U(t)f=\exp(t\Delta)f=\Phi(t)\ast f$, for each $f\in\mathcal{S}^{\prime}$,
where $\Phi(x,t)=\left(  4\pi t\right)  ^{-n/2}e^{-\frac{\left\vert
x\right\vert ^{2}}{4t}}$ is the heat kernel. For each $m\in\mathbb{N}_{0}$, it
follows that
\[
\left\vert \left(  \nabla_{x}^{m}\Phi\right)  (x,t)\right\vert \leq\bar
{h}\left(  x,t\right)  =h\left(  \left\vert x\right\vert ,t\right)  ,
\]
where $\bar{h}\left(  x,t\right)  =t^{-\frac{m+n}{2}}\bar{h}\left(  \frac
{x}{\sqrt{t}},1\right)  ,$ $h\left(  \rho,t\right)  \in\mathcal{S}\left(
\mathbb{R}\right)  $ and $\partial_{\rho}h\left(  \rho,t\right)
=t^{-\frac{m+n+1}{2}}\left(  \partial_{\rho}h\right)  \left(  \frac{\rho
}{\sqrt{t}},1\right)  ,$ for all $t>0$.

We are in position to present our estimates of the heat semigroup in the
context of $X$-blocks and $X$-Morrey spaces.

\begin{corollary}
\label{Lem: Estimativa nucleo del calor en Morrey} Let $N,Q\in\mathcal{H}$ and
$m\in\mathbb{N}_{0}$. Suppose that $\left(  L^{\infty},N\right)  _{\theta
,d}\hookrightarrow Q$ and $\sigma_{Q}=\theta\sigma_{N}$ for some $0<\theta<1$
and $1\leq l,d\leq\infty.$ Then,
\[
\left\Vert \nabla^{m}U(t)f\right\Vert _{\mathcal{M}\left[  Q\right]
^{l/\theta}}\leq Ct^{-\frac{1}{2}\left(  \frac{n}{l}-\theta\frac{n}{l}\right)
-\frac{m}{2}}\left\Vert f\right\Vert _{\mathcal{M}\left[  N\right]  ^{l}},
\]
for all $f\in\mathcal{M}[N]^{l}$, where the constant $C>0$ is independent of
$t>0.$
\end{corollary}

\textbf{Proof.} First, we have the inequalities
\[
\left\vert \left(  \nabla^{m}U(t)f\right)  \left(  \cdot\right)  \right\vert
\leq\left(  \bar{h}\left(  \cdot,t\right)  \ast\left\vert f\right\vert
\right)  \left(  \cdot\right)  \leq t^{-\frac{m+n}{2}}\left(  \bar{h}\left(
t^{-\frac{1}{2}}\cdot,1\right)  \ast\left\vert f\right\vert \right)  \left(
\cdot\right)  =t^{-\frac{m}{2}}\left(  \bar{h}\left(  \cdot,1\right)
\ast\left\vert f\right\vert \left(  t^{\frac{1}{2}}\cdot\right)  \right)
\left(  t^{-\frac{1}{2}}\cdot\right)  .
\]
So, Lemma
\ref{Lem: Limitacion de convolucion norma Morrey Q, por norma Morrey N}
yields
\begin{align*}
\left\Vert \nabla^{m}U(t)f\right\Vert _{\mathcal{M}\left[  Q\right]
^{l/\theta}}  &  \leq\left\Vert t^{-\frac{m}{2}}\left(  \bar{h}\left(
\cdot,1\right)  \ast\left\vert f\right\vert \left(  t^{\frac{1}{2}}%
\cdot\right)  \right)  \left(  t^{-\frac{1}{2}}\cdot\right)  \right\Vert
_{\mathcal{M}\left[  Q\right]  ^{l/\theta}}\\
&  \leq Ct^{\left(  -m+\frac{\theta n}{l}\right)  /2}\left\Vert \bar{h}\left(
\cdot,1\right)  \ast\left\vert f\right\vert \left(  t^{\frac{1}{2}}%
\cdot\right)  \right\Vert _{\mathcal{M}\left[  Q\right]  ^{l/\theta}}\\
&  \leq Ct^{\left(  -m+\frac{\theta n}{l}\right)  /2}\left\Vert f\left(
t^{\frac{1}{2}}\cdot\right)  \right\Vert _{\mathcal{M}\left[  N\right]  ^{l}%
}\\
&  \leq Ct^{\left(  -m-\frac{n}{l}(1-\theta)\right)  /2}\left\Vert
f\right\Vert _{\mathcal{M}\left[  N\right]  ^{l}}.
\end{align*}
\fin

\begin{corollary}
\label{Cor: Estimativa nucleo del calor en Blocos} Let $m\in\mathbb{N}_{0}$
and let $F$,$G\in\mathcal{H}$ be such that $F^{\prime},G^{\prime}%
\in\mathcal{H}$, and suppose that $\left(  G^{\prime},L^{\infty}\right)
_{\theta,d}\hookrightarrow F^{\prime}$ and $\sigma_{F^{\prime}}=\theta
\sigma_{G^{\prime}}$, for some $0<\theta<1$ and $1\leq l,d\leq\infty.$ Then,
for $\tilde{l}=\theta l/[1-(1-\theta)l]$ we have that%
\[
\left\Vert \nabla^{m}U(t)f\right\Vert _{\mathcal{PD}\left[  G\right]
^{\tilde{l}}}\leq Ct^{-\frac{1}{2}\left(  \frac{n}{l}-\frac{n}{\tilde{l}%
}\right)  -\frac{m}{2}}\left\Vert f\right\Vert _{\mathcal{PD}\left[  F\right]
^{l}},
\]

for all $f\in\mathcal{PD}\left[  F\right]  ^{l}$, where the constant $C>0$ is
independent of $t>0.$
\end{corollary}

\textbf{Proof. }Using duality, the previous lemma, and noting that $l^{\prime
}=\tilde{l}^{\prime}/\theta$, we can estimate
\begin{align*}
\left\Vert \nabla^{m}U(t)f\right\Vert _{\mathcal{PD}\left[  G\right]
^{\tilde{l}}}  &  =\mathop{\sup}\limits_{\left\Vert g\right\Vert
_{\mathcal{M}[G^{\prime}]^{\tilde{l}^{\prime}}}=1}\left\vert \left\langle
\nabla^{m}U(t)f,g\right\rangle \right\vert \\
&  =\mathop{\sup}\limits_{\left\Vert g\right\Vert _{\mathcal{M}[G^{\prime
}]^{\tilde{l}^{\prime}}}=1}\left\vert \left\langle f,\nabla^{m}%
U(t)g\right\rangle \right\vert \\
&  \leq\mathop{\sup}\limits_{\left\Vert g\right\Vert _{\mathcal{M}[G^{\prime
}]^{\tilde{l}^{\prime}}}=1}\left\Vert f\right\Vert _{\mathcal{PD}\left[
F\right]  ^{l}}\left\Vert \nabla^{m}U(t)g\right\Vert _{\mathcal{M}\left[
F^{\prime}\right]  ^{\tilde{l}^{\prime}/\theta}}\\
&  \leq C\mathop{\sup}\limits_{\left\Vert g\right\Vert _{\mathcal{M}%
[G^{\prime}]^{\tilde{l}^{\prime}}}=1}\left\Vert f\right\Vert _{\mathcal{PD}%
\left[  F\right]  ^{l}}t^{-\frac{1}{2}\left(  \frac{n}{\tilde{l}^{\prime}%
}-\theta\frac{n}{\tilde{l}^{\prime}}\right)  -\frac{m}{2}}\left\Vert
g\right\Vert _{\mathcal{M}\left[  G^{\prime}\right]  ^{\tilde{l}^{\prime}}%
}\leq Ct^{-\frac{1}{2}\left(  \frac{n}{l}-\frac{n}{\tilde{l}}\right)
-\frac{m}{2}}\left\Vert f\right\Vert _{\mathcal{PD}\left[  F\right]  ^{l}}.
\end{align*}
\fin

The following interpolation result will be useful for our ends.

\begin{lemma}
\label{Lem: Interpolacion de operadores Blocos} Let $X,X_{1},X_{2}%
\in\mathcal{H}$ be such that $\sigma_{X}=(1-\theta)\sigma_{X_{1}}+\theta
\sigma_{X_{2}}$ and $X\hookrightarrow\left(  X_{1},X_{2}\right)  _{\theta,d}$
for some $\theta\in(0,1)$ and $d\in\left[  1,\infty\right]  $. If $F_{1}$ and
$F_{2}$ are Banach spaces and $T$ is a continuous linear operator such that
\begin{equation}
T:\mathcal{PD}\left[  X_{1}\right]  ^{l_{1}}\longrightarrow F_{1}%
\,\,\,\mbox{and}\,\,\,T:\mathcal{PD}\left[  X_{2}\right]  ^{l_{2}%
}\longrightarrow F_{2}, \label{eq: aux operador Bloques}%
\end{equation}
with the operator norms $C_{1}$ and $C_{2}$, respectively. Then, for $\frac
{1}{l}=\frac{1-\theta}{l_{1}}+\frac{\theta}{l_{2}}$ we have
\[
T:\mathcal{PD}\left[  X\right]  ^{l}\longrightarrow\left(  F_{1},F_{2}\right)
_{\theta,d},
\]
with the operator norm bounded by $\tilde{C}=\left(  C_{1}\right)  ^{1-\theta
}\left(  C_{2}\right)  ^{\theta}.$
\end{lemma}

\textbf{Proof. }Note that, for every ball $D\left(  x,\rho\right)  ,$ we can
define
\[%
\begin{array}
[c]{cccc}%
T_{\left(  x,\rho\right)  }: & X_{1} & \longrightarrow & F_{1}\\
& f & \longmapsto & T\left(  1_{D\left(  x,\rho\right)  }f\right)
\end{array}
;\,\,\,\,%
\begin{array}
[c]{cccc}%
T_{\left(  x,\rho\right)  }: & X_{2} & \longrightarrow & F_{2}\\
& f & \longmapsto & T\left(  1_{D\left(  x,\rho\right)  }f\right)
\end{array}
,
\]
with norms $C_{1}\leq C\rho^{\left(  \frac{n}{l_{1}}+\sigma_{X_{1}}\right)  }$
and $C_{2}\leq C\rho^{\left(  \frac{n}{l_{2}}+\sigma_{X_{2}}\right)  }.$
Thus,
\[
T_{\left(  x,\rho\right)  }:X\longrightarrow\left(  F_{1},F_{2}\right)
_{\theta,d},
\]
with norm $C\leq C_{1}^{1-\theta}C_{2}^{\theta}\leq C\rho^{\left(
1-\theta\right)  \left(  \frac{n}{l_{1}}+\sigma_{X_{1}}\right)  }\rho
^{\theta\left(  \frac{n}{l_{2}}+\sigma_{X_{2}}\right)  }\leq C\rho^{\frac
{n}{l}+\sigma_{X}}$. For every block $b_{k}$ with $\mbox{supp}b_{k}\subset
D\left(  x_{k},\rho_{k}\right)  ,$ we have that
\[
\left\Vert T\left(  1_{D\left(  x_{k},\rho_{k}\right)  }b_{k}\right)
\right\Vert _{\left(  F_{1},F_{2}\right)  _{\theta,d}}=\left\Vert T_{\left(
x_{k},\rho_{k}\right)  }\left(  1_{D\left(  x_{k},\rho_{k}\right)  }%
b_{k}\right)  \right\Vert _{\left(  F_{1},F_{2}\right)  _{\theta,d}}\leq
C\rho^{\rho^{\frac{n}{l}+\sigma_{X}}}\left\Vert 1_{D\left(  x_{k},\rho
_{k}\right)  }b_{k}\right\Vert _{X}\leq C.
\]
Now, consider the space
\[
\mathcal{LD}\left[  X\right]  ^{l}=\left\{  f\in\mathcal{PD}\left[  X\right]
^{l};f=\sum\limits_{k=1}^{m}\alpha_{k}b_{k}\mbox{ where }m\in\mathbb{N}%
\mbox{ and }b_{k}\mbox{ is an }\left(  l,X\right)  \mbox{-block}\right\}  .
\]
Then, for $f\in\mathcal{LD}\left[  X\right]  ^{l}$ we obtain
\[
\left\Vert T\left(  f\right)  \right\Vert _{\left(  F_{1},F_{2}\right)
_{\theta,d}}\leq\sum\limits_{k=1}^{m}\left\vert \alpha_{k}\right\vert
\left\Vert T\left(  1_{D\left(  x_{k},\rho_{k}\right)  }b_{k}\right)
\right\Vert _{\left(  F_{1},F_{2}\right)  _{\theta,d}}\leq C\sum
\limits_{k=1}^{m}\left\vert \alpha_{k}\right\vert \leq C\left\Vert
f\right\Vert _{\mathcal{PD}\left[  X\right]  ^{l}}.
\]
Using that $\mathcal{LD}\left[  X\right]  ^{l}$ is dense in $\mathcal{PD}%
\left[  X\right]  ^{l},$ we can conclude the proof.

\fin

\begin{remark}
\label{remark:Interpolacion de operadores Blocos} Note that in the conditions
of the previous lemma we have that
\[
\mathcal{PD}\left[  X\right]  ^{l}\hookrightarrow\left(  \mathcal{PD}\left[
X_{1}\right]  ^{l_{1}},\mathcal{PD}\left[  X_{2}\right]  ^{l_{2}}\right)
_{\theta,d}.
\]

\end{remark}

With the above properties in hand, we are in a position to show the bilinear
estimate for the Navier-Stokes equations in $X$-Morrey spaces.

\begin{theorem}
\label{Teo:estim-bilin-cambio} Let $1\leq l<l_{0}<\infty$ and $X_{0}%
,X\in\mathcal{H}$ be such that $X_{0}^{\prime},X^{\prime}\in\mathcal{H}$,
$\sigma_{X}<\sigma_{X_{0}}$ and $l_{0}^{\prime}\sigma_{X_{0}^{\prime}%
}=l^{\prime}\sigma_{X^{\prime}}$. Assume the following conditions:

\begin{itemize}
\item[\textbf{(i)}] There exist Banach spaces $X_{1},X_{2}\in\mathcal{H}$ such
that
\begin{equation}
\theta_{1}<\frac{l_{0}^{\prime}}{l^{\prime}}<\theta_{2}<\min\left\{
\frac{l_{0}^{\prime}}{l^{\prime}}+\frac{2l_{0}^{\prime}}{n},1\right\}  ,
\label{eq:estim-bilin-weak-morrey-1}%
\end{equation}
where $\theta_{i}=\frac{\sigma_{X_{i}^{\prime}}}{\sigma_{X_{0}^{\prime}}}$,
$i=1,2$.

\item[\textbf{(ii)}] For $\theta\in(0,1)$ such that $\frac{l_{0}^{\prime}%
}{l^{\prime}}=\theta\theta_{2}+(1-\theta)\theta_{1}$, it follows that
$X\hookrightarrow\left(  X_{1},X_{2}\right)  _{\theta,1}$ and $\left(
L^{\infty},X_{0}^{\prime}\right)  _{\theta_{i},d_{i}}\hookrightarrow
X_{i}^{\prime}$, for some $1\leq d_{i}\leq\infty$, $i=1,2$.

\item[\textbf{(iii)}] $X_{0}^{\prime}\in\tilde{\mathcal{H}}$ and $\left\Vert
fg\right\Vert _{X_{0}^{\prime}}\leq C\left\Vert f\right\Vert _{X^{\prime}%
}\left\Vert g\right\Vert _{X^{\prime}}$.
\end{itemize}

Then, under \textbf{\textit{(i)}} and \textbf{\textit{(ii)}}, for
$E=\mathcal{PD}[X]^{l}$ and $E_{0}=\mathcal{PD}[X_{0}]^{l_{0}}$, conditions
\textbf{(H1)} and \textbf{(H4)} hold true. Moreover, if we assume
\textbf{\textit{(iii)}}, we obtain \textbf{(H2)} and \textbf{(H3)}. In
particular, for $l^{\prime}=n$, we obtain the bilinear estimate
\eqref{eq:Estimativa bilinear} in the space $E^{\prime}=\mathcal{M}[X^{\prime
}]^{n}$.
\end{theorem}

\noindent\textbf{Proof.} First, we assume \textbf{\textit{(i)}} and
\textbf{\textit{(ii)}}. Condition \textbf{(H1)} is a direct consequence of
Lemma \ref{Lem: espacio de bloques pertenece a la clase G}. To show
\textbf{(H4)}, we need only to check the conditions in Lemma
\ref{Lem: condiciones para estimativa integral}. In fact, define
$l_{i}^{\prime}=\frac{l_{0}^{\prime}}{\theta_{i}}$ and $E_{i}=\mathcal{PD}%
[X_{i}]^{l_{i}}$, for $i=1,2$. Using Corollary
\ref{Cor: Estimativa nucleo del calor en Blocos} with $G=X_{0}$, $F=X_{i}$ and
$\theta=\theta_{i}$, by item \textbf{\textit{(ii)}} we get
\[
\Vert\nabla U(t)f\Vert_{\mathcal{PD}[X_{0}]^{l_{0}}}\leq Ct^{\frac{1}%
{2}\left(  \frac{n}{l_{0}}-\frac{n}{l_{i}}\right)  -\frac{1}{2}}\Vert
f\Vert_{\mathcal{PD}[X_{i}]^{l_{i}}},\text{ for }i=1,2,
\]
which gives \textbf{(H4$^\prime$)} by observing that $\sigma_{E_{0}}=-\frac{n}{l_{0}}$
and $\sigma_{E_{i}}=-\frac{n}{l_{i}}$.

Next, by taking $\frac{1}{z_{i}}=\frac{1}{2}\left(  \frac{n}{l_{i}}-\frac
{n}{l}\right)  +1$, inequalities in (\ref{eq:estim-bilin-weak-morrey-1}) imply
that $\frac{1}{z_{1}}>1$ and $0<\frac{1}{z_{2}}<1$. Moreover, by the choice of
$\theta$, we arrive at
\begin{align*}
\frac{\theta}{z_{2}}+\frac{1-\theta}{z_{1}}  &  =\frac{n}{2}\left(  \frac
{1}{l^{\prime}}-\frac{\theta}{l_{2}^{\prime}}-\frac{1-\theta}{l_{1}^{\prime}%
}\right)  +1\\
&  =\frac{n}{2}\left(  \frac{1}{l^{\prime}}-\frac{\theta\theta_{2}%
+(1-\theta)\theta_{1}}{l_{0}^{\prime}}\right)  +1\\
&  =1,
\end{align*}
and then we obtain \textbf{(H4$^{\prime \prime}$)}.

In order to show \textbf{(H4$^{\prime \prime \prime}$ )}, we first note that the choice of $\theta$ in
item \textbf{\textit{(ii)}} is equivalent to $\sigma_{X}=\theta\sigma_{X_{2}%
}+(1-\theta)\sigma_{X_{1}}$. Then, Lemma
\ref{Lem: Interpolacion de operadores Blocos} and Remark
\ref{remark:Interpolacion de operadores Blocos} give
\[
\mathcal{PD}\left[  X\right]  ^{l}\hookrightarrow\left(  \mathcal{PD}\left[
X_{1}\right]  ^{l_{1}},\mathcal{PD}\left[  X_{2}\right]  ^{l_{2}}\right)
_{\theta,1},
\]
as required. Since we have showed \textbf{(H4$^\prime$)-(H4$^{\prime \prime \prime}$ )}, condition \textbf{(H4)}
follows from Lemma \ref{Lem: condiciones para estimativa integral}.

Assuming now \textbf{\textit{(iii)}}, a direct application of Lemma
\ref{Lem: Limitacion de operador de Riesz en Morrey} yields \textbf{(H2)}. On
the other hand, item \textbf{\textit{(iii)}} implies that $\frac{l^{\prime}%
}{l_{0}^{\prime}}=\frac{\sigma_{X_{0}^{\prime}}}{\sigma_{X^{\prime}}}=2$, and
thus Lemma \ref{Lem: producto en morrey} implies \textbf{(H3)}.

Finally, if $l^{\prime}=n$, the conditions in Theorem
\ref{Teo: Teo main. Estimativa integral en el dual.} are fulfilled by using the
spaces $E=\mathcal{PD}[X]^{l}$ and $E_{0}=\mathcal{PD}[X_{0}]^{l_{0}}$, and
then we obtain (\ref{eq:Estimativa bilinear}) with $E^{\prime}=\mathcal{M}%
[X]^{n},$ as desired. \fin

\begin{remark}
Note that the condition \textbf{(iii)} in Theorem \ref{Teo:estim-bilin-cambio}
and Remark \ref{remark:restriction-on-scaling} imply that $\sigma_{X^{\prime}%
}>-\frac{n}{2}$.
\end{remark}

\subsection{Application (bilinear estimate in weak-Morrey spaces)}

\label{eq:bilinear-estimate-weak-morrey-spaces} We wish to apply our result
for the weak-Morrey spaces $\mathcal{M}_{\left(  p^{\prime},\infty\right)
}^{n}=\mathcal{M}\left[  L^{(p^{\prime},\infty)}\right]  ^{n}$ and obtain the
corresponding bilinear estimate in \cite{Lucas,Lemarie2}. For that, take
$p,p_{0},l,l_{0}$ such that $l^{\prime}=n$, $l_{0}^{\prime}=n/2$,
$2<p^{\prime}\leq n$ and $p_{0}^{\prime}=p^{\prime}/2$, and consider
$X=L^{(p,1)}$, $X_{0}=L^{(p_{0},1)}.$ Note that we can take $l_{1},l_{2}%
,p_{1},p_{2}$ satisfying $1<p_{1}<p<p_{2}<p_{0}$, $p\leq l$, $\frac
{p_{i}^{\prime}}{l_{i}^{\prime}}=\frac{p^{\prime}}{l^{\prime}}$, for $i=1,2$,
and $\frac{1}{p}=\frac{1-\theta}{p_{1}}+\frac{\theta}{p_{2}}$ for some
$\theta\in\left(  0,1\right)  $.

Considering $X_{1}=L^{(p_{1},1)}$ and $X_{2}=L^{(p_{2},1)}$, it follows that
$L^{(p,1)}=\left(  L^{(p_{1},1)},L^{(p_{1},1)}\right)  _{\theta,1}$. Moreover,
we have that%
\[
\left(  L^{(p_{0}^{\prime},\infty)},L^{\infty}\right)  _{\frac{p_{0}^{\prime}%
}{p_{1}^{\prime}},\infty}=L^{(p_{1}^{\prime},\infty)}%
\,\,\,\mbox{and}\,\,\,\left(  L^{(p_{0}^{\prime},\infty)},L^{\infty}\right)
_{\frac{p_{0}^{\prime}}{p_{2}^{\prime}},\infty}=L^{(p_{2}^{\prime},\infty)}.
\]
Thus, we obtain the properties \textbf{\textit{(i)-(iii)}} in Theorem
\ref{Teo:estim-bilin-cambio} for the spaces $X_{0},X_{1}$ and then the
bilinear estimate (\ref{eq:Estimativa bilinear}) with $E^{\prime}%
=\mathcal{M}_{\left(  p^{\prime},\infty\right)  }^{n}$, as desired.

\section{Bilinear estimate in Besov-$X$-Morrey spaces}

\label{sect:bilinear-estimate-in-besov-morrey}

In this part we combine some ideas from Morrey and Besov-type spaces in order
to obtain the bilinear estimate (\ref{eq:Estimativa bilinear}) in
Besov-Morrey-type spaces.

\begin{lemma}
\label{Lem: Otras bernstein} Let $X\in\mathcal{H}$ be such that $X^{\prime}%
\in\mathcal{H}$, and suppose that:

\begin{description}
\item[\textbf{(BM1)}] There exists $H_{0}\in\mathcal{H}$ such that $\left\Vert
fg\right\Vert _{X^{\prime}}\leq C\left\Vert f\right\Vert _{H_{0}}\left\Vert
g\right\Vert _{H_{0}}$. Moreover, $\left(  L^{\infty},X^{\prime}\right)
_{\frac{1}{2},d_{0}}\hookrightarrow H_{0}$ for some $1\leq d_{0}\leq\infty.$

\item[\textbf{(BM2)}] There exist $H_{1}$ and $H_{2}$ such that $\left\Vert
fg\right\Vert _{H_{2}}\leq C\left\Vert f\right\Vert _{H_{1}}\left\Vert
g\right\Vert _{H_{1}}$. Moreover, $\sigma_{X^{\prime}}=\theta\sigma_{H_{2}}$
and $\left(  L^{\infty},H_{2}\right)  _{\theta,d_{2}}\hookrightarrow
X^{\prime}$ with $\frac{1}{2}<\theta<1$ and $1\leq d_{2}\leq\infty$.
Additionally, suppose that $\left(  L^{\infty},X^{\prime}\right)  _{\frac
{1}{2\theta},d_{1}}\hookrightarrow H_{1}$ with $1\leq d_{1}\leq\infty$.
\end{description}

Then, we obtain the condition \textbf{(B5)} for the spaces $F=\mathcal{P}%
\mathcal{D}[X]^{l},G_{0}=\mathcal{M}\left[  H_{0}\right]  ^{2l^{\prime}}%
,G_{1}=\mathcal{M}\left[  H_{1}\right]  ^{2l_{2}}$ and $G_{2}=\mathcal{M}%
\left[  H_{2}\right]  ^{l_{2}}$, where $l$ satisfies $\frac{n}{l}+\sigma
_{X}\geq0$ and $l^{\prime}\geq\frac{1}{\theta}$, and $l_{2}=\theta l^{\prime}$.
\end{lemma}

\textbf{Proof. }Using \textbf{(BM1)} and Lemma \ref{Lem: producto en morrey},
we get
\[
\left\Vert fg\right\Vert _{\mathcal{M}\left[  X^{\prime}\right]  ^{l^{\prime}%
}}\leq C\left\Vert f\right\Vert _{\mathcal{M}\left[  H_{0}\right]
^{2l^{\prime}}}\left\Vert g\right\Vert _{\mathcal{M}\left[  H_{0}\right]
^{2l^{\prime}}}.
\]
Moreover, using again \textbf{(BM1)} and Lemma
\ref{Lem: Limitacion de convolucion norma Morrey Q, por norma Morrey N}, we
arrive at
\[
\left\Vert \Delta_{0}f\right\Vert _{\mathcal{M}\left[  H_{0}\right]
^{2l^{\prime}}}\leq C\left\Vert f\right\Vert _{\mathcal{M}\left[  X^{\prime
}\right]  ^{l^{\prime}}}.
\]
Now, condition \textbf{(BM2)} and Lemma \ref{Lem: producto en morrey} yield
\[
\left\Vert fg\right\Vert _{\mathcal{M}\left[  H_{2}\right]  ^{l_{2}}}\leq
C\left\Vert f\right\Vert _{\mathcal{M}\left[  H_{1}\right]  ^{2l_{2}}%
}\left\Vert g\right\Vert _{\mathcal{M}\left[  H_{1}\right]  ^{2l_{2}}}.
\]
In turn, using again \textbf{(BM2)} and Lemma
\ref{Lem: Limitacion de convolucion norma Morrey Q, por norma Morrey N}, it
follows that
\[
\left\Vert \Delta_{0}f\right\Vert _{\mathcal{M}\left[  X^{\prime}\right]
^{l^{\prime}}}\leq C\left\Vert f\right\Vert _{\mathcal{M}\left[  H_{2}\right]
^{l_{2}}}\text{ and }\left\Vert \Delta_{0}f\right\Vert _{\mathcal{M}\left[
H_{1}\right]  ^{2l_{2}}}\leq C\left\Vert f\right\Vert _{\mathcal{M}\left[
X^{\prime}\right]  ^{l^{\prime}}}.
\]
\fin

In the sequel we prove the bilinear estimate in the case of Besov-Morrey-type spaces.

\begin{theorem}
\label{Teo-Est-BM-Spaces}Let $\frac{n}{2}<l^{\prime}<n$, $s=1-\frac
{n}{l^{\prime}}$ and let $X\in\mathcal{H}$ be such that $X^{\prime}%
\in\mathcal{H}.$ Moreover, assume \textbf{(BM1)} and \textbf{(BM2)}. Then, we
obtain the conditions \textbf{(B1)-(B6)} for the space $F=\mathcal{PD}\left[
X\right]  ^{l}$. As a consequence, the bilinear estimate
\eqref{eq:Estimativa bilinear} with $E^{\prime}=\dot{B}\mathcal{M}\left[
X^{\prime}\right]  _{\infty}^{l^{\prime},\frac{n}{l^{\prime}}-1}=\dot
{B}[\mathcal{M}[X^{\prime}]^{l^{\prime}}]_{\infty}^{\frac{n}{l^{\prime}}-1}$
holds true.
\end{theorem}

\textbf{Proof. }The property \textbf{(B1)} follows by Lemma
\ref{Lem: espacio de bloques pertenece a la clase G}. The conditions
\textbf{(B2)}, \textbf{(B3)} and \textbf{(B4)} follow respectively from
Lemma \ref{Cor: limitacion norma Linft deltaf con norma de bloques.},
Corollary \ref{Cor: convolucion en Bloques} and Theorem
\ref{Teo: Dualidad bloques-Morrey}. Also, we get \textbf{(B5)} by Lemma
\ref{Lem: Otras bernstein}. Finally, \textbf{(B6)} follows by taking $\theta$
in Lemma \ref{Lem: Otras bernstein} close enough to $1/2$ such that
$2l_{2}=2\theta l^{\prime}<n$. Finally, we conclude by using Theorem
\ref{Teo. Estimativa bilinear espacios tipo besov} with $E^{\prime}=\dot
{B}[\mathcal{M}[X^{\prime}]^{l^{\prime}}]_{\infty}^{\frac{n}{l^{\prime}}-1}$,
$F^{\prime}=\mathcal{M}[X^{\prime}]^{l^{\prime}}$ and $\sigma_{F^{\prime}%
}=-\frac{n}{l^{\prime}}.$\fin

\subsection{Applications (bilinear estimate in Besov-Lorentz-Morrey spaces)}

We can take $X=L^{(p,d)}$ with $2<p^{\prime}\leq l^{\prime}<n$, $\frac{n}%
{2}<l^{\prime}$ and $1\leq d\leq\infty$ in order to obtain the bilinear
estimate (\ref{eq:Estimativa bilinear}) with the Besov-Lorentz-Morrey
$E^{\prime}=\dot{B}\mathcal{M}_{(p^{\prime},d^{\prime}),\infty}^{l^{\prime
},\frac{n}{l^{\prime}}-1}=\dot{B}\mathcal{M}\left[  L^{(p^{\prime},d^{\prime
})}\right]  _{\infty}^{l^{\prime},\frac{n}{l^{\prime}}-1}$, which is a
generalization of the Besov-Morrey space $\dot{B}\mathcal{M}_{p^{\prime
},\infty}^{l^{\prime},\frac{n}{l^{\prime}}-1}=\dot{B}\mathcal{M}\left[
L^{p^{\prime}}\right]  _{\infty}^{l^{\prime},\frac{n}{l^{\prime}}-1}$ defined
in \cite{KoYa}. This result recovers and extends the bilinear estimate
obtained in \cite{LucJhe2}. In the case $d=p,$ we obtain $X=L^{(p,p)}=L^{p}$
and (\ref{eq:Estimativa bilinear}) in the Besov-Morrey space $E^{\prime}%
=\dot{B}\mathcal{M}_{p^{\prime},\infty}^{l^{\prime},\frac{n}{l^{\prime}}-1}.$

\section{Uniqueness for Navier-Stokes equations} \label{Sec:unicidad}

We finish the study of Navier-Stokes equations by proving a uniqueness result, regardless of the size of the initial data and solutions. For that, we employ the bilinear estimate \eqref{eq:Estimativa bilinear} and consider some further standard properties for our framework.

\begin{theorem} \label{Teo: condiciones para Unicidad} Let $E\in\mathcal{G}$
	be a Banach space with $\sigma_{E^{\prime}}=-1$, such that the bilinear
	estimate \eqref{eq:Estimativa bilinear} holds true. Moreover, assume
	that there exist Banach spaces $E_{3},E_{4}\in\mathcal{G}$ verifying
	the following conditions:
	\begin{enumerate}
		\item[\textbf{(U1)}]  $-1<\sigma_{E_{3}^{\prime}}<1$;
		\item[\textbf{(U2)}]  $\left\Vert U\left(s\right)f\right\Vert _{E_{3}^{\prime}}\leq C\left\Vert f\right\Vert _{E_{3}^{\prime}}$;
		\item[\textbf{(U3)}] $\left\Vert U\left(s\right)f\right\Vert _{E_{3}^{\prime}}\leq Cs^{\frac{1}{2}\left(\sigma_{E^{\prime}}-\sigma_{E_{3}^{\prime}}\right)}\left\Vert f\right\Vert _{E^{\prime}}$;
		\item[\textbf{(U4)}]  $\left\Vert f\cdot g\right\Vert _{E_{4}^{\prime}}\leq C\left\Vert f\right\Vert _{E^{\prime}}\left\Vert f\right\Vert _{E_{3}^{\prime}}$;
		\item[\textbf{(U5)}]  $\left\Vert \nabla U\left(s\right)\mathbb{P}f\right\Vert _{E^{\prime}}\leq Cs^{\frac{1}{2}\left(\sigma_{E^{\prime}_{4}}-\sigma_{E^{\prime}}\right)-\frac{1}{2}}\left\Vert f\right\Vert _{E_{4}^{\prime}}$.
	\end{enumerate}
		Then, if $u$ and $v$ are two mild solutions of \eqref{eq:N-S-equation1}-\eqref{eq:N-S-equation3}
		in $C\left([0,T);E^{\prime}\right)$ with the same initial data $u_{0}\in\tilde{E^{\prime}}$,
		it follows that $u(\cdot,t)=v(\cdot,t)$ in $E^{\prime}$ for all $t\in[0,T)$.
		Here, $\tilde{E^{\prime}}$ stands for the maximal closed subspace
		of $E^{\prime}$ in which the heat semigroup $\{U(t)\}_{t\geq0}$
		is strongly continuous.
	
\end{theorem}

\textbf{Proof. }
The proof is based on an argument due to Meyer \cite[p. 188]{Meyer}. First, we prove that there exists $0<T_{1}<T$ such that $u(\cdot,t)=v(\cdot,t)$ in $E^{\prime}$ for all $t\in\left[0,T_{1}\right)$. Denoting $w=u-v,w_{1}=U(t)u_{0}-u$ and $w_{2}=U(t)u_{0}-v$, we have
that
\[
\begin{aligned}u\otimes u-v\otimes v & =w\otimes u+v\otimes w\\
	& =w\otimes U(t)u_{0}+U(t)u_{0}\otimes w-w\otimes w_{1}-w_{2}\otimes w.
\end{aligned}
\]

We can estimate the difference $w(t)$ in $E^{\prime}$ as follows
\[
\begin{aligned}\|w(t)\|_{E^{\prime}} & =\left\Vert \int_{0}^{t}\nabla U(t-s)\mathbb{P}(u\otimes u-v\otimes v)ds\right\Vert _{E^{\prime}}\\
	& \leq\left\Vert \int_{0}^{t}\nabla U(t-s)\mathbb{P}\left(w\otimes w_{1}+w_{2}\otimes w\right)ds\right\Vert _{E^{\prime}}\\
	& \qquad+\left\Vert \int_{0}^{t}\nabla U(t-s)\mathbb{P}\left(w\otimes U(s)u_{0}+U(s)u_{0}\otimes w\right)ds\right\Vert _{E^{\prime}}\\
	& :=I_{1}(t)+I_{2}(t).
\end{aligned}
\]

From the bilinear estimate \eqref{eq:Estimativa bilinear}, we conclude
that
\[
I_{1}(t)\leq K\sup_{0<t<T_{1}}\|w\|_{E^{\prime}}\left(\sup_{0<t<T_{1}}\left\Vert w_{1}\right\Vert _{E^{\prime}}+\sup_{0<t<T_{1}}\left\Vert w_{2}\right\Vert _{E^{\prime}}\right).
\]

Now we estimate $I_{2}(t)$. First we use \textbf{(U5)} and then \textbf{(U4)}
to obtain

\[
\begin{aligned}I_{2}(t) & \leq C\int_{0}^{t}\left\Vert \nabla U(t-s)\mathbb{P}\left(w\otimes U(s)u_{0}+U(s)u_{0}\otimes w\right)\right\Vert _{E^{\prime}}ds\\
	& \le C\int_{0}^{t}(t-s)^{\frac{1}{2}\left(\sigma_{E_{4}^{\prime}}-\sigma_{E^{\prime}}\right)-\frac{1}{2}}\left\Vert w\otimes U(s)u_{0}+U(s)u_{0}\otimes w\right\Vert _{E_{4}^{\prime}}ds\\
	& \le C\int_{0}^{t}(t-s)^{\frac{1}{2}\left(\sigma_{E_{4}^{\prime}}-\sigma_{E^{\prime}}\right)-\frac{1}{2}}\left\Vert w\otimes U(s)u_{0}\right\Vert _{E_{4}^{\prime}}ds\\
	& \le C\int_{0}^{t}(t-s)^{\frac{1}{2}\left(\sigma_{E_{4}^{\prime}}-\sigma_{E^{\prime}}\right)-\frac{1}{2}}\left\Vert w\right\Vert _{E^{\prime}}\left\Vert U(s)u_{0}\right\Vert _{E_{3}^{\prime}}ds\\
	& \le C\int_{0}^{t}(t-s)^{\frac{1}{2}\left(\sigma_{E_{4}^{\prime}}-\sigma_{E^{\prime}}\right)-\frac{1}{2}}\left\Vert w\right\Vert _{E^{\prime}}s^{-\frac{1}{2}\left(\sigma_{E_{3}^{\prime}}-\sigma_{E^{\prime}}\right)}s^{\frac{1}{2}\left(\sigma_{E_{3}^{\prime}}-\sigma_{E^{\prime}}\right)}\left\Vert U(s)u_{0}\right\Vert _{E_{3}^{\prime}}ds\\
	& \le C\sup_{0<t<T_{1}}\left\Vert w\right\Vert _{E^{\prime}}\sup_{0<t<T_{1}}t^{\frac{1}{2}\left(\sigma_{E_{3}^{\prime}}-\sigma_{E^{\prime}}\right)}\left\Vert U(t)u_{0}\right\Vert _{E_{3}^{\prime}}\int_{0}^{t}(t-s)^{\frac{1}{2}\left(\sigma_{E_{4}^{\prime}}-\sigma_{E^{\prime}}\right)-\frac{1}{2}}s^{-\frac{1}{2}\left(\sigma_{E_{3}^{\prime}}-\sigma_{E^{\prime}}\right)}ds\\
	& \le C\sup_{0<t<T_{1}}\left\Vert w\right\Vert _{E^{\prime}}\sup_{0<t<T_{1}}t^{\frac{1}{2}\left(\sigma_{E_{3}^{\prime}}-\sigma_{E^{\prime}}\right)}\left\Vert U(t)u_{0}\right\Vert _{E_{3}^{\prime}}\int_{0}^{1}(1-\rho)^{\frac{1}{2}\left(\sigma_{E_{4}^{\prime}}-\sigma_{E^{\prime}}\right)-\frac{1}{2}}\rho^{-\frac{1}{2}\left(\sigma_{E_{3}^{\prime}}-\sigma_{E^{\prime}}\right)}d\rho\\
	& \le C\sup_{0<t<T_{1}}\left\Vert w\right\Vert _{E^{\prime}}\sup_{0<t<T_{1}}t^{\frac{1}{2}\left(\sigma_{E_{3}^{\prime}}-\sigma_{E^{\prime}}\right)}\left\Vert U(t)u_{0}\right\Vert _{E_{3}^{\prime}}.
\end{aligned}
\]
In the above computations, the finiteness of the last integral is guaranteed because condition \textbf{(U4) }implies that $\sigma_{E_{4}^{\prime}}=\sigma_{E^{\prime}}+\sigma_{E_{3}^{\prime}}$ (see Remark \ref{remark:scaling-condit-to-holder-type-inequa}) and then condition \textbf{(U1)} implies that $-1<\frac{1}{2}\left(\sigma_{E_{4}^{\prime}}-\sigma_{E^{\prime}}\right)-\frac{1}{2}<0$ and $-1<-\frac{1}{2}\left(\sigma_{E_{3}^{\prime}}-\sigma_{E^{\prime}}\right)<0$. With all the above developments, we arrive at

\[
\sup_{0<t<T_{1}}\|w(t)\|_{E^{\prime}}\leq CZ\left(T_{1}\right)\sup_{0<t<T_{1}}\|w\|_{E^{\prime}},
\]
where
\[
Z\left(T_{1}\right)=\sup_{0<t<T_{1}}\left\Vert w_{1}\right\Vert _{E^{\prime}}+\sup_{0<t<T_{1}}\left\Vert w_{2}\right\Vert _{E^{\prime}}+\sup_{0<t<T_{1}}t^{\frac{1}{2}\left(\sigma_{E_{3}^{\prime}}-\sigma_{E^{\prime}}\right)}\left\Vert U(t)u_{0}\right\Vert _{E_{3}^{\prime}}.
\]

Note that from the conditions for $u,v$ and $u_{0},$ we have that $U(t)u_{0},u,v\rightarrow u_{0}$ in $E^{\prime}$ as $t\rightarrow0^{+}$, and then
\[
\lim_{t\rightarrow0^{+}}\left\Vert w_{1}\right\Vert _{E^{\prime}}=\lim_{t\rightarrow0^{+}}\left\Vert w_{2}\right\Vert _{E^{\prime}}=0.
\]
Now, we prove that
\[
\limsup_{t\rightarrow0^{+}}t^{\frac{1}{2}\left(\sigma_{E_{3}^{\prime}}-\sigma_{E^{\prime}}\right)}\left\Vert U(t)u_{0}\right\Vert _{E_{3}^{\prime}}=0.
\]
To do so, we define $u_{0k}=U\left(\frac{1}{k}\right)u_{0}$ for all
$k\in\mathbb{N}$. It follows by \textbf{(U3)} that $u_{0k}\in E_{3}^{\prime}$,
and from \textbf{(U2) }that $\left\Vert U(t)u_{0k}\right\Vert _{E_{3}^{\prime}}\le C\left\Vert u_{0k}\right\Vert _{E_{3}^{\prime}}$.
Using again that $u_{0k}\rightarrow u_{0}$ in $E^{\prime}$ as $k\rightarrow\infty$, we can proceed as follows:
\[
\begin{aligned} & \limsup_{t\rightarrow0^{+}}t^{\frac{1}{2}\left(\sigma_{E_{3}^{\prime}}-\sigma_{E^{\prime}}\right)}\left\Vert U(t)u_{0}\right\Vert _{E_{3}^{\prime}}\\
	& \qquad\qquad\leq\limsup_{t\rightarrow0^{+}}t^{\frac{1}{2}\left(\sigma_{E_{3}^{\prime}}-\sigma_{E^{\prime}}\right)}\left\Vert U(t)\left(u_{0}-u_{0k}\right)\right\Vert _{E_{3}^{\prime}}+\limsup_{t\rightarrow0^{+}}t^{\frac{1}{2}\left(\sigma_{E_{3}^{\prime}}-\sigma_{E^{\prime}}\right)}\left\Vert U(t)u_{0k}\right\Vert _{E_{3}^{\prime}}\\
	& \qquad\qquad\leq C\left\Vert u_{0}-u_{0k}\right\Vert _{E^{\prime}}+C\left\Vert u_{0k}\right\Vert _{E_{3}^{\prime}}\limsup_{t\rightarrow0^{+}}t^{\frac{1}{2}\left(\sigma_{E_{3}^{\prime}}-\sigma_{E^{\prime}}\right)}\\
	& \qquad\qquad\leq C\left\Vert u_{0}-u_{0k}\right\Vert _{E^{\prime}}\rightarrow0,\text{ as }k\rightarrow\infty.
\end{aligned}
\]
Consequently, there exists $T_{1}>0$ small enough such that $CZ\left(T_{1}\right)<1$
and then $w(t)=0$ for all $t\in\left[0,T_{1}\right)$. The remainder
of the proof is to show that in fact $T_{1}\in(0,T]$ can be
arbitrary. Define
\[
T_{*}=\sup\left\{ \tilde{T};0<\tilde{T}<T,u(t)=v(t)\text{ in }E^{\prime}\text{ for all }t\in[0,\tilde{T})\right\} .
\]
If $T_{*}=T$ we are done; otherwise, we have that $u(t)=v(t)$
for $t\in\left[0,T_{*}\right)$ which implies that $u\left(T_{*}\right)=v\left(T_{*}\right)$ due to the time-continuity of $u$ and $v$. It follows from the
first part of the proof that there exists $T_{2}>0$ small enough
such that $u(t)=v(t)$ for $t\in\left[T_{*},T_{*}+T_{2}\right)$.
Therefore $u(t)=v(t)$ for $t\in\left[0,T_{*}+T_{2}\right)$ which contradicts the definition of $T_{*}$.
\fin

\subsection{Applications}
In this section we present some examples of spaces $E^{\prime}$ where
the bilinear estimate (\ref{eq:Estimativa bilinear}) holds true and then the uniqueness follows in a direct way by applying Theorem \ref{Teo: condiciones para Unicidad}.
\begin{itemize}
	\item Consider the space $E^{\prime}=L^{n,\infty}$. In this case the uniqueness
	result follows by choosing $E_{3}^{\prime}=L^{q,\infty}$ and $E_{4}^{\prime}=L^{r,\infty}$
	where $n<q<\infty$ and $\frac{1}{r}=\frac{1}{n}+\frac{1}{q}.$
	\item Similarly, if we consider $E^{\prime}=W\dot{K}_{n,\infty}^{0}$, the
	uniqueness result follows by choosing $E_{3}^{\prime}=W\dot{K}_{q,\infty}^{0}$
	and $E_{4}^{\prime}=W\dot{K}_{r,\infty}^{0}$ where $n<q<\infty$
	and $\frac{1}{r}=\frac{1}{n}+\frac{1}{q}.$
	\item For the space $E^{\prime}=\dot{B}\left[W\dot{K}_{p,\infty}^{0}\right]_{\infty}^{\frac{n}{p}-1},$ where
	$\frac{n}{2}<p<n$, we can consider $E_{3}^{\prime}=\dot{B}\left[W\dot{K}_{p,\infty}^{0}\right]_{\infty}^{\frac{n}{p}-1+\rho}$
	and $E_{4}^{\prime}=\dot{B}\left[W\dot{K}_{p,\infty}^{0}\right]_{\infty}^{\frac{n}{p}-2+\rho}$
	with $\rho>0$ small enough. All the conditions in Theorem \ref{Teo: condiciones para Unicidad}
	follow from Lemma \ref{Prop: Estimativa para el nucleo del calor en Besov-Lorentz-Morrey y besov-Lorentz-Bolcos} and Lemma \ref{Lem: Estimativa producto}.
\end{itemize}


\section{Other nonlinearities and PDEs}

\label{section:related-models}

The present section is devoted to give some other examples of nonlinearities
and PDEs to which the theory developed in the previous sections can be
adapted. \vspace{0.3cm}

\noindent\textit{\underline{Reaction-diffusion system}}: we can consider the
bilinear form
\begin{equation}
B_{1,i}(u,v)(t)=\int_{0}^{t}U(t-s)\left\langle A_{i}\cdot
u(s),v(s)\right\rangle \ ds\ \text{ for }i=1,\ldots,n,
\label{eq:bilin-form-reac-diff}%
\end{equation}
which is related to the system of reaction-diffusion equations
\begin{equation}
u_{it}=\Delta u_{i}+\left\langle A_{i}\cdot u,u\right\rangle ,\text{ for
}i=1,\ldots,n, \label{eq:reac-difus-syst}%
\end{equation}
where $A_{i}$ denotes constant $n\times n$-matrices and $u=(u_{1},\ldots
,u_{n}):\mathbb{R}^{n}\times\lbrack0,\infty)\rightarrow\mathbb{R}^{n}$.

\vspace{0.3cm}

\noindent\textit{\underline{Nonlocal advection-diffusion system}}: consider
the bilinear form
\begin{equation}
B_{2}(u,v)(t)=-\int_{0}^{t}\nabla U(t-s)\left(  u(s)\otimes S(v)\right)  \ ds,
\label{eq:bilin-nonlocal-form}%
\end{equation}
which is related to the nonlocal parabolic equation
\begin{equation}
u_{t}=\Delta u+\nabla\cdot\left(  u\otimes S(u)\right)  ,
\label{eq:non-local-model}%
\end{equation}
where $S(u)$ is a vector field given by the convolution of $u$ with a vector
function whose components are $b$-homogeneous singular kernels $\psi$, that
is,
\begin{equation}
S(u)=\psi\ast u, \label{aux-nonlocal-1}%
\end{equation}
where $\psi(\lambda x)=\lambda^{b}\psi(x)$ and $|\psi(x)|\leq|x|^{b}$ for some
$b<0$. In this case we are concerned with values $b=-n+\delta$ for
$0<\delta<n$. For the reader convenience, examples of such nonlocal operators
are given below. \vspace{0.3cm}

\begin{itemize}
\item \textit{Models appearing in nonequilibrium statistical mechanics}:
$\psi(x)=\nabla G_{n}(x)$ ($\delta=1$), where $G_{n}$ is the fundamental
solution of the Laplacian in $\mathbb{R}^{n}$, see \cite{Biler}.
\vspace{0.3cm}

\item \noindent\textit{The Debye system from the theory of electrolytes}
presents the form
\begin{equation}%
\begin{array}
[c]{ccc}%
u_{1t} & = & \Delta u_{1}-\nabla\cdot\left(  u_{1}\psi\right) \\
u_{2t} & = & \Delta u_{2}+\nabla\cdot\left(  u_{2}\psi\right)
\end{array}
, \label{aux-Debye-1}%
\end{equation}
where $\psi=\nabla G_{n}\ast(u_{1}-u_{2})$ ($\delta=1$). \vspace{0.3cm}

\item \noindent\textit{Aggregation diffusion equations for multi-species} is a
system with the form
\[
\frac{\partial u_{i}}{\partial t}=\Delta u_{i}-\nabla\left[  u_{i}\nabla
\sum_{j=1}^{n}h_{ij}K\ast u_{j}\right]  ,\ i=1,2,\ldots,n,
\]
where $u_{i}(x,t)$ stands for the density of the population $i$ and $h_{ij}$
are constants denoting attraction ($h_{ij}>0$) or repulsion ($h_{ij}<0$) of
the population $i$ to population $j$.
\end{itemize}

Recall the nonlocal operator%

\begin{equation}
(-\Delta)^{-\delta/2}(f)=c_{\delta}\int_{\mathbb{R}^{n}}|x-y|^{-n+\delta
}f(y)\ dy, \label{eq:laplac-inv-frac}%
\end{equation}
where $c_{\delta}>0$ is a constant and $0<\delta<n$. For the bilinear forms
$B_{1}$ and $B_{2}$ related to the reaction-diffusion equation and the
nonlocal parabolic equation, we have an analogous result to Theorem
\ref{Teo:estim-bilin-cambio}. Let us state the result for the former case.

\begin{theorem}
\label{Theorem:bilinear-estim-B1} Let $1\leq l<l_{0}<\infty$ and $X_{0}%
,X\in\mathcal{H}$ be such that $\sigma_{X}<\sigma_{X_{0}}$, $l_{0}^{\prime
}\sigma_{X_{0}^{\prime}}=l^{\prime}\sigma_{X^{\prime}}$ and $X_{0}^{\prime
},X^{\prime}\in\mathcal{H}$. Assume that the conditions \textbf{\textit{(i)}}
and \textbf{\textit{(ii)}} in Theorem \ref{Teo:estim-bilin-cambio} and the
H\"{o}lder-type inequality $\Vert fg\Vert_{X_{0}^{\prime}}\leq C\Vert
f\Vert_{X^{\prime}}\Vert g\Vert_{X^{\prime}}$ hold true. If $n>4$, then we
have the bilinear estimate
\begin{equation}
\sup_{t>0}\Vert B_{1}(u,v)(t)\Vert_{\mathcal{M}[X^{\prime}]^{n/2}}\leq
K\sup_{t>0}\Vert u(t)\Vert_{\mathcal{M}[X^{\prime}]^{n/2}}\sup_{t>0}\Vert
v(t)\Vert_{\mathcal{M}[X^{\prime}]^{n/2}}. \label{eq:estim-bilin-react-difuss}%
\end{equation}

\end{theorem}

\noindent\textbf{Proof.} We only give some steps. First, from items
\textbf{\textit{(i)}} and \textbf{\textit{(ii)}} in Theorem
\ref{Teo:estim-bilin-cambio} we can perform the same arguments as in Lemma
\ref{Lem: condiciones para estimativa integral} in order to get
\begin{equation}
\int_{0}^{\infty}s^{\frac{1}{2}\left(  \frac{n}{l}-\frac{n}{l_{0}}\right)
-1}\Vert U(s)f\Vert_{\mathcal{PD}[X_{0}]^{l_{0}}}\ ds\leq C\Vert
f\Vert_{\mathcal{PD}[X]^{l}}. \label{eq:H4-for-eact-difuss}%
\end{equation}
Next, the H\"{o}lder-type inequality and condition $\frac{n}{l}-\frac{n}%
{l_{0}}-2=0$ imply that $l^{\prime}=n/2$ and $l_{0}^{\prime}=n/4$. Thus, one
can obtain \eqref{eq:estim-bilin-react-difuss} by proceeding as in Lemma
\ref{Lem: Estimativa integral espacio dual} and Theorem
\ref{Teo: Teo main. Estimativa integral en el dual.}.\fin

Next, we perform the bilinear estimate for the form $B_{2}(\cdot,\cdot)$
defined in (\ref{eq:bilin-nonlocal-form}). This part employs some ideas of
\cite[Chapter 7]{adams-libro} in order to control the nonlocal operator
(\ref{aux-nonlocal-1}). Similar estimates for fractional operators can be
found in \cite{Lemarie2,Hatano}.

\begin{definition}
\label{definition:norm-power-rescaled}Given $\kappa>1$, we define the family
$\mathcal{H}_{\kappa}$ of all spaces $X\in\mathcal{H}$ such that the function
\[
f\mapsto\Vert|f|^{\kappa}\Vert_{X}^{1/\kappa}%
\]
satisfies the triangle inequality and the norm of $X$ has the following
monotonicity property:
\begin{equation}
|f(x)|\leq|g(x)|,\ \text{ for a.e. }x\in\mathbb{R}^{n}\Rightarrow\Vert
f\Vert_{X}\leq\Vert g\Vert_{X}. \label{eq:monotonicity-porperty-norm}%
\end{equation}

\end{definition}

For $X\in\mathcal{H}_{\kappa}$, it is easy to see that $X_{\kappa}=\{f\in
L_{loc}^{1}:\ |f|^{\kappa}\in X\}$ is a Banach space and satisfies
\textbf{(M1)}, \textbf{(M2)}, \textbf{(M3)}, \textbf{(M5)} and \textbf{(M6)}.
Now, let us recall the centered maximal operators
\begin{equation}
M_{\alpha}f(x)=\sup_{R>0}\frac{1}{R^{n-\alpha}}\int_{D(x,R)}|f(y)|\ dy,\text{
for }f\in L_{loc}^{1}\text{,} \label{eq:maximal-operators}%
\end{equation}
where $0\leq\alpha\leq n$. We start with the following preliminary result.

\begin{lemma}
\label{lema:Maximal-H-L} Let $X\in\mathcal{H}$ and $l>1$ be such that
$\frac{n}{l}+\sigma_{X}\leq0,$ the monotonicity property
\eqref{eq:monotonicity-porperty-norm} is satisfied, and the maximal operator
$M_{0}:X\rightarrow X$ is bounded. Then, $M_{0}:\mathcal{M}[X]^{l}%
\rightarrow\mathcal{M}[X]^{l}$ is bounded.
\end{lemma}

\textbf{Proof.} Take $x_{0}\in\mathbb{R}^{n}$ and $R_{0}>0,$ and consider
$f(x)=f(x)1_{D(x_{0},2R_{0})}(x)+f(x)1_{D(x_{0},2R_{0})^{c}}(x)=f_{1}%
(x)+f_{2}(x)$. Since $M_{0}$ is subadditive, i.e.
\[
M_{0}f\leq M_{0}f_{1}+M_{0}f_{2},
\]
we can estimate
\begin{align*}
\Vert M_{0}f_{1}\cdot1_{D(x_{0},R_{0})}\Vert_{X}\leq &  \ \Vert M_{0}%
f_{1}\Vert_{X}\\
\leq &  \ C\Vert f1_{D(x_{0},2R_{0})}\Vert_{X}\\
\leq &  \ CR_{0}^{-\frac{n}{l}-\sigma_{X}}\Vert f\Vert_{\mathcal{M}[X]^{l}},
\end{align*}
where we have used the boundedness of $M_{0}$ on $X$. On the other hand, using
Proposition \ref{prop:conseq-condition-morrey-space}, for $x\in D(x_{0}%
,R_{0})$ and $r>0$, we have that
\begin{align*}
\frac{1}{r^{n}}\int_{D(x,r)}|f_{2}(y)|\ dy=  &  \ \frac{1}{r^{n}}%
\int_{D(x,r)\cap D(x_{0},2R_{0})^{c}}|f(y)|\ dy1_{(R_{0},+\infty)}(r)\\
\leq &  \ Cr^{\sigma_{X}}\Vert f1_{D(x,r)}\Vert_{X}1_{(R_{0},+\infty)}(r)\\
\leq &  \ CR_{0}^{-\frac{n}{l}}\Vert f\Vert_{\mathcal{M}[X]^{l}}%
1_{(R_{0},+\infty)}(r).
\end{align*}
This last inequality leads us to
\begin{align*}
\Vert M_{0}f_{2}\cdot1_{D(x_{0},R_{0})}\Vert_{X}\leq &  \ CR_{0}^{-\frac{n}%
{l}}\Vert f\Vert_{\mathcal{M}[X]^{l}}\Vert1_{D(x_{0},R_{0})}\Vert_{X}\\
\leq &  \ CR_{0}^{-\frac{n}{l}-\sigma_{X}}\Vert f\Vert_{\mathcal{M}[X]^{l}}.
\end{align*}
Combining all the above inequalities, we are done.\fin

The $X$-Morrey spaces provide a suitable framework to control the maximal
operators $M_{\alpha},$ for $0<\alpha\leq n$. In fact, using again Proposition
\ref{prop:conseq-condition-morrey-space}, note that
\begin{align*}
\frac{1}{R^{n-\alpha}}\int_{D(x,R)}|f(y)|\ dy\leq &  \ \frac{C}{R^{-\sigma
_{X}-\alpha}}\Vert f1_{D(x,R)}\Vert_{X}\\
\leq &  \ C\Vert f\Vert_{\mathcal{M}[X]^{n/\alpha}},
\end{align*}
for all $f\in\mathcal{M}[X]^{n/\alpha}$. It follows that
\begin{equation}
M_{\alpha}f(x)\leq C\Vert f\Vert_{\mathcal{M}[X]^{n/\alpha}},
\label{eq:maximal-op-morrey-spaces}%
\end{equation}
provided that $\alpha+\sigma_{X}\leq0$.

With the above developments in hand, we are in a position to show the
boundedness of the operator $(-\Delta)^{-\delta/2}$ in the framework of
$X$-Morrey spaces.

\begin{proposition}
\label{prop:boundedness_fractional_operator} Consider $X\in\mathcal{H}%
_{\frac{\alpha}{\alpha-\delta}}$ with $0<\delta<\alpha\leq n$ and such that
$\alpha+\sigma_{X}\leq0$, and suppose that the maximal operator $M_{0}$ is
bounded in $X$. Then, the fractional operator $(-\Delta)^{-\delta
/2}:\mathcal{M}[X]^{\frac{n}{\alpha}}\rightarrow\mathcal{M}\left[
X_{\frac{\alpha}{\alpha-\delta}}\right]  ^{\frac{n}{\alpha-\delta}}$ is bounded.
\end{proposition}

\textbf{Proof.} We employ the so-called Hedberg trick (see \cite[Thm.
7.1]{adams-libro}) to obtain
\begin{align*}
|(-\Delta)^{-\delta/2}f(x)|\leq &  \ C(M_{\alpha}f(x))^{\frac{\delta}{\alpha}%
}(M_{0}f(x))^{1-\frac{\delta}{\alpha}}\\
\leq &  \ C(\Vert f\Vert_{\mathcal{M}[X]^{n/\alpha}})^{\frac{\delta}{\alpha}%
}(M_{0}f(x))^{\frac{\alpha-\delta}{\alpha}}.
\end{align*}
Multiplying by $1_{D(x,R)}$, using \textbf{(M2)}, and taking the norm of
$X_{\frac{\alpha}{\alpha-\delta}}$, we arrive at
\[
\Vert(-\Delta)^{-\delta/2}f\cdot1_{D(x,R)}\Vert_{X_{\frac{\alpha}%
{\alpha-\delta}}}\leq\ C(\Vert f\Vert_{\mathcal{M}[X]^{n/\alpha}}%
)^{\frac{\delta}{\alpha}}\Vert M_{0}f\cdot1_{D(x,R)}\Vert_{X}^{\frac
{\alpha-\delta}{\alpha}}.
\]
Finally, multiplying the last inequality by $R^{(\alpha-\delta)(1+\sigma
_{X}/\alpha)}$, taking the supremum, and applying Lemma \ref{lema:Maximal-H-L}%
, we conclude the proof. \fin

\begin{remark}
\label{remark:Alternative_condition-fractional-operator} The previous result
also works well if, instead of assuming that $X\in\mathcal{H}_{\frac{\alpha
}{\alpha-\delta}}$, we assume that there exists $Y\in\mathcal{H}$ satisfying
the compatibility condition%
\[
\Vert f^{\frac{\alpha-\delta}{\alpha}}\Vert_{Y}\leq C\Vert f\Vert_{X}%
^{\frac{\alpha-\delta}{\alpha}},
\]
for all $f\in X$ such that $f\geq0.$ In this case, we have that $(-\Delta
)^{-\delta/2}:\mathcal{M}[X]^{\frac{n}{\alpha}}\rightarrow\mathcal{M}\left[
Y\right]  ^{\frac{n}{\alpha-\delta}}$ is bounded.
\end{remark}

Finally, we give the statement of the bilinear estimate for $B_{2}$ in the
framework of $X$-Morrey spaces.

\begin{theorem}
\label{Theorem:bilinear-estim-B2} Let $n>\delta+2,$ $1\leq l<l_{0}<\infty$ and
let $X_{0},X\in\mathcal{H}$ be such that $\frac{n}{l}-\frac{n}{l_{0}}=1$,
$X_{0}^{\prime},X^{\prime}\in\mathcal{H}$, $\sigma_{X}<\sigma_{X_{0}}$, and
$l_{0}^{\prime}\sigma_{X_{0}^{\prime}}=l^{\prime}\sigma_{X^{\prime}}$. Assume
the conditions \textbf{\textit{(i)}} and \textbf{\textit{(ii)}} in Theorem
\ref{Teo:estim-bilin-cambio}. Moreover, assume that the H\"{o}lder-type
inequality $\Vert f\cdot g\Vert_{X_{0}^{\prime}}\leq C\Vert f\Vert_{X^{\prime
}}\Vert|g|^{\kappa}\Vert_{X^{\prime}}^{\frac{1}{\kappa}}$, holds true for
$\kappa=\frac{\alpha}{\alpha-\delta}$ with $0<\delta<\alpha\leq n$. Then,
\begin{equation}
\sup_{t>0}\Vert B_{2}(u,v)(t)\Vert_{\mathcal{M}[X^{\prime}]^{l^{\prime}}}\leq
K\sup_{t>0}\Vert u(t)\Vert_{\mathcal{M}[X^{\prime}]^{l^{\prime}}}\sup
_{t>0}\Vert S(v)(t)\Vert_{\mathcal{M}[X_{\kappa}^{\prime}]^{l_{1}}},
\label{eq:estim-bilin-advect-difuss_1}%
\end{equation}
for $\frac{n}{l_{1}}+\frac{\sigma_{X^{\prime}}}{\kappa}\leq0$. In addition,
suppose that $X^{\prime}$ satisfies the hypotheses of Proposition
\ref{prop:boundedness_fractional_operator} and $l_{1}=\frac{n}{\alpha-\delta}%
$, then we obtain $\alpha=\delta+1$ and, in particular, by taking $l^{\prime
}=\frac{n}{\delta+1}$ we have the bilinear estimate
\begin{equation}
\sup_{t>0}\Vert B_{2}(u,v)(t)\Vert_{\mathcal{M}[X^{\prime}]^{n/(\delta+1)}%
}\leq K\sup_{t>0}\Vert u(t)\Vert_{\mathcal{M}[X^{\prime}]^{n/(\delta+1)}}%
\sup_{t>0}\Vert v(t)\Vert_{\mathcal{M}[X^{\prime}]^{n/(\delta+1)}}.
\label{eq:estim-bilin-advect-difuss_2}%
\end{equation}

\end{theorem}

\textbf{Proof.} It was already showed that conditions \textbf{\textit{(i)}}
and \textbf{\textit{(ii)}} in Theorem \ref{Teo:estim-bilin-cambio} imply
\textbf{(H4)}. So, Lemma \ref{Lem: Estimativa integral espacio dual} and the
H\"{o}lder-type inequality yield
\begin{align*}
\Vert B_{2}(u,v)\Vert_{\mathcal{M}[X^{\prime}]^{l^{\prime}}}\leq &
C\sup_{t>0}\Vert uS(v)\Vert_{\mathcal{M}[X_{0}^{\prime}]^{l_{0}^{\prime}}}\\
\leq &  \ C\sup_{t>0}\Vert u\Vert_{\mathcal{M}[X^{\prime}]^{l^{\prime}}}\Vert
S(v)\Vert_{\mathcal{M}\left[  X_{\kappa}^{\prime}\right]  ^{l_{1}}}.
\end{align*}
By taking $l_{1}=\frac{n}{\alpha-\delta}$, using the H\"{o}lder-type
inequality and the identity $\frac{n}{l_{0}^{\prime}}-\frac{n}{l^{\prime}}=1$,
it follows that $\alpha=\delta+1$. Finally, noting that $|S(v)(x)|\leq
(-\Delta)^{-\delta/2}|v|(x)$ and applying Proposition
\ref{prop:boundedness_fractional_operator}, we obtain
(\ref{eq:estim-bilin-advect-difuss_2}). \fin

\begin{remark}
Note that in Theorem \ref{Theorem:bilinear-estim-B2} we are forced to take
$n>\delta+2$ for the $X_{0}^{\prime}$-Morrey space to be well defined. Also,
Theorem \ref{Theorem:bilinear-estim-B2} works well under the condition stated
in Remark \ref{remark:Alternative_condition-fractional-operator} with the
obvious modifications.
\end{remark}

\subsection{Applications}

\subsubsection{Reaction-diffusion system}

As a case that realizes Theorem \ref{Theorem:bilinear-estim-B1}, we take
$X^{\prime}=L^{(p^{\prime},\infty)}$ and $X_{0}^{\prime}=L^{(p_{0}^{\prime
},\infty)}$ with $2p_{0}^{\prime}=p^{\prime}$. Also, we can choose
$p_{1},p_{2}>1$ such that
\[
\frac{p_{0}^{\prime}}{p_{1}^{\prime}}<\frac{1}{2}<\frac{p_{0}^{\prime}}%
{p_{2}^{\prime}}<1.
\]
The rest of parameters can be chosen as in Example
\ref{eq:bilinear-estimate-weak-morrey-spaces}. Then, we obtain
(\ref{eq:estim-bilin-react-difuss}) in the weak-Morrey space $\mathcal{M}%
_{\left(  p^{\prime},\infty\right)  }^{n/2}=$ $\mathcal{M}[L^{(p^{\prime
},\infty)}]^{n/2}$ for $n>4$ and $p^{\prime}\leq n/2$. As a consequence of
(\ref{eq:estim-bilin-react-difuss}), by adapting an argument by \cite{Meyer},
we can get uniqueness of mild solutions in the class $C(\left[  0,T\right)
;\tilde{Z})$ with $Z=\mathcal{M}_{\left(  p^{\prime},\infty\right)  }^{n/2}$
for system (\ref{eq:reac-difus-syst}).

\subsubsection{Nonlocal advection-diffusion system}

We conclude this work with an application for the nonlocal parabolic equation.
From Theorem \ref{Theorem:bilinear-estim-B2}, we obtain the bilinear estimate
for $B_{2}$ in weak-Morrey spaces by taking $X^{\prime}=L^{(p^{\prime}%
,\infty)}$, $X_{0}^{\prime}=L^{(p_{0}^{\prime},\infty)}$, $\alpha=\delta+1$
and $p_{0}^{\prime}=\left(  \frac{\delta+1}{\delta+2}\right)  p^{\prime}$.
Recall that the boundedness of the maximal operator $M_{0}$ is well
established in the context of weak-$L^{r}$ spaces, for all $r>1$. The
necessary H\"{o}lder-type inequality reduces to
\[
\Vert f\cdot g\Vert_{L^{(p_{0}^{\prime},\infty)}}\leq C\Vert f\Vert
_{L^{(p^{\prime},\infty)}}\Vert g\Vert_{L^{(p^{\prime}(\delta+1),\infty)}}.
\]
Therefore, the bilinear estimate for $B_{2}$ is verified in the weak-Morrey
space $\mathcal{M}_{\left(  p^{\prime},\infty\right)  }^{n/(1+\delta
)}=\mathcal{M}\left[  L^{(p^{\prime},\infty)}\right]  ^{n/(1+\delta)}$ with
$\frac{\delta+2}{\delta+1}<p^{\prime}\leq\frac{n}{\delta+1}$. As an
application of that estimate, we also obtain the uniqueness of mild solutions
in the class $C(\left[  0,T\right)  ;\tilde{Z})$ with $Z=\mathcal{M}_{\left(
p^{\prime},\infty\right)  }^{n/(1+\delta)}$ for equation
(\ref{eq:non-local-model}).

\end{document}